\def\tex#1#2{#2} 

\tex
{}
{
\documentclass[10pt]{article}
\usepackage[a4paper, top=2.7cm, left=2.5cm, right=2.5cm, bottom=2.6cm]
{geometry}
\usepackage{amscd}
\usepackage{amssymb}
\usepackage{amsmath}
\usepackage{amsthm}
\usepackage{graphicx}
\usepackage{mathrsfs}
\usepackage{hyperref}
\def\cal{\mathcal}
\begin{document}
}

\tex
{
\def\[{[}
\def\]{]}
\def\cite{{}}   
}
{
\def\[{}
\def\]{}
}
%


\tex 
{
\input amssym.tex
\input graphicx.tex
\input color.tex
\input miniltx.tex
}
{}
\hoffset=.4truecm
\voffset=-0.5truecm
\hsize=15truecm
\vsize=24truecm
\baselineskip 16pt
\parindent=0cm

\def\litem{\par\noindent\hangindent=\parindent\ltextindent}
\def\ltextindent#1{\hbox to \hangindent{#1\hss}\ignorespaces}
\long\def\ignore#1\recognize{}
\def\no#1{}
\def\comment#1{}
\def\label#1{}

\def\big{\bigskip}
\def\med{\medskip}
\def\hs{\hskip}
\def\vs{\vskip}
\def\qed{\hfill$\circlearrowleft$\med}
\def\gb{\goodbreak}
\def\cl{\centerline}
\def\ds{\displaystyle}
\def\ol{\overline}
\def\ul{\underline}
\def\sm{\setminus}
\def\ld{\ldots}
\def\cd{\cdot}
\def\to{,\ldots,}
\def\sub{\subseteq}
\def\{{\lbrace}
\def\}{\rbrace}
\def\congruent{\equiv}
\def\isom{\cong}
\def\map{\rightarrow}
\def\inv{^{-1}}
\def\abs#1{\vert#1\vert}
\def\wt{\widetilde}
\def\wh{\widehat}

\def\N{{\Bbb N}}
\def\R{{\Bbb R}}
\def\C{{\Bbb C}}
\def\Z{{\Bbb Z}}
\def\Q{{\Bbb Q}}
\def\P{{\Bbb P}}
\def\K{{\bf K}}
\def\G{{\bf G}}
\def\H{{\bf H}}

\def\6{\partial}
\def\min{{\rm min}}
\def\max{{\rm max}}
\def\ord{{\rm ord}}
\def\Id{{\rm Id}}
\def\supp{{\rm supp}}
\def\Ker{{\rm Ker}}
\def\Im{{\rm Im}}
\def\Re{{\rm Re}}
\def\Grass{{\rm Grass}}
\def\a{\alpha}
\def\b{\beta}
\def\ai{{\alpha_i}}
\def\aj{{\alpha_j}}
\def\bi{{\beta_i}}
\def\bj{{\beta_j}}
\def\t{\theta}
\def\del{\delta}
\def\eps{\varepsilon}
\def\O{{\cal O}}
\def\RR{{\cal R}}
\def\A{{\Bbb A}}
\def\T{{\Bbb T}}
\def\Aut{{\rm Aut}}
\def\Autc{{\rm Aut}_\C}
\def\Autr{{\rm Aut}_\R}
\def\GL{{\rm GL}}
\def\Spec{{\rm Spec}}
\def\semidirect{\ \vert\!\!\!\times\  }
\def\1{{\bf 1}}
\def\inin{{\rm in}}
\def\co{{\rm co}}
\def\FF{{\cal F}}
\def\NN{{\cal N}}
\def\PP{{\cal P}}
\def\QQ{{\cal Q}}
\def\DD{{\cal D}}
\def\TT{{\cal T}}
\def\UU{{\cal U}}
\def\VV{{\cal V}}
\def\XX{{\cal X}}
\def\YY{{\cal Y}}
\def\ZZ{{\cal Z}}
\def\GG{{\cal G}}
\def\HH{{\cal H}}
\def\MM{{\cal M}}
\def\AA{{\cal A}}
\def\BB{{\cal B}}
\def\CC{{\cal C}}
\def\SS{{\cal S}}
\def\II{{\cal I}}
\def\IP{{\cal IP}}
\def\DP{{\cal DP}}
\def\PGL{{\rm PGL}}
\def\cross{{c\!\!c}}
\def\pp{{ \rm p\!\!\!p }}
\def\qq{{ \rm q\!\!\!q }}
\def\xx{{ \rm x\!\!x }}
\def\yy{{ \rm y\!\!y}}
\def\zz{{ \rm z\!\!z}}
\def\kk{K}
\def\CR{{\rm CR}}

\def\sk{{\rm sk}}
\def\representative{representative }
\def\projection{projection }
\def\triples{{\rm triples}}
\def\qua{{\rm quad}}
\def\sym{{\rm string}}
\def\proj{{\rm proj}}
\def\jj{{\it J$\,{}^2$}}


\def\ooo{\red}


\font\Times=ptmr at 10pt
\font\smallTimes=ptmr at 8pt
\font\smallit=ptmri at 8pt
\font\bf=ptmb at 10pt
\font\Bf=ptmb at 12pt
\font\BF=ptmb at 14pt
\font\Bfit=ptmbi at 12pt
\font\BFit=ptmbi at 14pt

\Times



\def\red{\includegraphics[width=.02\hsize]{reddots.pdf} }


\no{1}\def\darwinsdrawing{{\includegraphics[width=.5\hsize]{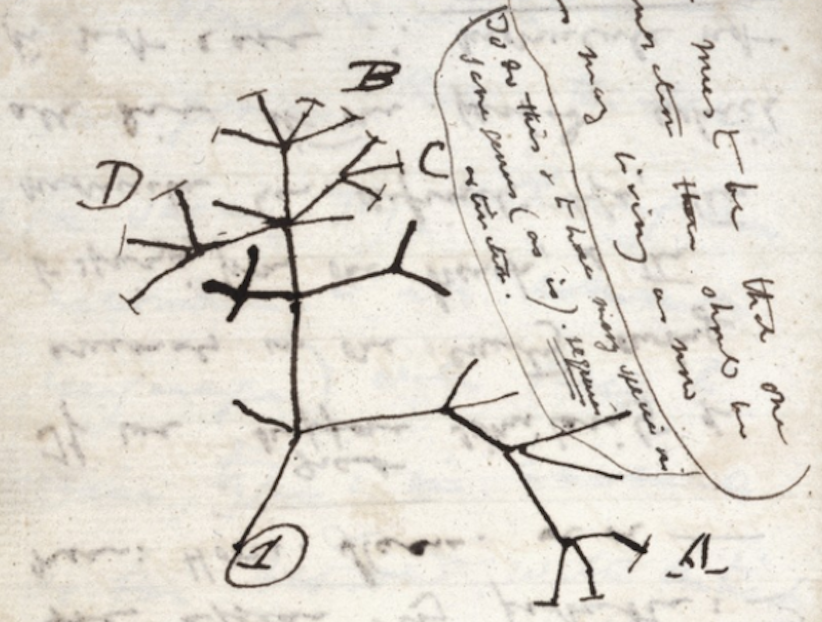}}}

\no{10}\def\fourlabelstree{{\includegraphics[width=.55\hsize]{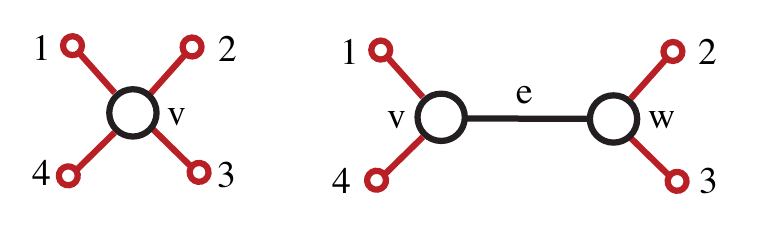}}}

\no{5}\def\stablecurve{{\includegraphics[width=.4\hsize]{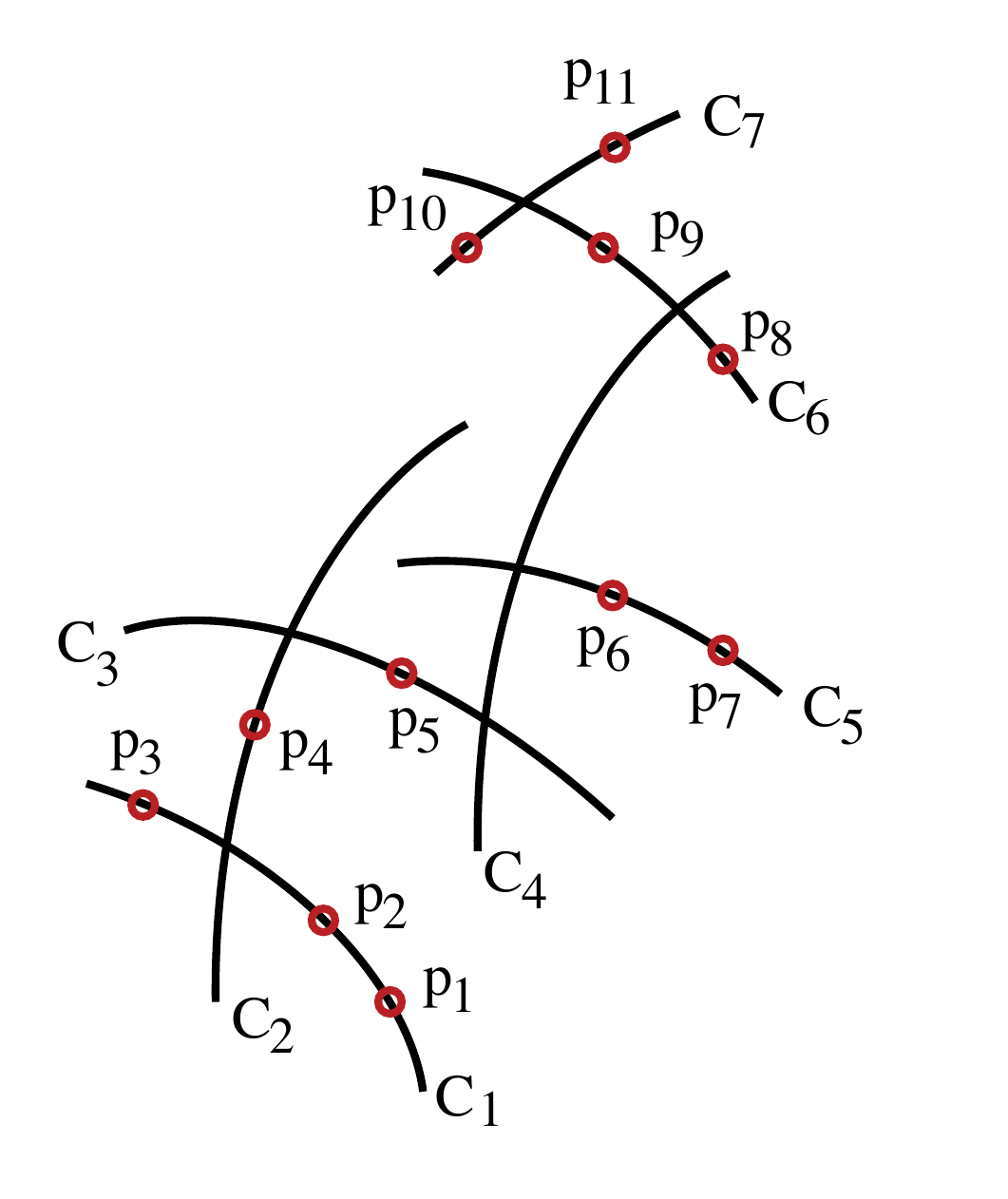}}}

\no{6}\def\phylogenetictreeexample{{\includegraphics[width=.5\hsize]{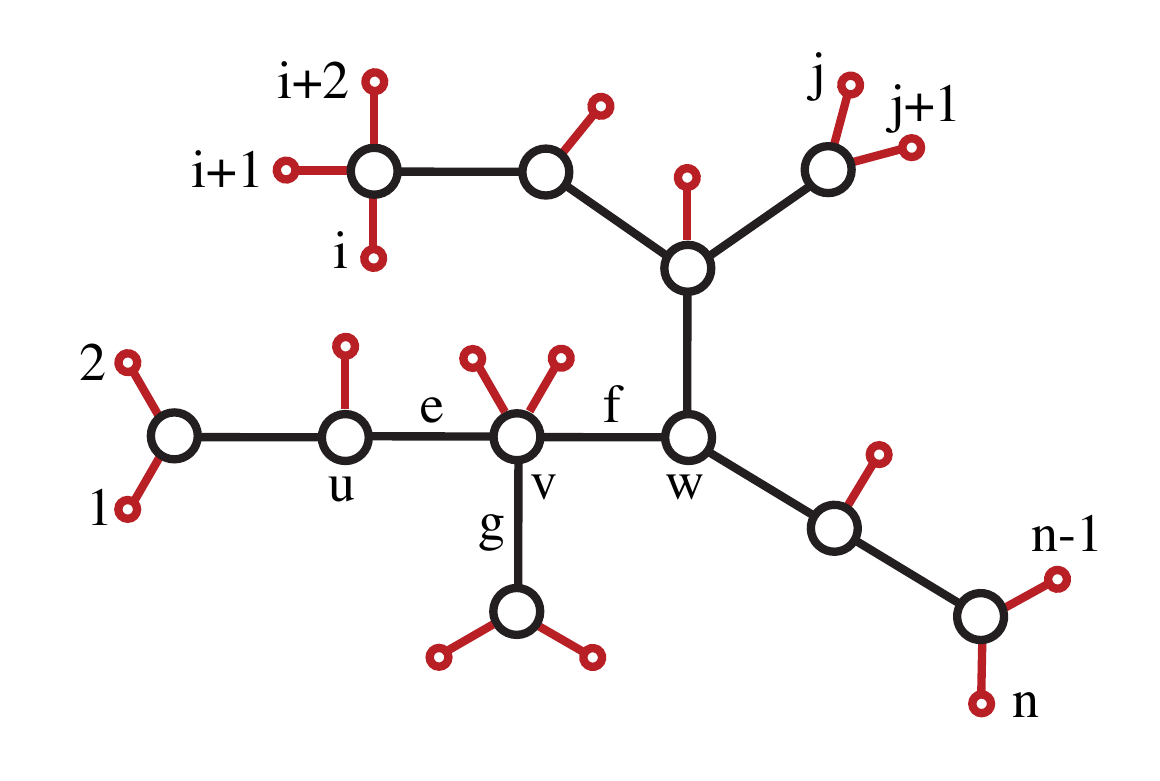}}}

\no{10}\def\fivelabelstreegeneric{{\includegraphics[width=.23\hsize]{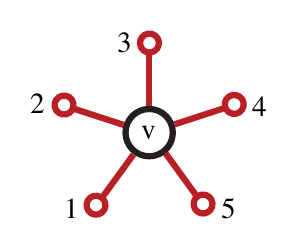}}}

\no{100 + 13}\def\fivelabelstreesspecial{{\includegraphics[width=.75\hsize]{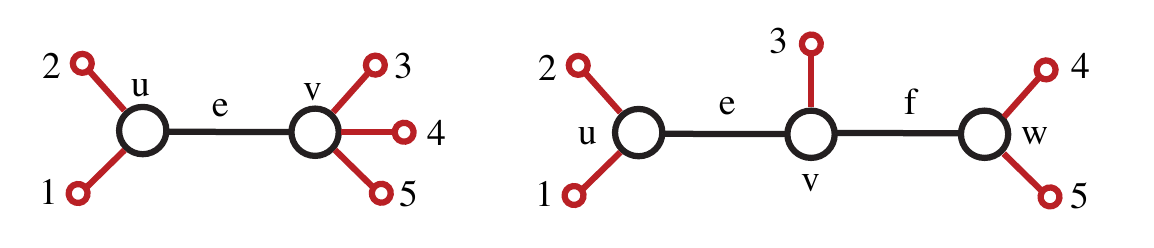}}}

\no{12}\def\phylogenetictrees{{\includegraphics[width=.9\hsize]{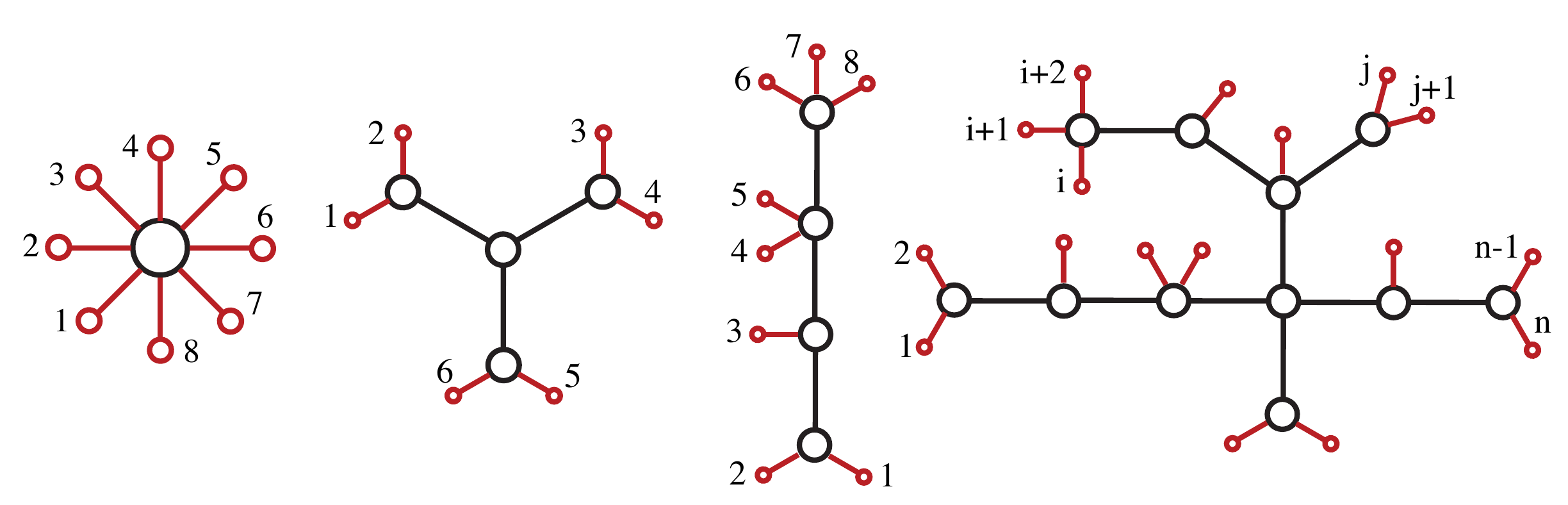}}}

\no{120}\def\meetingpoints{{\includegraphics[width=1\hsize]{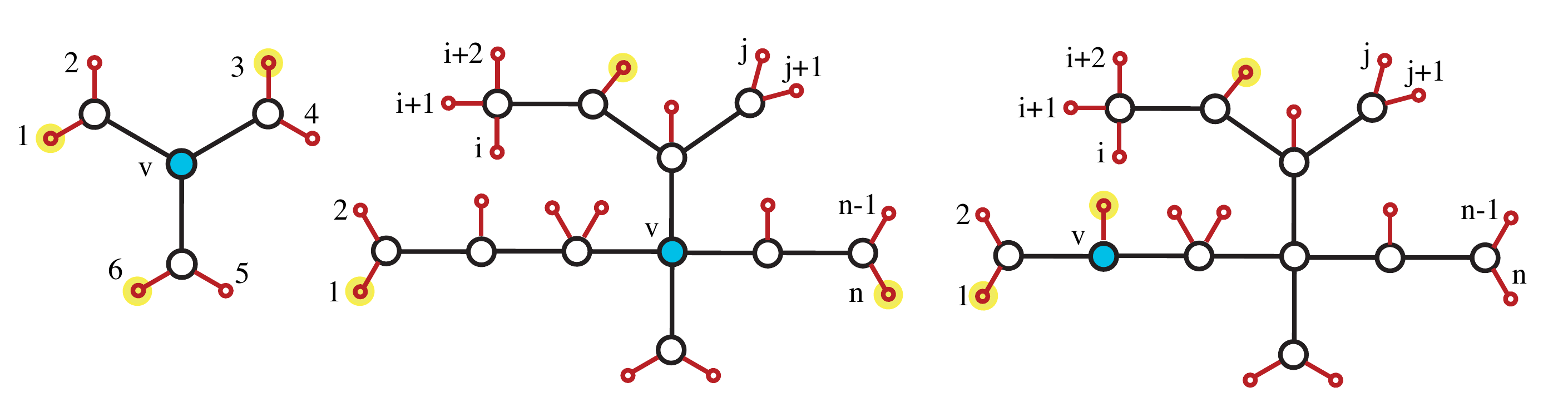}}}

\no{121}\def\destinationsets{{\includegraphics[width=.95\hsize]{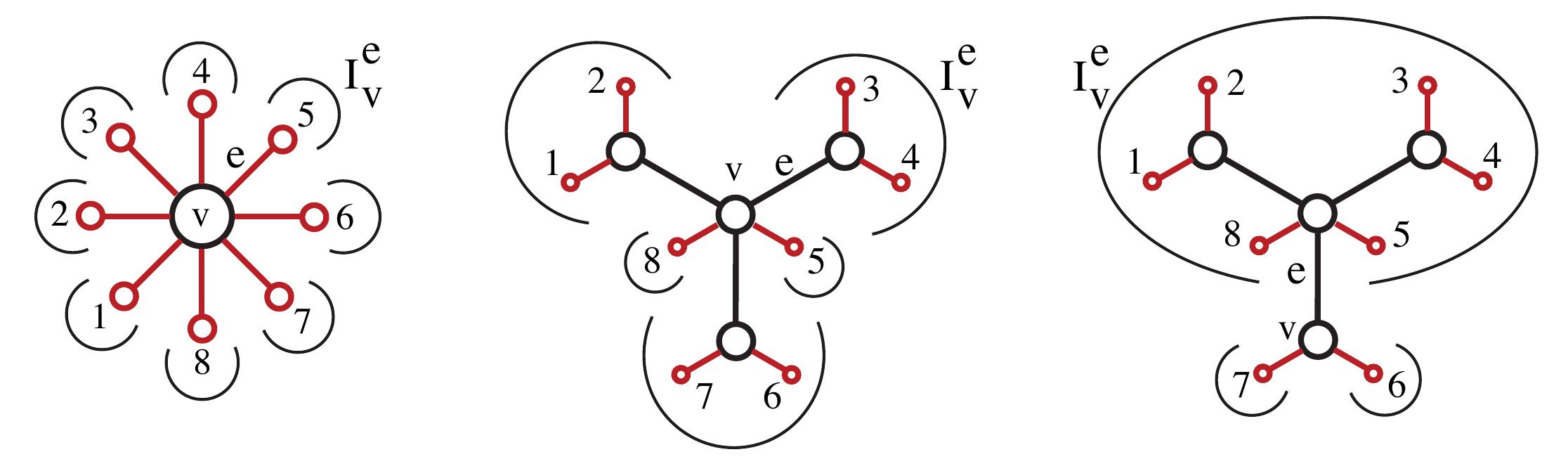}}}

\no{122}\def\destinationsetscomplementary{{\includegraphics[width=.45\hsize]{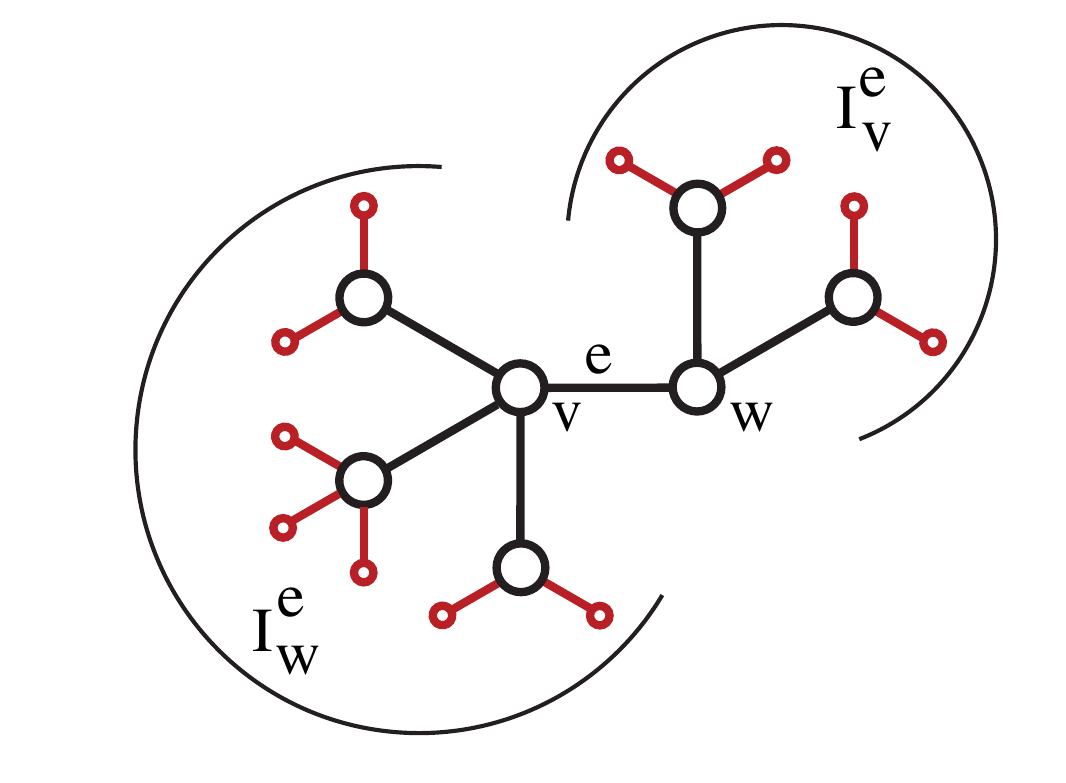}}}

\no{123}\def\clippingleavesA{{\includegraphics[width=.95\hsize]{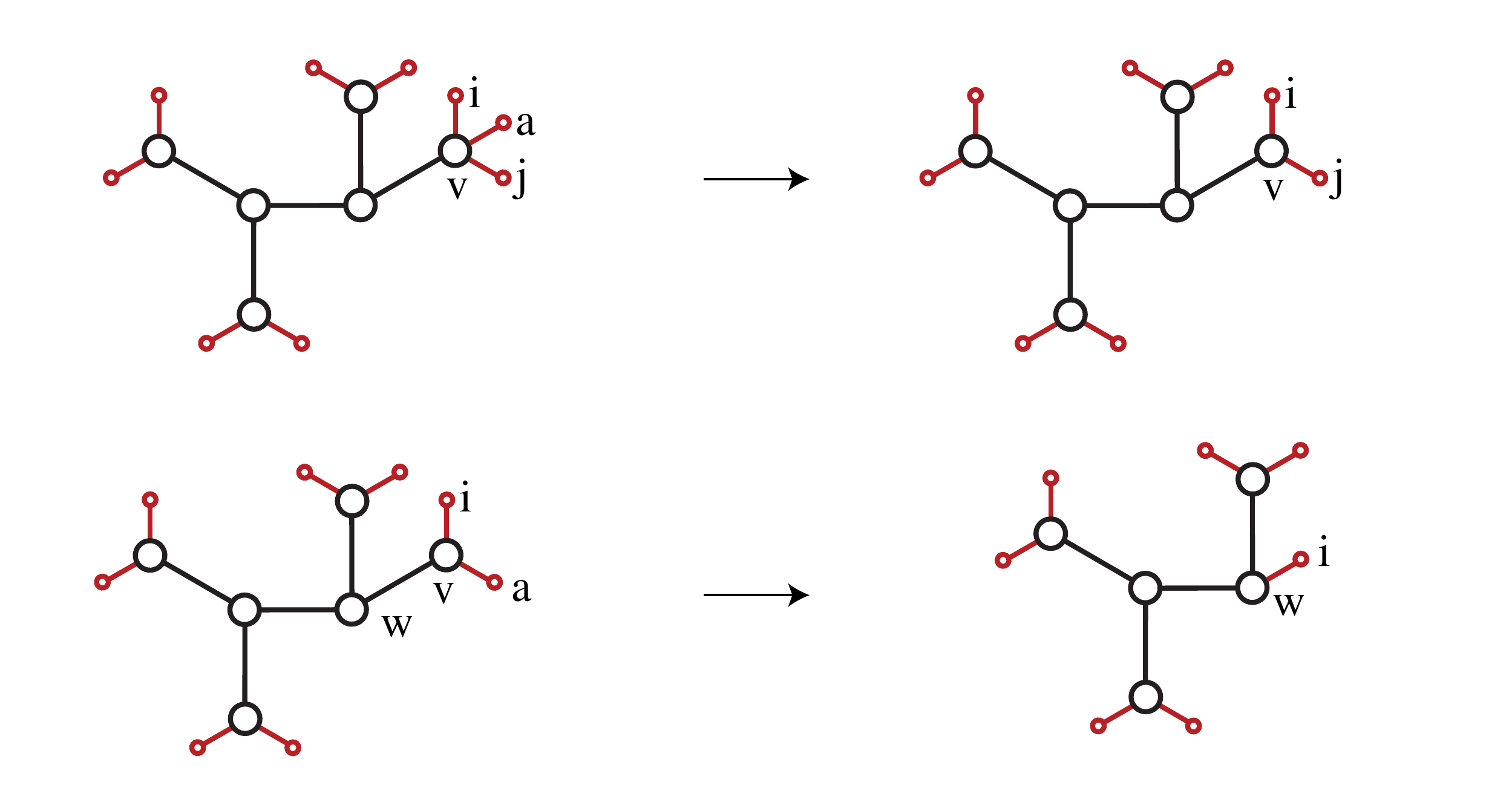}}}

\no{123}\def\clippingleavesB{{\includegraphics[width=.95\hsize]{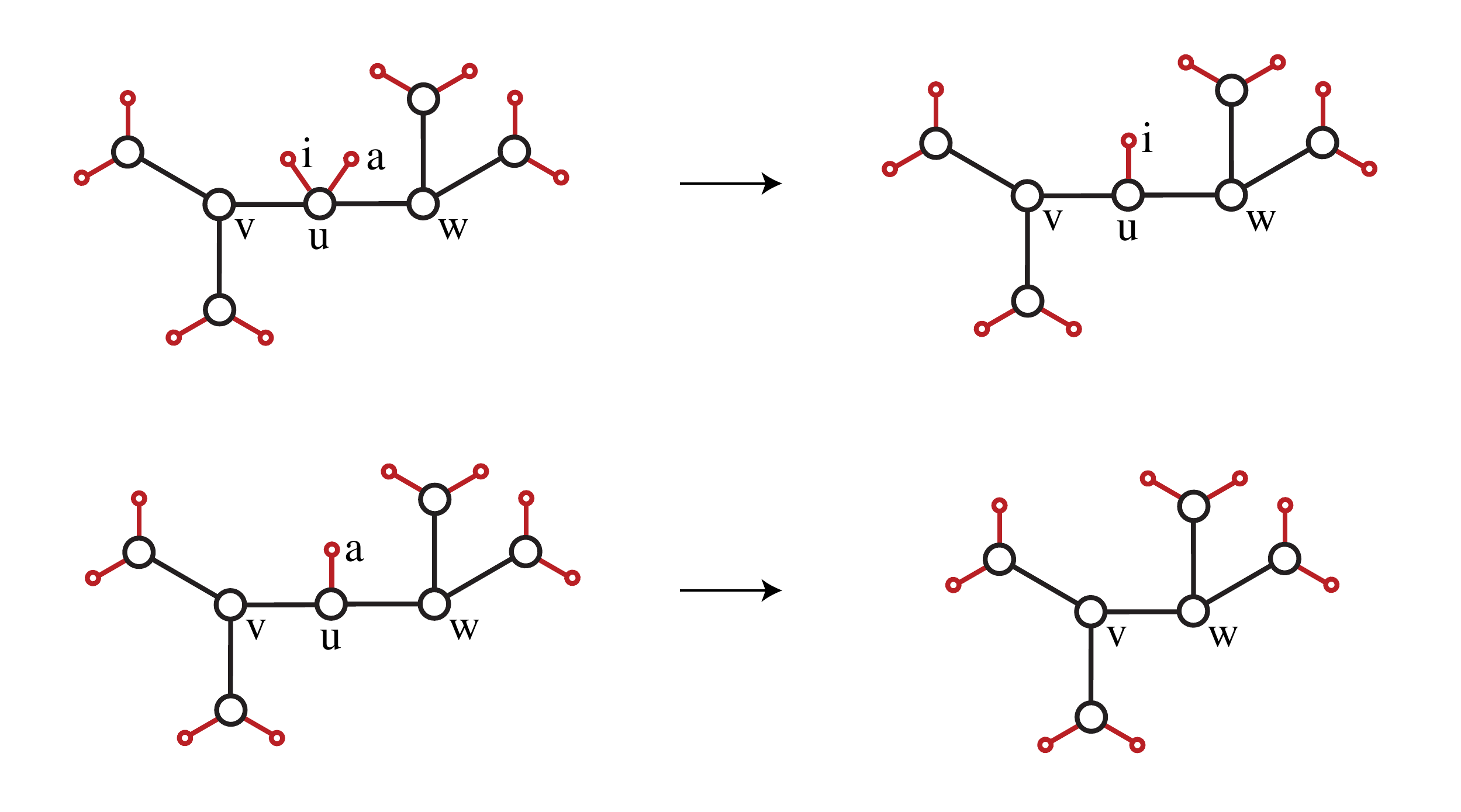}}}

\no{124}\def\deletingedgeA{{\includegraphics[width=.8\hsize]{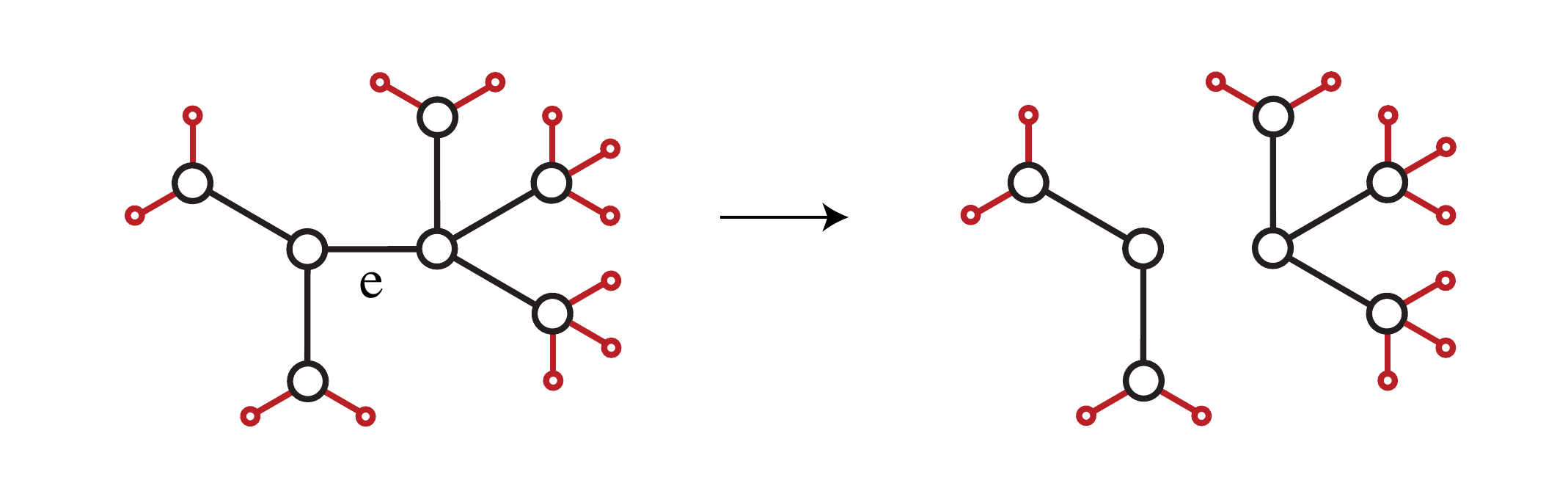}}}

\no{124}\def\deletingedgeB{{\includegraphics[width=.8\hsize]{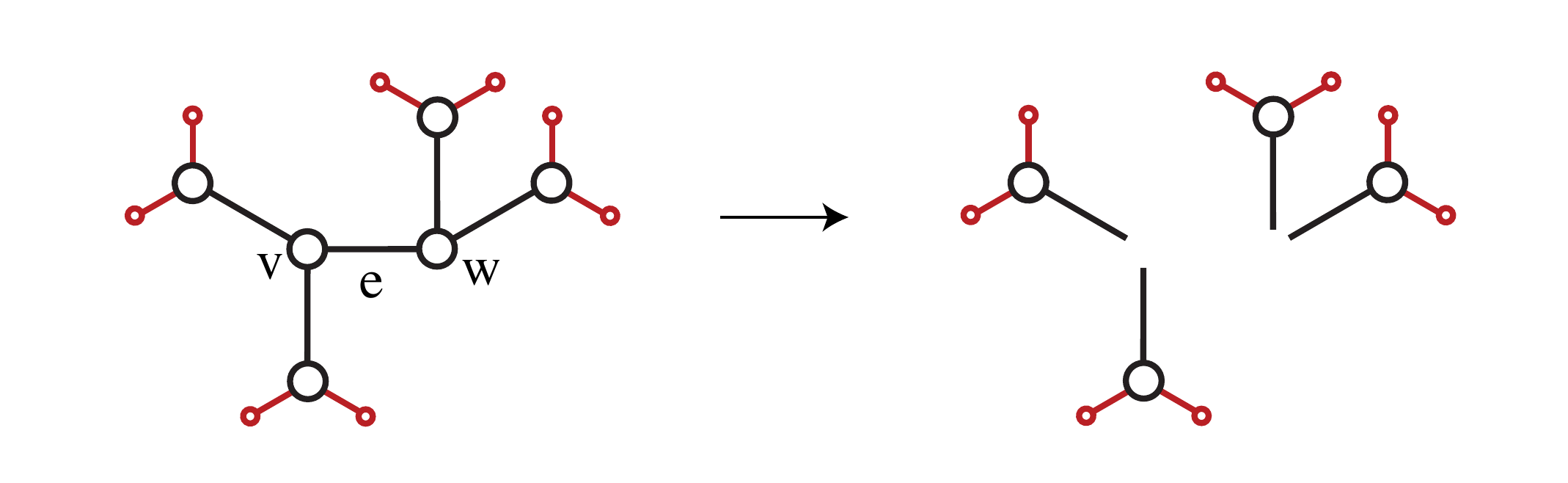}}}

\no{125}\def\contractinginsertingedge{{\includegraphics[width=.8\hsize]{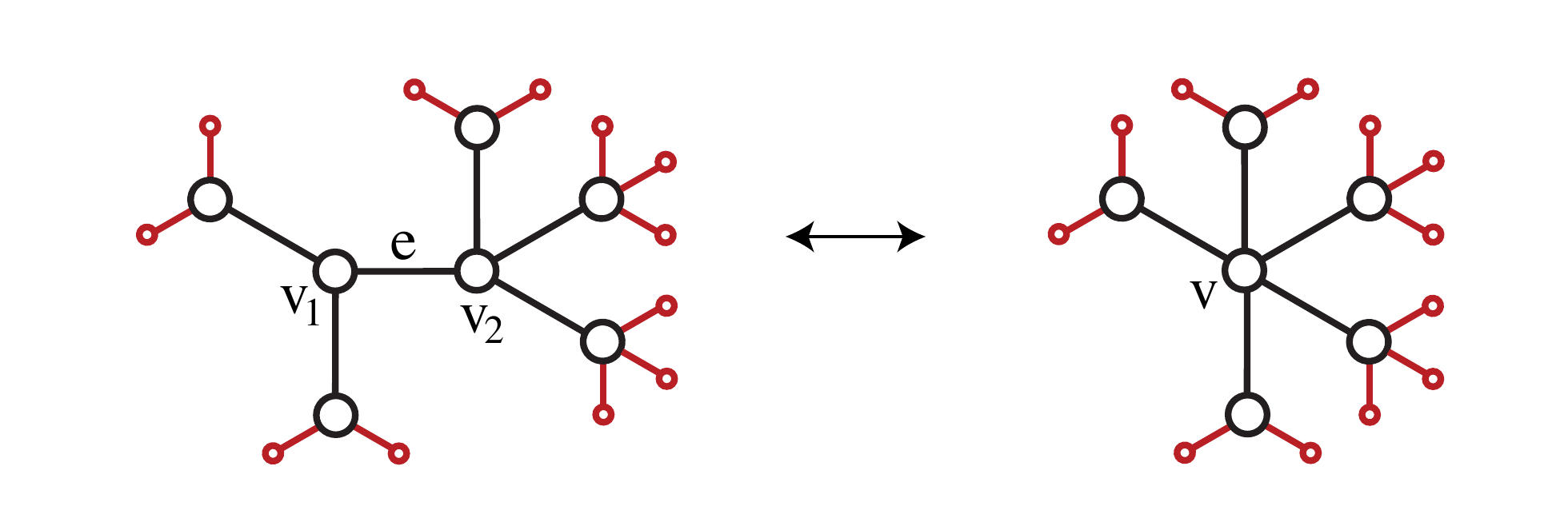}}}

\no{13 + 100}\def\fivelabelstreesspecial{{\includegraphics[width=.75\hsize]{fivelabelstreesspecial.pdf}}}

\no{130}\def\lemmacomplementaryincidence{{\includegraphics[width=.7\hsize]{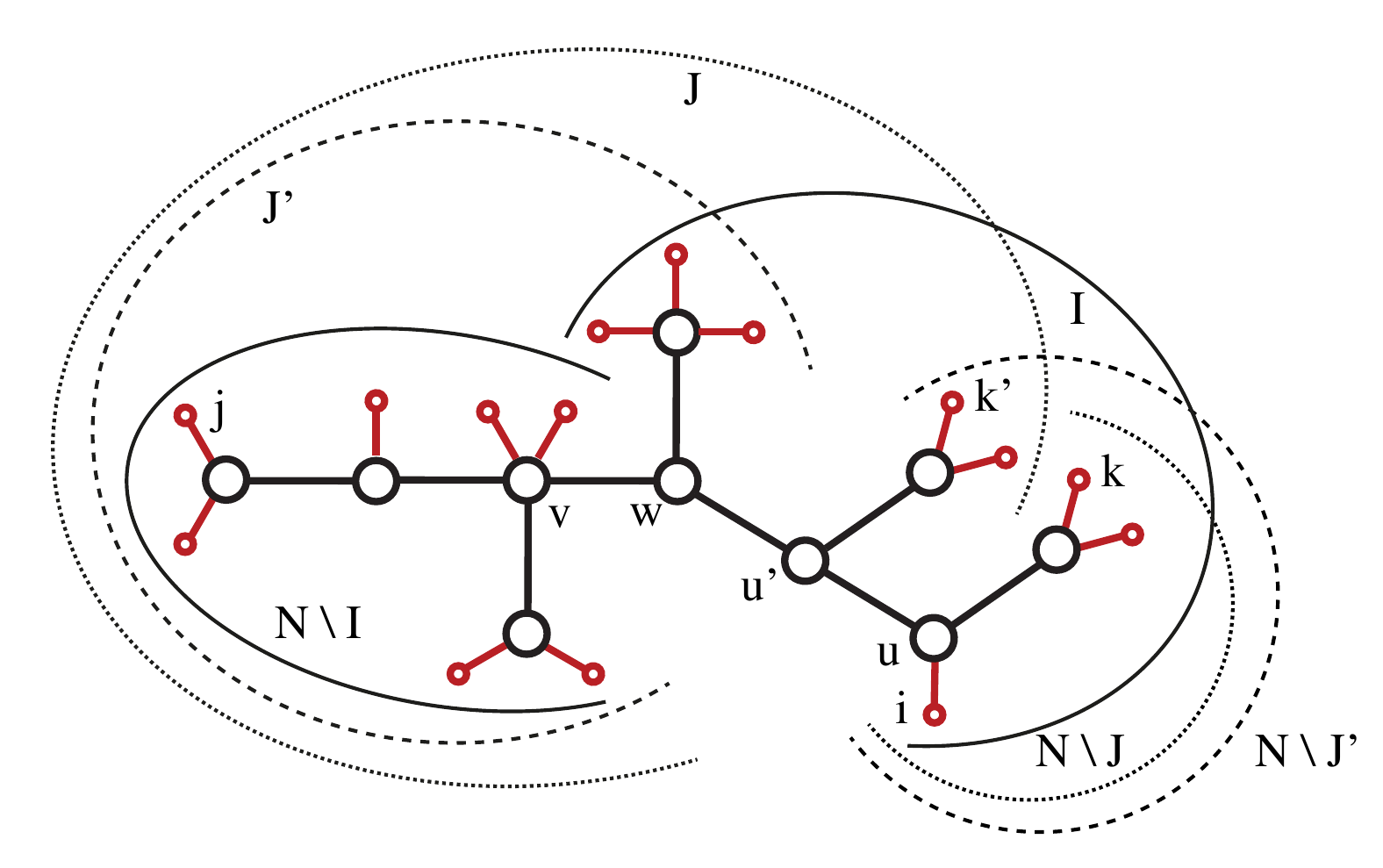}}}

\no{131}\def\clippingleaf{{\includegraphics[width=.8\hsize]{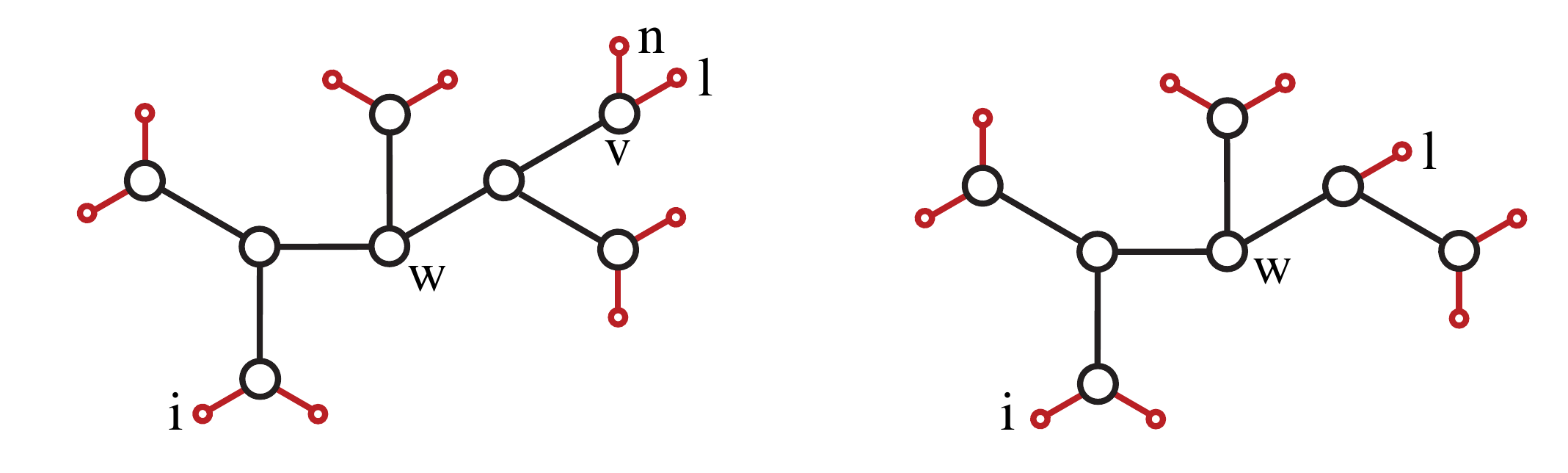}}}

\no{14} \def\stablecurvesandtree{{\includegraphics[width=.8\hsize]{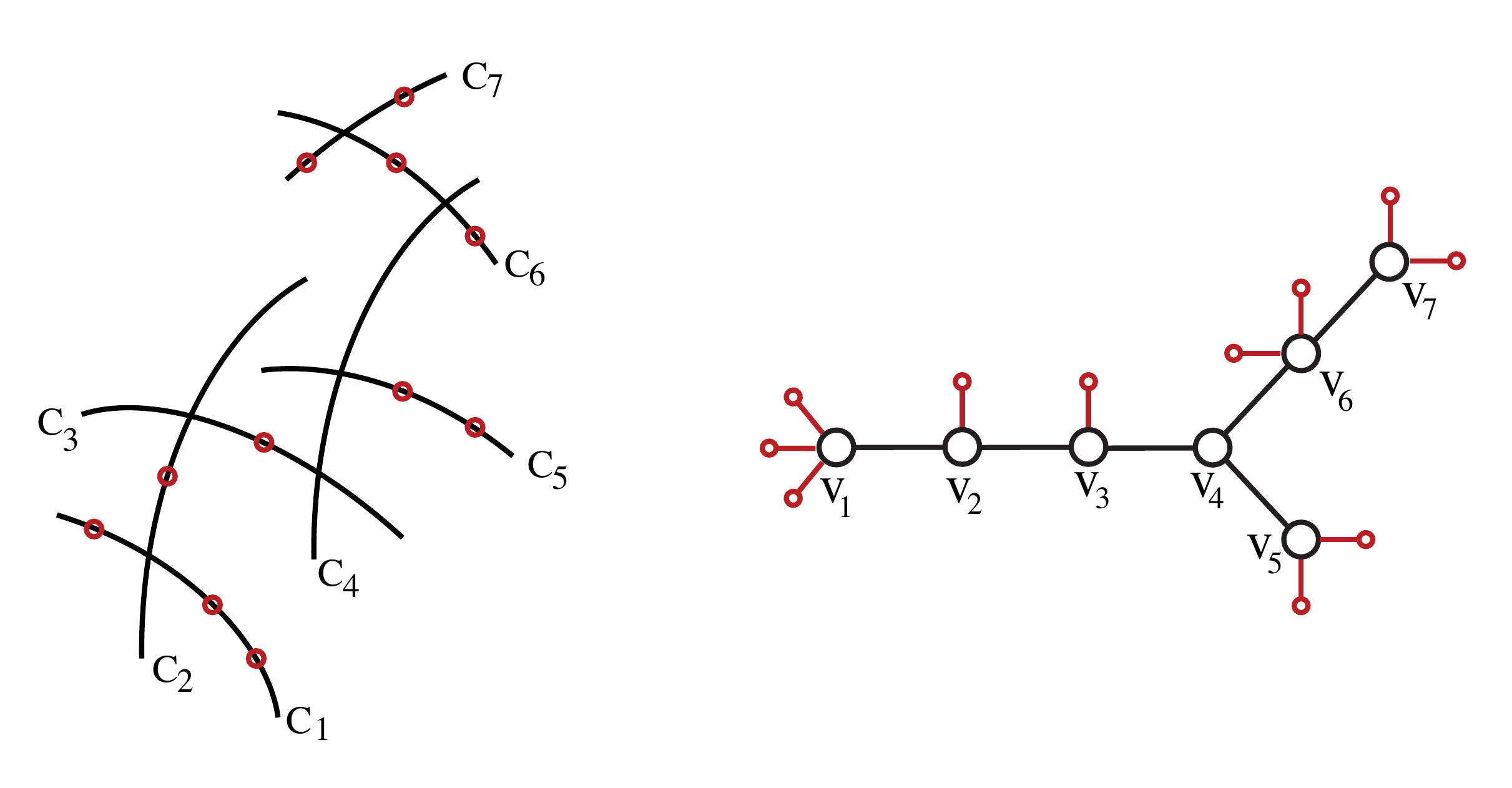}}}

\no{140} \def\medianprojection{{\includegraphics[width=.5\hsize]{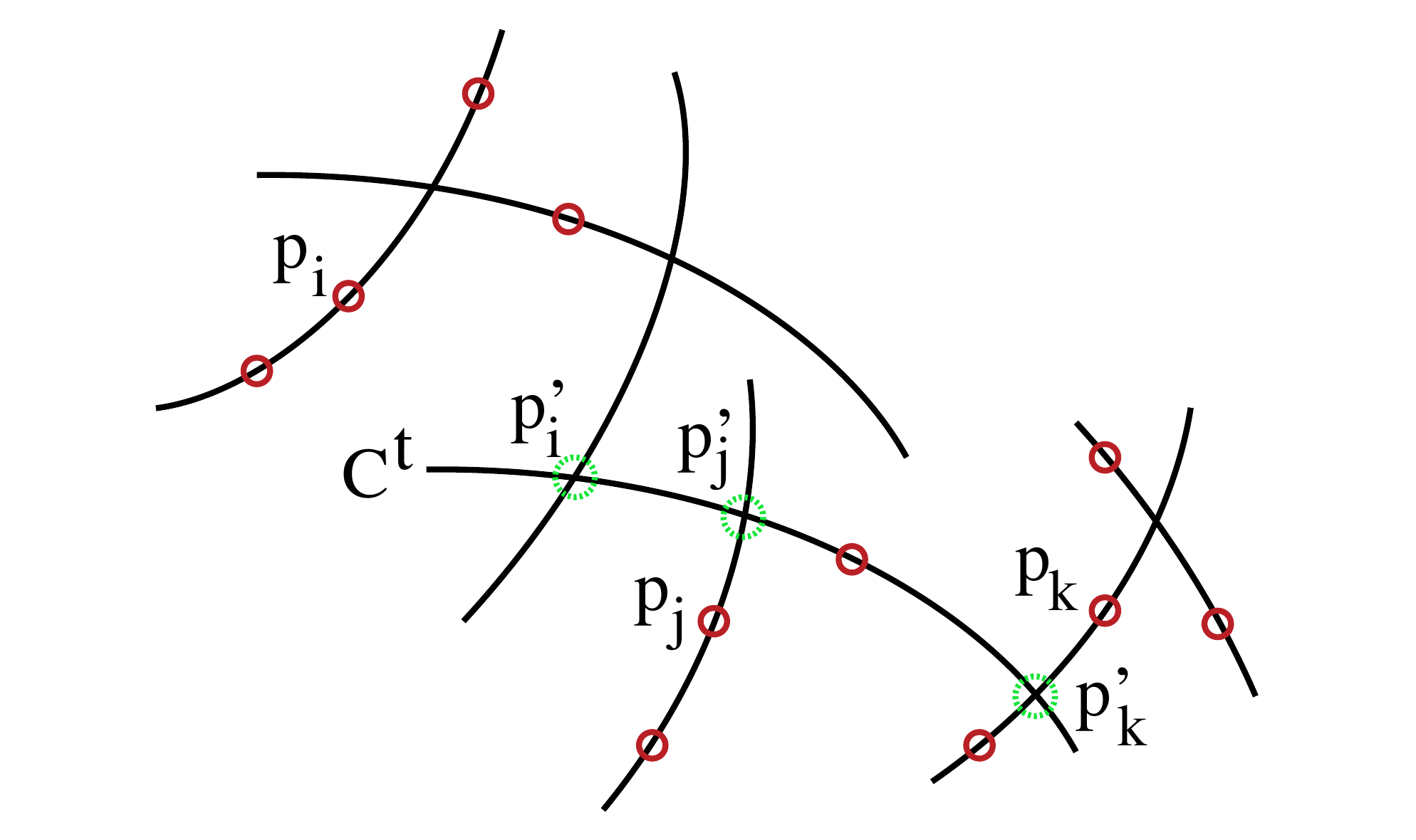}}}

\no{17}\def\vertexsplitting{{\includegraphics[width=.7\hsize]{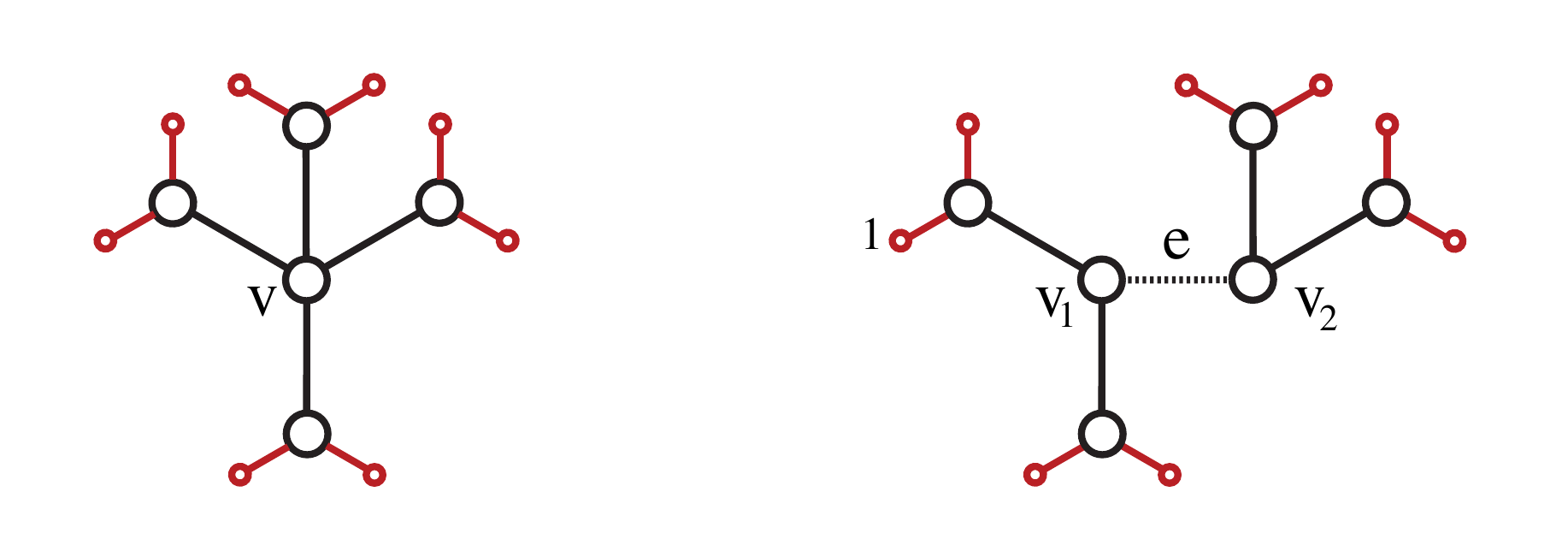}}}

\no{170}\def\fivelabelstreeextremal{{\includegraphics[width=.45\hsize]{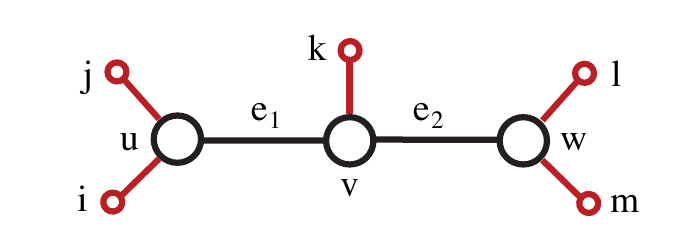}}}

\no{171}\def\Htree{{\includegraphics[width=.7\hsize]{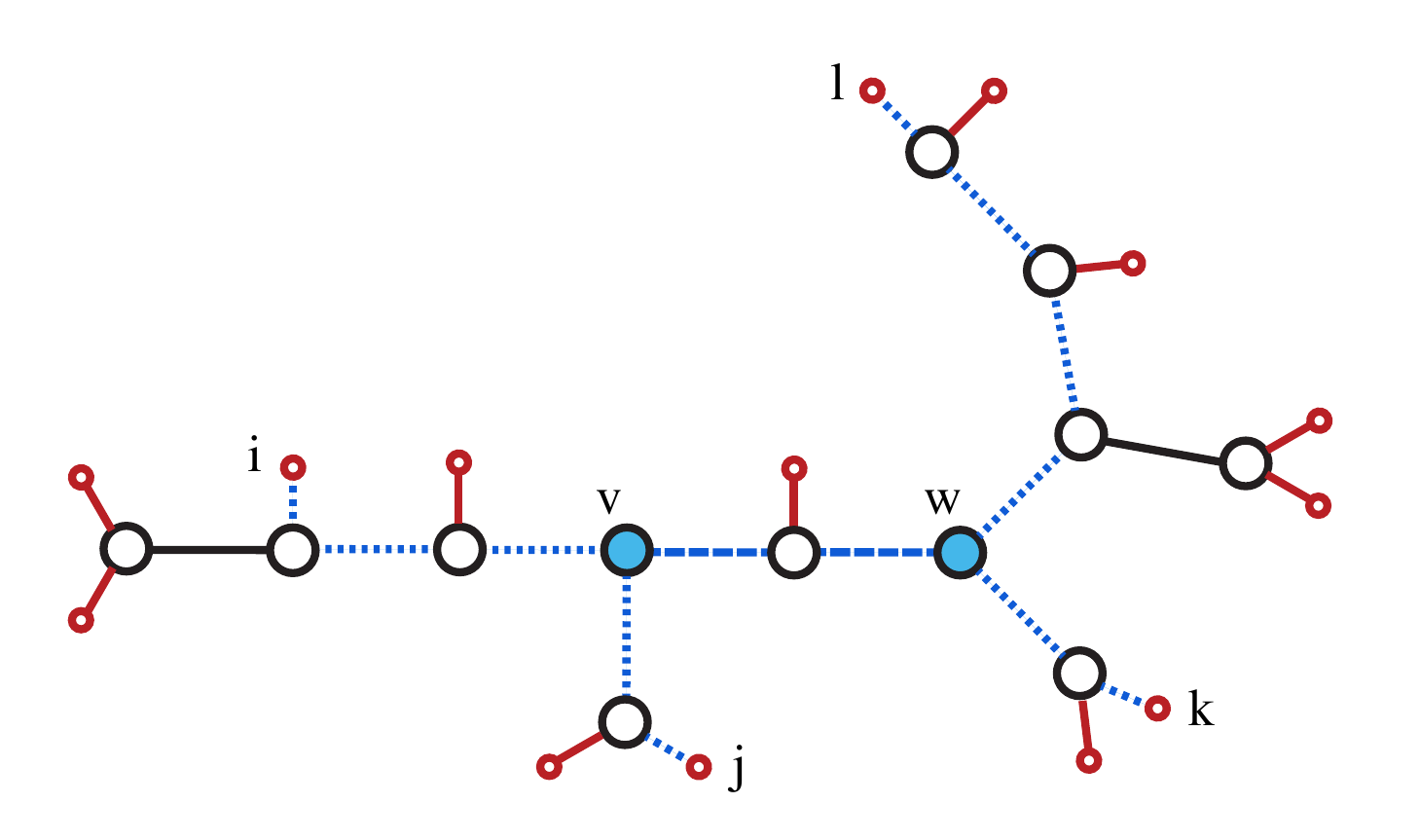}}}

\no{172}\def\Htreeinduction{{\includegraphics[width=.7\hsize]{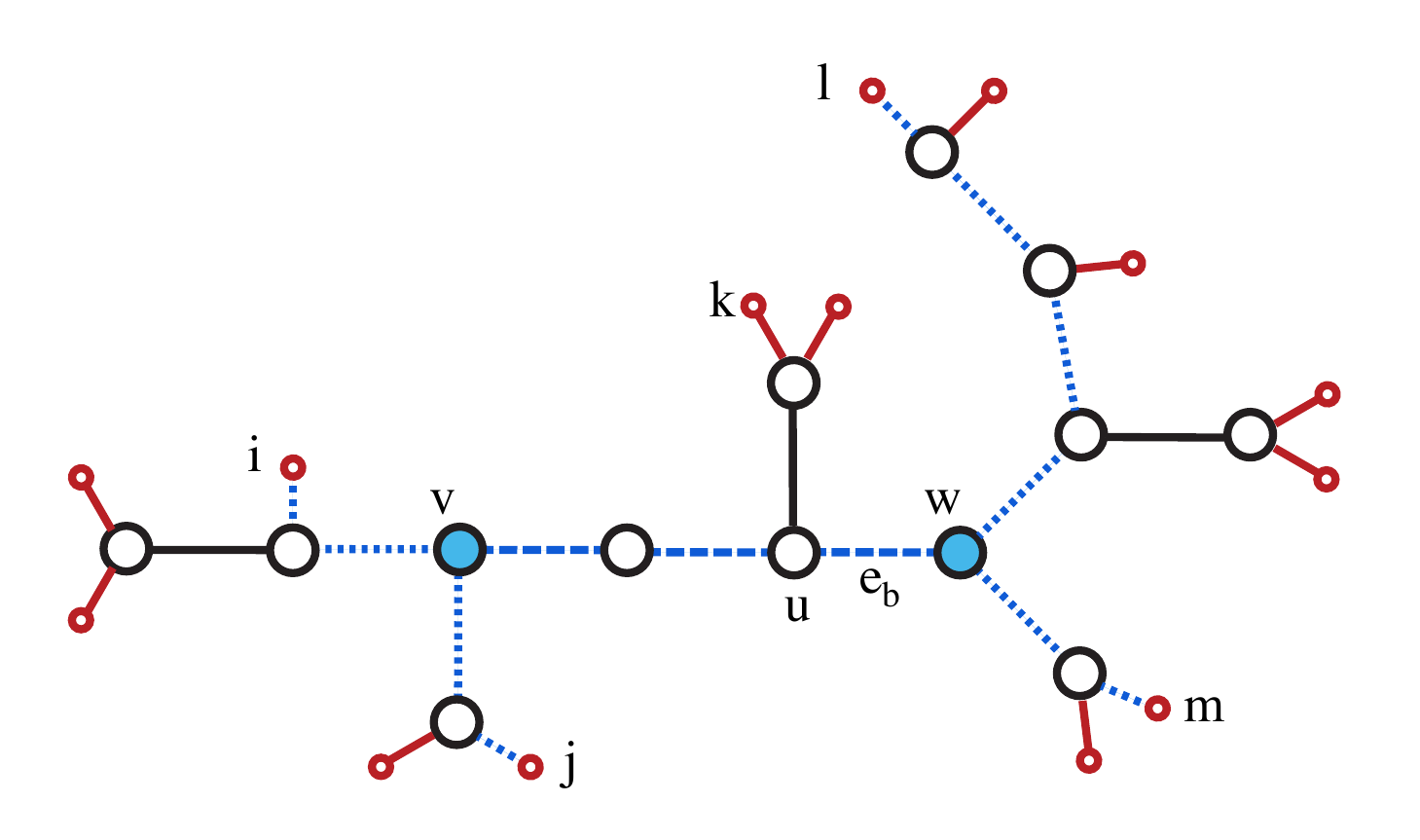}}}

\no{173} \def\clippingleafinduction{{\includegraphics[width=.75\hsize]{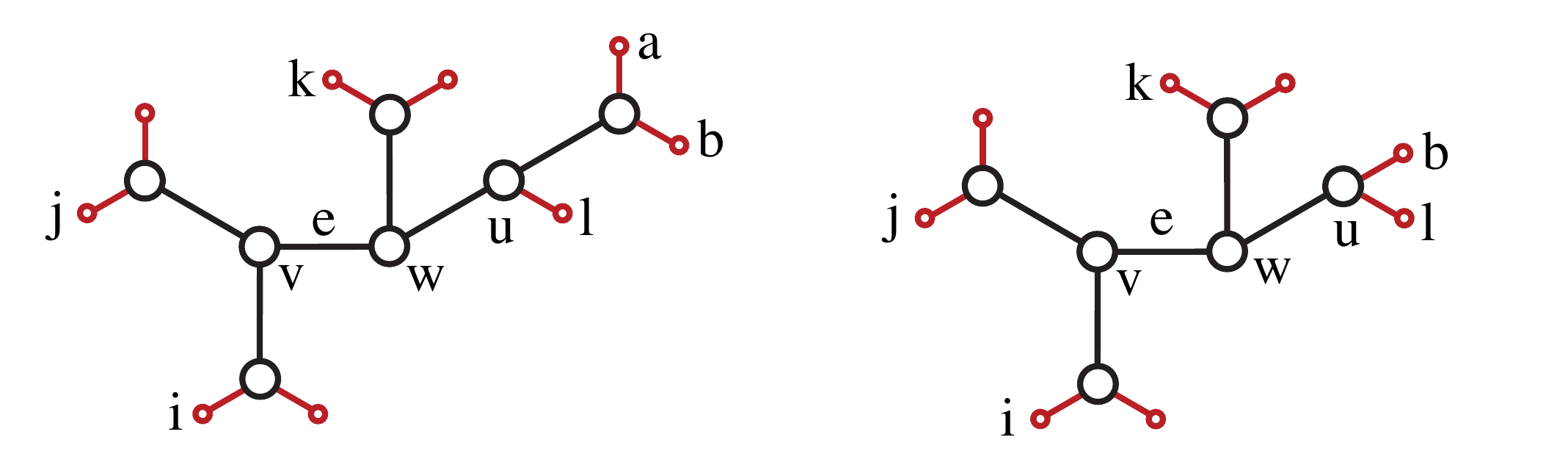}}}

\no{174}\def\clippingleafinductionii{{\includegraphics[width=.75\hsize]{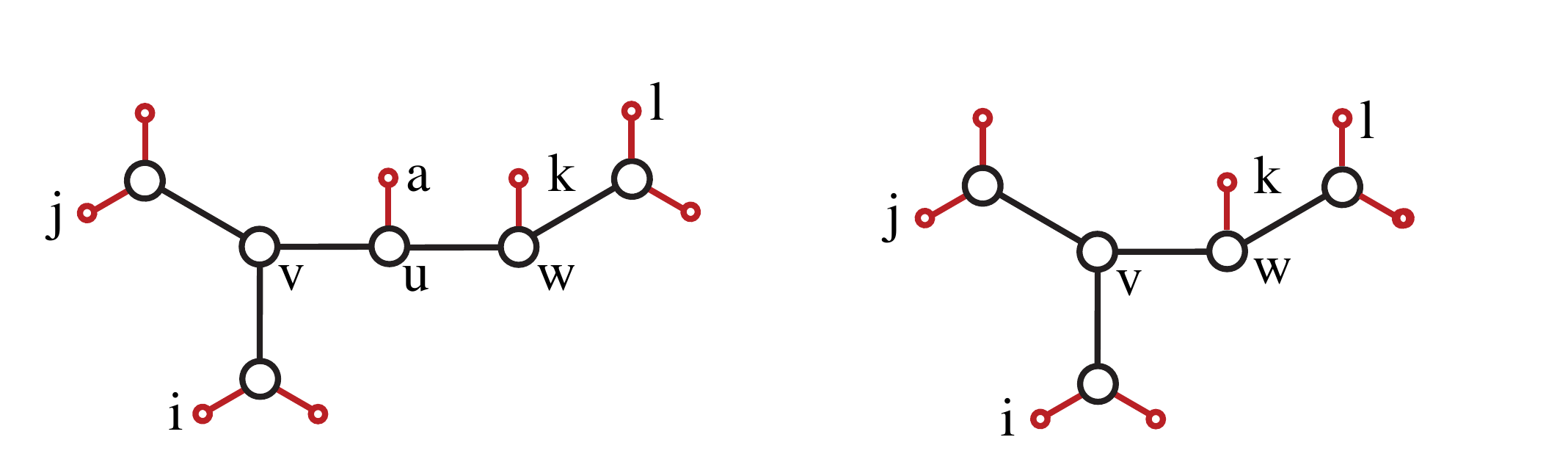}}}

\no{19} \def\twoverticestree{{\includegraphics[width=.5\hsize]{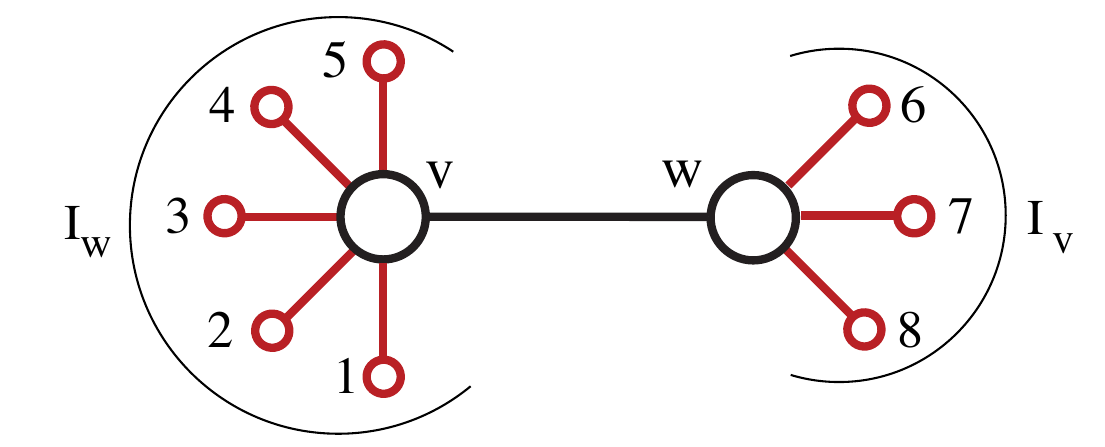}}}

\no{20} \def\destinationsimplepaths{{\includegraphics[width=.4\hsize]{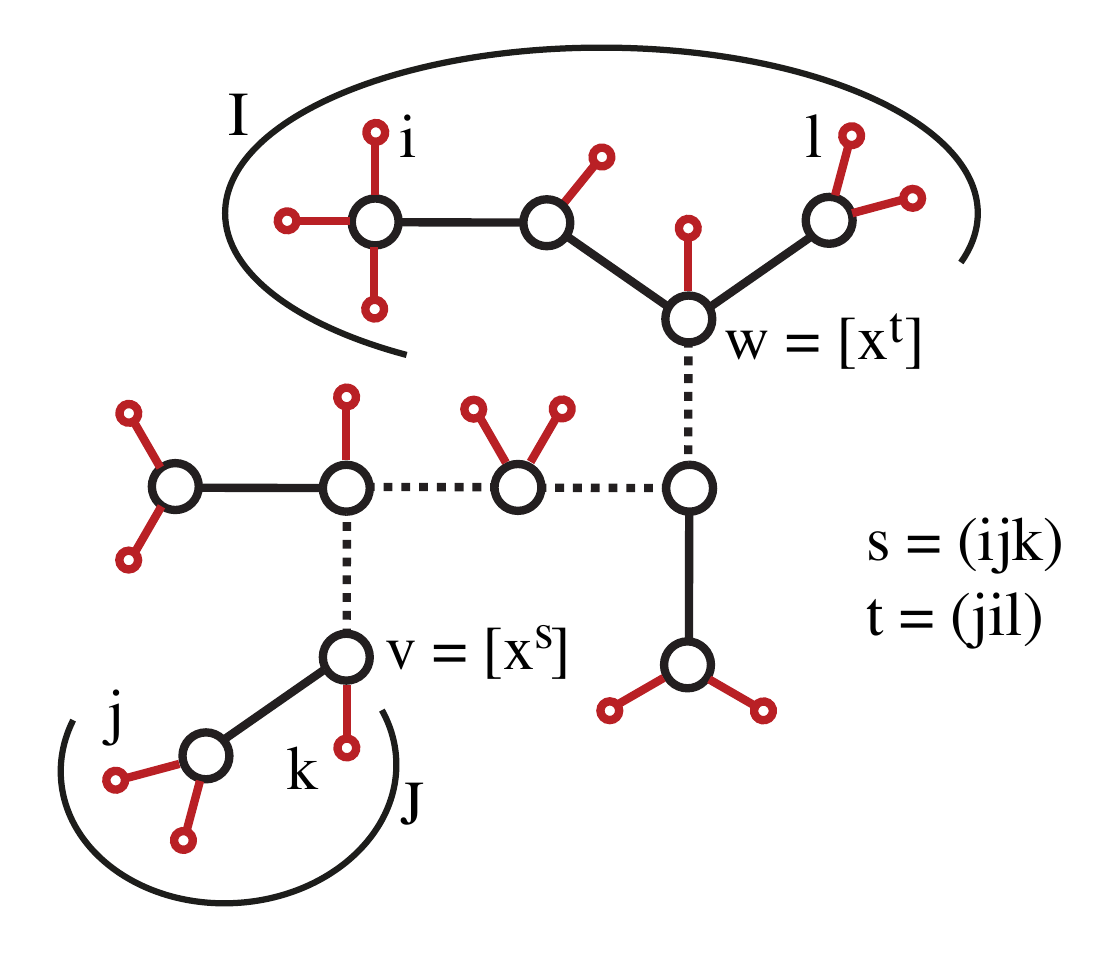}}}

\no{200} \def\ribisltree{{\includegraphics[width=.33\hsize]{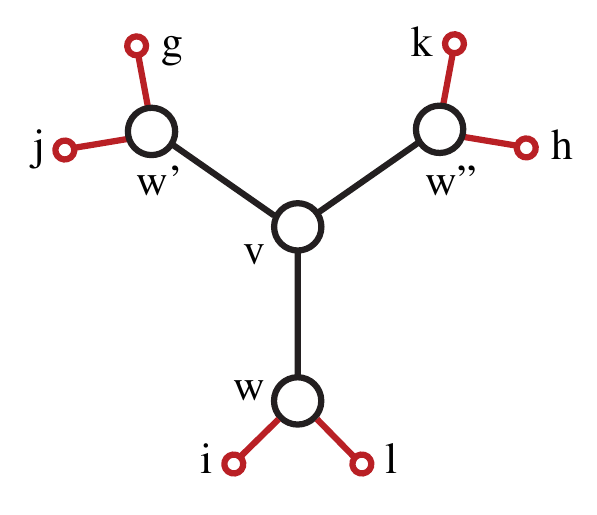}}}

\no{201} \def\viertelvorachttreelemmaseven{{\includegraphics[width=.33\hsize]{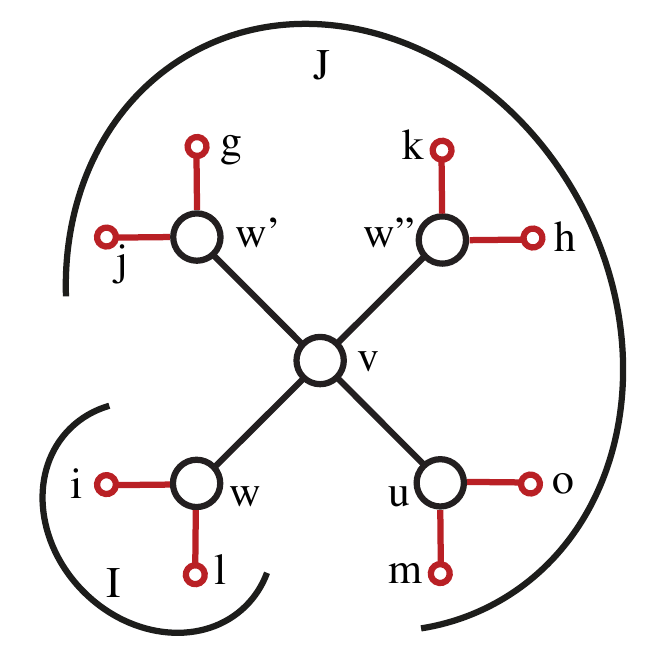}}}

\no{202} \def\choicetriples{{\includegraphics[width=.55\hsize]{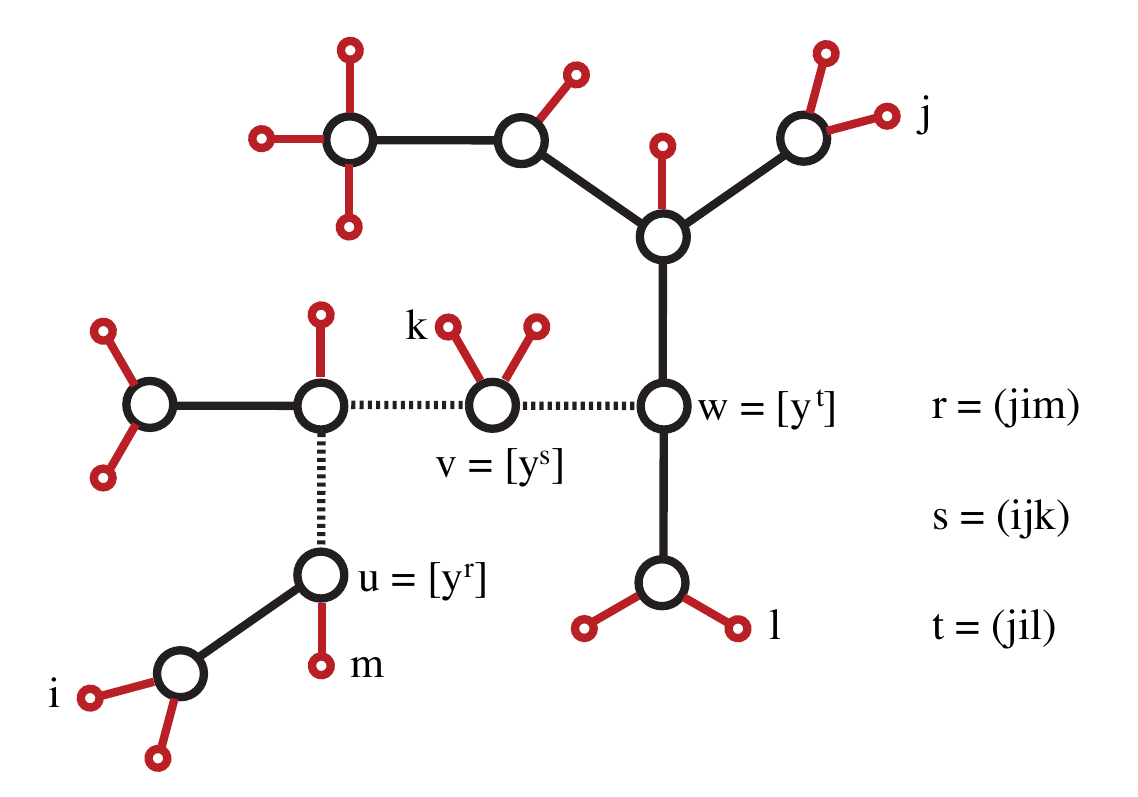}}}

\no{21} \def\viertelvorneuntree{{\includegraphics[width=.32\hsize]{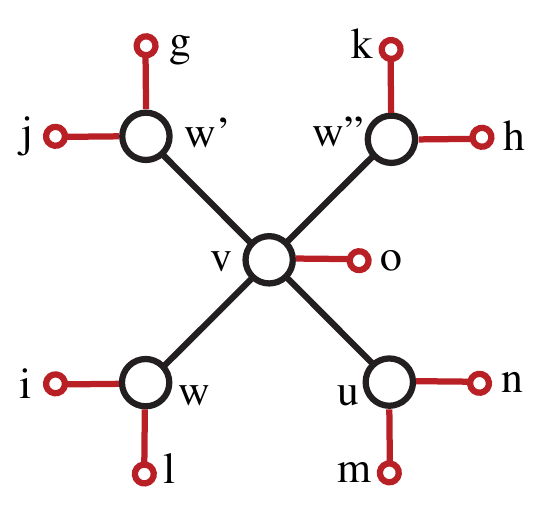}}}

\ignore

\def\neighbours{{\includegraphics[width=.6\hsize]{neighbours.pdf}}}

\def\transitivityequations{{\includegraphics[width=.47\hsize]{transitivityequations.pdf}}}

\def\viertelvorachttree{{\includegraphics[width=.32\hsize]{viertelvorachttree.pdf}}}

\def\phylogenetictreefiveleaves{{\includegraphics[width=.8\hsize]{phylogenetictreefiveleaves.pdf}}}

\def\projections{{\includegraphics[width=.5\hsize]{projections.pdf}}}

\def\setpartitions{{\includegraphics[width=.5\hsize]{setpartitions.pdf}}}

\def\arboralcovering{{\includegraphics[width=.4\hsize]{arboralcovering.pdf}}}

\def\orientedtree{{\includegraphics[width=.45\hsize]{orientedtree.pdf}}}

\def\orbitgraph{{\includegraphics[width=.6\hsize]{pictureorbitgraph.jpg}}}

\def\clustergraph{{\includegraphics[width=.6\hsize]{pictureclustergraph.pdf}}}

\def\picturesprays{{\includegraphics[width=.6\hsize]{picturesprays.pdf}}}

\def\phylogenetictreeA{{\includegraphics[width=.4\hsize]{picturephylogenetictreeA.pdf}}}

\def\phylogenetictreeB{{\includegraphics[width=.6\hsize]{picturephylogenetictreeB.pdf}}}

\def\phylogenetictreetwoverticesA{{\includegraphics[width=.5\hsize]{phylogenetictreetwoverticesA.pdf}}}

\def\phylogenetictreetwoverticesB{{\includegraphics[width=.6\hsize]{phylogenetictreetwoverticesB.pdf}}}

\def\phylogenetictreethreeverticesA{{\includegraphics[width=.7\hsize]{phylogenetictreethreeverticesA.pdf}}}

\def\phylogenetictreevertices{{\includegraphics[width=.6\hsize]{phylogenetictreevertices.pdf}}}

\def\complementarysprays{{\includegraphics[width=.6\hsize]{picturecomplementarysprays.pdf}}}

\def\deletingedgeA{{\includegraphics[width=.6\hsize]{picturedeletingedgeA.pdf}}}

\def\deletingedgeB{{\includegraphics[width=.6\hsize]{picturedeletingedgeB.pdf}}}

\def\edgesandclusters{{\includegraphics[width=.6\hsize]{pictureedgesandclusters.pdf}}}

\def\gluingedges{{\includegraphics[width=.7\hsize]{picturegluingedges.pdf}}}

\def\complementarysets{{\includegraphics[width=.3\hsize]{picturecomplementarysets.jpg}}}

\def\cuttingoffvertex{{\includegraphics[width=.7\hsize]{picturecuttingoffvertex.pdf}}}

\def\phylogenetictreethin{{\includegraphics[width=.7\hsize]{phylogenetictreethin.pdf}}}

\def\phylogenetictreeminimal{{\includegraphics[width=.6\hsize]{phylogenetictreeminimal.pdf}}}

\def\phylogenetictreenumeration{{\includegraphics[width=.7\hsize]{phylogenetictreenumeration.pdf}}}

\def\phylogenetictreenumerationspecial{{\includegraphics[width=.7\hsize]{phylogenetictreenumerationspecial.pdf}}}

\def\extremaltree{{\includegraphics[width=.8\hsize]{extremaltree.pdf}}}

\def\extremaltreeedgeselected{{\includegraphics[width=.8\hsize]{extremaltreeedgeselected.pdf}}}

\def\extremaltreeedgedeleted{{\includegraphics[width=.8\hsize]{extremaltreeedgedeleted.pdf}}}

\def\startree{{\includegraphics[width=.5\hsize]{startree.pdf}}}

\def\treeinductionA{{\includegraphics[width=.8\hsize]{treeinductionA.pdf}}}

\def\treeinductionB{{\includegraphics[width=.8\hsize]{treeinductionB.pdf}}}

\def\treeribislcontracted{{\includegraphics[width=.4\hsize]{treeribislcontracted.pdf}}}

\def\treelargeribisl{{\includegraphics[width=.4\hsize]{treelargeribisl.pdf}}}

\def\contractionoftrees{{\includegraphics[width=.4\hsize]{contractionoftrees.pdf}}}

\def\fibercombinatoricssixfive{{\includegraphics[width=.3\hsize]{fibercombinatoricssixfive.pdf}}}

\def\fibertreessixfive{{\includegraphics[width=.35\hsize]{fibertreessixfive.pdf}}}


\def\genericspecialtreefiveleaves{{\includegraphics[width=1\hsize]{genericspecialtreefiveleaves.pdf}}}

\def\compositiontreesA{{\includegraphics[width=.5\hsize]{compositiontreesA.jpg}}}

\def\compositiontreesB{{\includegraphics[width=.5\hsize]{compositiontreesB.jpg}}}

\def\fourcomponents{{\includegraphics[width=.5\hsize]{fourcomponents.jpg}}}

\def\phylogenetictreetriples{{\includegraphics[width=.7\hsize]{phylogenetictreetriples.pdf}}}

\def\phylogenetictreepairingtriples{{\includegraphics[width=.5\hsize]{phylogenetictreepairingtriples.jpg}}}

\def\destinationssimplepathslong{{\includegraphics[width=.9\hsize]{destinationssimplepathslong.png}}}

\recognize




\vs 1cm

\cl{\Bf PHYLOGENETIC TREES AND THE MODULI SPACE OF }\med

\cl{\Bf {\BFit n} POINTS ON THE PROJECTIVE LINE
%
\tex
{\baselineskip 12pt \footnote{}{\smallTimes MSC 2020: 14-02, 14D20, 14D22, 14H10, 05C05.
Supported by the Austrian Science Fund FWF within the projects P-31338 and DOI 10.55776/P34765. The second named author was additionally supported by the Austrian Science Fund within the project W1214-N15, DK9. We thank Michel Brion, Francis Brown, \'Eric Fusy, Matteo Gallet, Mahyar Hasani, and Gregor Kemper for valuable information and feedback.}}
{\footnote{
Supported by the Austrian Science Fund FWF within the projects P-31338 and DOI 10.55776/P34765. The second named author was additionally supported by the Austrian Science Fund within the project W1214-N15, DK9. We thank Michel Brion, Francis Brown, \'Eric Fusy, Matteo Gallet, Mahyar Hasani, and Gregor Kemper for valuable information and feedback.}}
} 
\big
\centerline{ Herwig HAUSER, Jiayue QI, Josef SCHICHO}\med
\vs .7cm

\centerline{\bf Abstract}\med

\centerline{{\baselineskip 12pt \vbox{\hsize 10.3cm \smallTimes  Choose $n$ pairwise distinct points $ x_1,...,x_n$ on the projective line $\P^1$ and submit them to the action of M\"obius transformations. The orbit space consisting of isomorphism classes of such $n$-tuples is classically denoted by $\MM_{0,n}$. Deligne, Mumford and Knudsen constructed in a series of celebrated papers a natural compactification $\ol\MM_{0,n}$ of $\MM_{0,n}$ by introducing the concept of $n$-pointed stable curves of genus zero. They then proved that $\ol\MM_{0,n}$ is a smooth projective variety and in fact represents a fine moduli space for isomorphism classes of such curves. \\
The present text offers an alternative approach to these constructions by using the geometric combinatorics of phylogenetic trees. The decisive clue for this method is to embed $\MM_{0,n}$ suitably into a large projective variety $(\P^1)^{n{n\choose 3}}$ and to take as its compactification simply the Zariski-closure $\XX_n$ of the image of $\MM_{0,n}$ therein. The points of $\XX_n$ are strings \med
\cl{{$\xx=(x^t)_t\in (\P^1)^{n{n\choose 3}}$}}
%
%
of $n$-tuples $x^t\in(\P^1)^n$, where $t$ runs over the a triples $t=(ijk)$ in $\{1,...,n\}$, subject to the condition that all these $n$-tuples have the same cross-ratio $\cross_q(x^t)=\cross_q(\xx)$, for every quadruple $q$ in $\{1,...,n\}$. One then associates to each string $\xx\in\XX_n$ a phylogenetic tree $\Gamma_\xx$, i.e., a finite graph without loops and no vertices of degree $2$. This tree carries precise combinatorial information about $\xx$; it allows one to prove that $\XX_n$ is smooth and irreducible, and it stratifies $\XX_n$ naturally in locally closed strata. There is a natural projection map $\pi:\XX_{n+1}\map \XX_n$ given by forgetting the components of strings which involve the index $n+1$. The fibers of this projection turn out to be $n$-pointed stable curves of genus zero, thus reproducing from scratch the concept of Deligne, Mumford and Knudsen. Actually, $\pi$ represents a universal family of such curves. This allows us to see that $\XX_n$ has all expected properties of a fine moduli space of $n$-pointed stable curves. It is thus canonically isomorphic to the Deligne-Mumford-Knudsen compactification $\ol\MM_{0,n}$.}}}\big\big


\cl{\Bf Part I. Introduction} \label{part:Introduction}\big


{\Bf 1. Four stories to start with} \label{section:fourstories}\med


(1) The group $\PGL_2$ acts on the projective line $\P^1$ by M\"obius transformations $z\map {az+b\over cz+d}$ and hence, for any $n\geq 1$, componentswise on the cartesian product $(\P^1)^n$. Restrict the action to the open subvariety $(\P^1)^n\sm\Delta_n$ consisting of {\it $n$-gons} $x=(x_1,...,x_n)$ with pairwise distinct entries $x_i\neq x_j$. The orbit space $\UU_n:=((\P^1)^n\sm \Delta_n)/\PGL_2$ of $\PGL_2$-equivalence classes of $n$-gons embeds {\it symmetrically} into the projective variety $(\P^1)^{n{n\choose 3}}$: To do so, choose for an orbit $[x]=[(x_1,...,x_n)]$ of an $n$-gon $x\in (\P^1)^n\sm\Delta_n$ and for every triple $t=(ijk)$ in $\{1,...,n\}$ a distinguished representative $x^t\in (\P^1)^n$ of $[x]$. This $n$-gon $x^t$ is uniquely prescribed by requiring to have at positions $i,j,k$ values $0$, $1$, $\infty$. The {\it string} $\xx=(x^t)_{t\, {\rm triple}}$ of all such representative $n$-gons of $[x]$ then defines an embedding $\Sigma_n:\UU_n\hookrightarrow (\P^1)^{n{n\choose 3}},\, [x]\map \xx,$ into a high-dimensional projective ambient variety. The Zariski-closure $\XX_n=\ol{\Sigma_n(\UU_n)}$ of the image in $(\P^1)^{n{n\choose 3}}$ is the object of interest we will study in this article.\med


(2) Denote by $\TT_n$ the closed subvariety of $(\P^1)^{n{n\choose 3}}$ of strings $\xx=(x^t)_t$ with $x^t_i=0$, $x^t_j=1$ and $x^t_k=\infty$ if $t=(ijk)$. Define $\YY_n\subseteq \TT_n$ as the closed subvariety of strings $\xx=(x^t)_t$ with equal cross-ratios $\cross(a,b,c,d)={(a-c)(b-d)\over (a-d)(b-c)}$ (see Section 9), i.e., such that $\cross_q(x^s)=\cross_q(x^t)$ for all quadruples $q$ in $\{1,...,n\}$ and all triples $s$ and $t$. By continuity, $\XX_n\subseteq \YY_n$. We prove that $\XX_n=\YY_n$. This describes the Zariski-closure of ${\Sigma(\UU_n)}$ in $(\P^1)^{n{n\choose 3}}$ by equations.\med

(3) To every $n$-gon $x^t$ of a string $\xx\in\YY_n$, associate its {\it incidence sets} $I^t\subset\{1,...,n\}$. Each $I^t$ collects the indices $i$ with equal entries $x^t_i$ of $x^t$. This gives a partition $\II^t$ of $\{1,...,n\}$. It only depends on the orbit $[x^t]$ of $x^t$. From the collection of all incidence partitions $\II^t$ of strings $\xx\in\YY_n$ one constructs the {\it incidence graph} $\Gamma_\xx$ of $\xx$. This is a finite planar tree. Its (inner) nodes (or: vertices) are the orbits $[x^t]$ of $n$-gons, its leaves (= outer nodes) are {\it singleton} incidence sets $I=\{i\}$. Edges between inner nodes correspond to complementary incidence sets $I\sqcup J=\{1,....,n\}$, and a leave $i$ is attached to an inner node $[x^t]$ if $\{i\}$ is a singleton incidence set for $x^t$ (see Section 13 for the precise construction). We prove that $\Gamma_\xx$ is a {\it phylogenetic tree} (Fig.~1). \med



\cl{\darwinsdrawing}

\cl{{Figure 1.} Darwin's famous drawing from 1837 of a phylogenetic tree.}\med\gb


{\it Example.} \label{example:stringfourgon} The string\med

\cl{$\xx=\left((0,1,\infty,a), (0,1,1-a,\infty),(0,{1\over 1-a},1,\infty), ({1\over a},0,1,\infty)\right)\in(\P^1)^{16}$,}\med

belongs to $\YY_4$ and has, for $a\neq 0,1,\infty$, only singleton incidence sets $I^t=\{i\}$, $i=1,...,4$. The orbits of all four $4$-gons are the same, and the phylogenetic tree $\Gamma_\xx$ consists of one vertex with four leaves, see Fig.~10, left. If, however, $a=0,1,$ or $\infty$, there are two singleton incidence sets and one with two elements. The four $4$-gons of $\xx$ define two orbits. Let us take for instance $a=0$, then\med

\cl{$\xx=\left((0,1,\infty,0), (0,1,1,\infty),(0,1,1,\infty), (\infty,0,1,\infty)\right)\in \YY_4$,}\med

with distinct orbits $[(0,1,\infty,0)]\neq [(0,1,1,\infty)]$, though both $4$-gons have the same cross-ratio. The phylogenetic tree is depicted in Fig.~10, right.\med


\cl{\fourlabelstree}

\cl{{\it Figure 10.} The two phylogenetic trees with four leaves.}\med\gb


(4) The last story concerns the appearance of stable curves in this context. Indeed, while the natural projection $\pi:\XX_{n+1}\map \XX_n$ is defined simply in terms of strings of $(n+1)$- and $n$-gons -- it forgets all entries involving the index $n+1$ (see Section 11) --, the more sophisticated concept of $n$-pointed stable genus zero curves appears a posteriori as the fibers $\pi\inv(\xx)$ of $\pi$, the marked points on them given by well chosen sections $\sigma_p:\XX_n\map\XX_{n+1}$ of $\pi$. In this way, $\pi$ becomes the universal family of $n$-pointed stable curves. This closes the circle of thoughts: Classically, one starts with the definition of $n$-pointed stable curves of genus zero and then proceeds by constructing a fine moduli space for them. In the present paper, we go in the opposite direction, starting with the space $\XX_n$ of strings as a natural compactification of $\UU_n$, and then recover naturally the notion of $n$-pointed stable curves together with their moduli space $\ol{\MM_{0,n}}$. \med\gb

These four stories shall serve as a gentle approach to the theory of Deligne-Mumford and Knudsen-Mumford about the existence and shape of $\ol{\MM_{0,n}}$, while avoiding the use of advanced tools from algebraic geometry. Our method is to exploit the combinatorial geometry of phylogenetic trees. Doing so, an astonishing phenomenon happens: Imagine for instance that one wants to show the smoothness of the variety $\XX_n$, or that the fibers of $\pi:\XX_{n+1}\map\XX_n$ are stable curves. For these tasks, it turns out that the geometric features of the associated phylogenetic trees $\Gamma_\xx$ of strings $\xx$ in $\XX_n$ serve like an {\it instruction manual} which tells one  how to design the proof: Performing elementary operations with the tree - deleting leaves, contracting and inserting edges, cutting the tree in different components, travelling along paths in the tree - one can develop almost instinctly a strategy for proving the required assertions.%
\tex
{{\baselineskip 12pt \footnote{$\!\!{}^1$}{\smallTimes The authors concede to have rarely experienced this form of ``exterior assistance'' when trying to prove theorems.}}}
{\footnote{The authors concede to have rarely experienced this form of ``exterior assistance'' when trying to prove theorems.}}
%
%
One emphasis of the present article will be to explain how to ``read'' this manual and how to profit from it.\med\big


\hs 4cm\vbox {\smallTimes \baselineskip 12pt

{\bf Part I: Introduction}

1. Four stories to start with 

2. The symmetrization of $n$-gons by strings

3. Upshot of main result

4. Recap on Deligne-Mumford and Knudsen-Mumford

5. Stable curves and moduli spaces

6. The dual graph of an $n$-pointed stable curve\med}

\hs 4cm\vbox {\smallTimes \baselineskip 12pt

{\bf Part II: Constructions}

7. Cross-ratios

8. Strings and $n$-gons

9. The varieties $\TT_n$, $\XX_n$ and $\YY_n$

10. Limits of orbits

11. The projection map $\pi_a:\XX_{n+1}\map \XX_n$

12. Phylogenetic trees

13. The phylogenetic tree of a string

14. The phylogenetic tree of a stable curve 

15. Constructing a string from a stable curve 

16. The main theorem about $\XX_n$\med}

\hs 4cm\vbox {\smallTimes \baselineskip 12pt

{\bf Part III: Proofs}

17. The smoothness of $\XX_n$

18. The stratification of $\XX_n$

19. The boundary divisor $\BB_n$

20. Constructing a stable curve from a string


21. Sections of $\pi_a:\XX_{n+1}\map \XX_n$

22. The equality of $\XX_n$ with $\YY_n$

23. The isomorphism between $\XX_n$ and $\ol{\MM_{0,n}}$.

24. Proof of the main theorem\med

Index of notation

References

}
\med \med\gb


{\Bf 2. The symmetrization of {\Bfit n}-gons by strings}\label{section:symmetrization}\med

When compactifying the moduli space $\MM_{0,n}$ of $n$ pairwise distinct points $x_1,...,x_n$ in $\P^1$ under the $\PGL_2$-action, one has to develop a suitable concept of limit as some of the points come together and coalesce. Following the concept proposed by Grothendieck in {\cite{\[SGA7-I\]}}, see the citation in {\cite{\[FMP\]}}, p.~189, Deligne and Mumford used stable curves to define such limits {\cite{\[DM\]}}: instead of just letting move the points on $\P^1$ while they come closer, let also move and vary $\P^1$ itself, that is, take a family of rational curves whose generic member is irreducible and smooth and hence isomorphic to $\P^1$, but which may specialize to a union of $\P^1$'s at certain points (think of a family of hyperbolas in $\R^2$ degenerating to the two coordinate axes). This means to consider morphisms $X\map S$ whose generic fiber is isomorphic to $\P^1$, but whose special fibers may be a union of transversal rational curves, together with $n$ disjoint sections $\sigma_1,...,\sigma_n:S\map X$ selecting in each fiber the $n$ marked points. To make the whole concept work, some technical modifications and assumptions are necessary. They will be described in a later section.\med

The concept of $n$-pointed stable curves of genus zero then defines the required compactification $\ol\MM_{0,n}$ of $\MM_{0,n}$. This space is a smooth, irreducible projective variety and a fine moduli space for isomorphism classes of $n$-pointed stable curves {\cite{\[DM, Knu1\]}} (see Section 5 for the definition of a fine moduli space). In the present paper, we propose another limit construction for $n$ pairwise distinct points in $\P^1$. At the end, stable curves will show up, and we will see that our compactification $\XX_n$ is in fact isomorphic to $\ol\MM_{0,n}$.\med\gb

The idea is very simple and beautiful; it has appeared in the literature in various disguises, mostly using cross-ratios, see for instance {\cite{\[GHP\]}}, (1.3), p.~133, {\cite{\[Br\]}}, Sec.~2, p.~381, {\cite{\[SchT\]}}, p.~1, {\cite{\[Sin, HKT2\]}}. We will propose and pursue a slightly different procedure than the ones in these references, compelling the symmetry of all our constructions:  The group $\PGL_2$ acts on $\P^1$ by M\"obius transformations, \med

\cl{$\ds z\map A\cdot z ={az+b\over cz+d}$}
for 
\tex
{$\ds A=\left(\matrix {a & b\cr c & d}\right)\in \PGL_2$}
{$\ds A=\left(\begin{array}{rr} a & b\\ c & d\end{array}\right)\in \PGL_2$}
and $z\in \P^1$, with the obvious rules of calculus for the value $\infty$. The action is sharply $3$-transitive: For any two $3$-gons $x=(x_1,x_2,x_3)$ and $y=(y_1,y_2,y_3)$ with pairwise distinct entries there exists a unique matrix $A\in\PGL_2$ sending $x$ to $y$. Let $\Delta_n$ denote the big diagonal in $(\P^1)^n$ consisting of $n$-gons with at least two equal entries. Then, for every $n$-gon $x=(x_1,...,x_n)\in (\P^1)^n\sm\Delta_n$ with pairwise distinct entries, the induced action of $\PGL_2$ on $(\P^1)^n$ (acting componentswise) allows one to transform $x$ into an $n$-gon $y=(y_1,...,y_n)$ whose entries at three specified places, say, $i,j,k$, have prescribed values, for instance, $y_i=0$, $y_j=1$, and $y_k=\infty$. This $n$-gon $y$ is then unique and a {\it distinguished representative} of the $\PGL_2$-orbit $[x]$ of $x$. It clearly depends on the choice of the triple $t=(ijk)\in {N\choose 3}$, with $N=\{1,...,n\}$ and $N\choose 3$ the set of triples in $N$. We write $y=x^t=(x^t_1,...,x^t_n)\in (\P^1)^n$, with\med

\cl{$x^t_i=0$, $x^t_j=1$, $x^t_k=\infty$.}\med\gb

Then define the map

\cl{$\Sigma_n:(\P^1)^n\sm\Delta_n\map (\P^1)^{n{n\choose 3}}$,}\med

\cl{$x\map \xx=(x^t)_{t\in {N\choose 3}}$,}\med

by sending an $n$-gon $x$ to its {\it string $\xx$} listing all distinguished representatives $x^t$ (in an arbitrary order). This map is constant on $\PGL_2$-orbits, and injective on the set of orbits, sending different orbits to different strings. Therefore, by construction, the map $\Sigma_n$ passes to the quotient, thus defining an injective map on the space of $\PGL_2$-orbits,\med

\cl{$\Sigma_n: \UU_n=((\P^1)^n\sm\Delta_n)/\PGL_2\hookrightarrow (\P^1)^{n{n\choose 3}}$,}\med

\cl{$[x]\map \xx=(x^t)_{t\in {N\choose 3}}$,}\med\gb

called the {\it symmetrization map of $n$-gons}. It sends orbits to strings -- and thus represents equivalence classes of $n$-gons by points of the projective variety $(\P^1)^{n{n\choose 3}}$. The advantage of this construction instead of picking just a single triple $t$, typically $t=(123)$, and the $n$-gon $x^t$ as a representative of $[x]$, lies in the invariance of $\Sigma_n$ under permutations of $1,...,n$. This symmetry will become very beneficial later on.\med

We may now define the ``limit'' of a family of $n$-gons $x$, or of orbits $[x]$, as the respective limit of the images $\xx$ in $(\P^1)^{n{n\choose 3}}$. Said differently, we define the compactification $\XX_n$ of $\UU_n$ as the {\it Zariski-closure} of the image of $\UU_n$ in $(\P^1)^{n{n\choose 3}}$,\med

\cl{$\XX_n=\ol{\Sigma_n(\UU_n)}\subset (\P^1)^{n{n\choose 3}}$.}\med

For every $n\geq 3$, this is by construction an irreducible projective variety of dimension $n-3$. It comes with projection maps $\pi:\XX_{n+1}\map\XX_n$, given as the restriction of the map $(\P^1)^{(n+1){n+1\choose 3}}\map(\P^1)^{n{n\choose 3}}$ which forgets the components of strings $\xx$ involving the index $n+1$.\med

\med\med\gb


{\Bf 3. Upshot of main result}\label{section:upshottheorem}\med

Here is the outline of what will be explained in this text, see Section 16 for a detailed statement and compare with {\cite{\[DM\]}}, Thm.~2.7, {\cite{\[Knu1\]}}, Thm.~5.2. Set $N=\{1,...,n\}$ and let $\Delta_n$ be the big diagonal in $(\P^1)^n$.\med\gb


{\bf Theorem.} \label{theorem:upshottheorem} {\it Let $\PGL_2$ act on $(\P^1)^n$ componentswise by M\"obius transformations. Denote by $\XX_n$ the Zariski closure of the image of $\UU_n=((\P^1)^n\sm \Delta_n)/\PGL_2$ in $(\P^1)^{n{n\choose 3}}$ under the symmetrization map $\Sigma_n:\UU_n\map (\P^1)^{n{n\choose3}}$ sending $\PGL_2$-orbits $[x]$ to their string $\xx=(x^t)_{t\in {N\choose 3}}$ of distinguished $n$-gons $x^t\in (\P^1)^n$. Associate to each $\xx$ a phylogenetic tree $\Gamma_\xx$ with $n$ leaves, given as the incidence graph of $\xx$.}

(1) {\it The variety $\XX_n$ is a smooth, closed and irreducible subvariety of $(\P^1)^{n{n\choose 3}}$ of dimension $n-3$. A set of defining equations for $\XX_n$ in $(\P^1)^{n{n\choose 3}}$ is given by the equality of cross-ratios $\cross_q(x^s)=\cross_q(x^t)$ between the $n$-gons $x^t$ of a string $\xx\in(\P^1)^{n{n\choose 3}}$, for $q=(ijk\ell)$ a quadruple of numbers in $N$. }
 
(2) {\it The variety $\XX_n$ carries a natural stratification whose locally closed strata $\SS_T$ consist of strings $\xx$ with the same phylogenetic tree $\Gamma_\xx=T$; the open dense stratum is $\UU_n$ and corresponds to the generic phylogenetic tree $T_*$ with $n$ leaves; a stratum $\SS_{T'}$ lies in the closure of a stratum $\SS_T$ if and only if the tree $T$ can be obtained from $T'$ by the contraction of edges.}

(3) {\it The boundary $\BB_n=\XX_n\sm\UU_n$ of $\XX_n$ is a simple normal crossings divisor.}

(4) {\it The projection $\pi: \XX_{n+1} \map \XX_n$ given by forgetting all entries involving the last index $n+1$ is a flat projective morphism of algebraic varieties with one-dimensional reduced fibers.}

(5) {\it There is a natural isomorphism between the Deligne-Mumford compactification $\ol\MM_{0,n}$ of $\MM_{0,n}$ and $\XX_n$: It sends an $n$-pointed stable curve $C$ to a string of $n$-gons $\xx=(x^t)_{t\in{N\choose 3}}$ given as the contraction of $C$ with respect to a median component $C_t$ of $C$, and, conversely, a string $\xx\in\XX_n$ to the stable curve given as the fiber $\pi\inv(\xx)$ under $\pi:\XX_{n+1}\map\XX_n$. Under this isomorphism, the augmented dual graph $\Gamma_C$ of a stable curve $C$ equals the phylogenetic tree $\Gamma_\xx$ of its image string $\xx$.}

(6) {\it The variety $\XX_n$ represents a fine moduli space for isomorphism classes of $n$-pointed stable curves of genus zero.}

(7) {\it The projection map $\pi: \XX_{n+1} \map \XX_n$ together with $n$ suitably chosen sections $\sigma_p$ represents a universal family for isomorphims classes of $n$-pointed stable curves.} \med


{\it Remarks.} (a) Various definitions and constructions of a moduli space of stable curves as a Zariski-closure appear in the literature, e.g.~{\cite{\[Br, GHP, Sin\]}}. It is shown in {\cite{\[GHP\]}}, using methods inspired by Knudsen's arguments and different from ours, that a space $B_n$ constructed from the $n$ points of a stable curve is isomorphic to $\ol\MM_{0,n}$ and thus a fine moduli space of $n$-pointed stable curves of genus zero. In {\cite{\[Br\]}}, Section 2.1 and 2.8, an a priori asymmetric construction  associates to equivalence classes of orbits of $n$-gons just one distinguished $n$-gon (with respect to the triple $t=(123)$), and embeds it then into $(\P^1)^{n\choose 4}$ by taking all its cross-ratios, thus getting a locally closed subvariety ${\frak M}_{0,n}$. The smoothness of its Zariski-closure $\ol{\frak M}_{0,n}$ and further properties are proven in Thm.~2.25 and its corollary Cor.~2.32. The isomorphism of $\ol{\frak M}_{0,n}$ with $\ol\MM_{0,n}$ is only mentioned {\cite{\[Br\]}} Section 2.9. See also {\cite{\[GML, KT1, KT2\]}} for further descriptions of $\ol\MM_{0,n}$.

(b) The systematic use of phylogenetic trees associated to strings $\xx$ in $\XX_n$ seems to have been largely neglected in the literature so far. They mostly appear as the dual graph of stable curves, without pursuing their combinatorial structure (but see {\cite{\[Kapr1, GML\]}} for another connection). In the present text, in contrast, they play a central role and are defined from scratch, just using the incidence relations between the entries of a string $\xx$: as mentioned earlier, they are a valuable source of information to design proofs, and, moreover, they naturally stratify $\XX_n$.

(c) Experimental studies suggest that graphs similar to -- but more complicated than  -- phylogenetic trees can also be associated to strings of $n$-gons in the projective plane $\P^2$. Analogous phenomena as in the case of points on the projective line are observed in various examples but still lack a deeper understanding. It seems that projective geometry and duality find here a combinatorial expression which might be helpful for studying moduli problems and configuration spaces in the spirit of Fulton-MacPherson {\cite{\[FMP\]}}.\med

\med\gb

{\Bf 4. Recap on Deligne-Mumford and Knudsen-Mumford} \label{section:recapdkm}\med

Let us briefly revise the main aspects of the abundant literature about moduli spaces of $n$-gons in $\P^1$ (see {\cite{\[Behr, Cap, Cav, Kol, Maz\]}} for general information about moduli spaces): Deligne and Mumford proved that the moduli space $\ol \MM_{0,n}$ of $n$-pointed stable curves is irreducible and a smooth stack {\cite{\[DM\]}}. Knudsen then showed that $\ol \MM_{0,n}$ is even a smooth projective variety {\cite{\[Knu1\]}}. For points in higher dimensional projective spaces, Gelfand and MacPherson associated to $\PGL_{d+1}$-orbits of $n$-gons in $\P^d$ matroids and the matroid polytope {\cite{\[GMP\]}}, see {\cite{\[Kapr1\]}} for a succinct description. Kapranov mentions in the introduction that the action of the maximal torus in $\GL_n$ on the Grassmannian $\Grass(d+1,n)$ of $(d+1)$-planes in $n$-space is equivalent to the study of $\PGL_{d+1}$-orbits of $n$-gons in $\P^d$. In Section 1.3, he defines phylogenetic trees and proves in Thm.~1.3.6 that these are in bijection with tilings of the simplicial polytope by matroid polytopes.\med

 
In \tex{[DM]}{\cite{\[DM\]}}, def.~1.1, p.~76, Deligne and Mumford recall the definition of a (family of) stable curves of genus $g$, originally proposed by Grothendieck in {\cite{\[SGA7-I\]}}, as a proper flat morphism $C\map S$ of algebraic varieties whose geometric fibers are reduced, connected normal crossings curves of genus $g$ and where each irreducible rational component of a fiber $C_s$ meets at least three other irreducible components (these items will be explained in the course of the article). 
%
\ignore
[\ooo we don't understand condition (ii) there; the geometric fiber is the fiber $C_s=C\times_S\Spec\, \ol \kk$, where $\Spec\, \kk\map S$ defines the point $s\in S$ and $\ol \kk$ denotes the algebraic closure of $\kk$, see https://mathoverflow.net/questions/121932/geometric-fibers-of-schemes].
\recognize
%
Note that Deligne-Mumford do not talk about {\it $n$-pointed} stable curves, say, families of curves $ C\map S$ with sections of $S\map C$, since the considered curves are supposed to have genus $g\geq 2$ (for which stability is ensured without marking points).\med

On page 86 of {\cite{\[DM\]}}, the authors associate the {\it dual graph} $\Gamma$ to a stable curve $C$. Its vertices are the irreducible components of $C$, and two vertices are connected by an edge if the two components intersect. They only use this graph once in the paper, and only to a very small amount. They don't consider {\it labels} and {\it leaves} as we will do later on.\med\gb

In the introduction of {\cite{\[DM\]}}, the irreducibility of the moduli space $\MM_g$ of curves of genus $g$ in characteristic $0$ is claimed to be classical, proven by Enriques-Chisini. Deligne and Mumford, however, affirm to be closer to the (incomplete) proof of Severi {\cite{\[Sev\]}}. They complete the gap and base the proof on the {\it Stable Reduction Theorem}.\med\gb
\ignore
[DM, Lemma 1.16, p.~86]: Any automorphism of a stable curve which fixes the dual graph is the identity. Whereas at the beginning [DM, p.~76] only curves of genus $g\geq 2$ are considered, the proof of the lemma considers curves of all genera.\med\recognize

The Knudsen-Mumford compactification $\ol{\MM}_{0,n}$ of $\MM_{0,n}$ of stable $n$-pointed genus zero curves is proven in {\cite{\[Knu1\]}} to be a smooth projective variety whose boundary $\ol{\MM}_{0,n}\sm\MM_{0,n}$ is a normal crossings divisor and represents a fine moduli space, see also {\cite{\[KM, Kapr1, Hass, KT1, KT2, Keel, MR, HKT1, Kol, SchT\]}}. The boundary divisor $\ol{\MM}_{0,n}\sm\MM_{0,n}$ admits a canonical stratification whose strata are locally closed and consist of $n$-pointed stable curves with the same (augmented) dual graph. \med

Keel gives a quite concise and understandable account of Knudsen's paper {\cite{\[Keel\]}}. Kapranov {\cite{\[Kapr1\]}}, Thm.~4.3.3, Keel {\cite{\[Keel\]}}, and Fulton-MacPherson {\cite{\[FMP\]}} describe $\ol{\MM_{0,n}}$ as certain blowups of $(\P^1)^{n-3}$. Gerritzen, Herrlich and van der Put have given in {\cite{\[GHP\]}} an interpretation of $\ol{\MM_{0,n}}$ in terms of spaces of cross-ratios, see also {\cite{\[FMP\]}}, p.~189. 
Further useful references are {\cite{\[Br, GP, Hass, HK, Kol, MR, SchT, Tem, Vor\]}}.\med

\med\gb


{\Bf 5. Stable curves and moduli spaces}\label{section:stablecurves}\med

In this article, an {\it $n$-pointed stable curve of genus} $0$ (over a field) is a one-dimensional, reduced and connected but possibly reducible variety $C$ defined over a fixed ground field $\kk$ with $n$ {\it marked points} $p_1,...,p_n$ on it such that the following holds: the irreducible components $C_i$ of $C$ are smooth rational curves (i.e., isomorphic to the projective line $\P^1$) any two of which meet transversally (i.e., like coordinate axes) and such that no cycles of pairwise intersections are created; the points $p_i$ lie outside the intersection points of the components and are hence smooth points of $C$; and each component of $C$ has at least three {\it special points}, that is, either intersection points with other components (= the singular points of $C$) or marked points $p_i$. See Fig.~5.\med


For $n\geq 4$, there are infinitely many isomorphism classes of $n$-pointed stable curves, since moving a marked point on $\P^1$ while keeping the other points fixed changes the cross-ratio -- but cross-ratios are an invariant of stable curves under isomorphism.\med


\cl{\stablecurve}\vs-.3cm

\cl{{\it Figure 5.} The example of a stable curve with $7$ components and $11$ marked points.}\med\gb


A {\it family of $n$-pointed stable curve of genus zero}  over a base variety $S$ is a flat and proper morphism $\varphi:\CC\map S$ such that each fiber $C_s=\varphi\inv(s)$ is an $n$-pointed stable curve with marked points $\sigma_1(s),...,\sigma_n(s)$ given by $n$ disjoint sections $\sigma_p:S\map\CC$ of $\varphi$ avoiding the singular points of the fibers.
Two families $\varphi:\CC\map S$ and $\psi:\DD\map S$ of $n$-pointed stable curves over the same base variety $S$ and with sections $\sigma_1,\dots,\sigma_n:S\map \CC$ and $\tau_1,\dots,\tau_n:S\map\DD$, respectively, are {\it isomorphic} if and only if there is an
isomorphism $ \chi :\CC\map\DD$ such that $\varphi=\psi\circ\chi$ and $\chi\circ\sigma_p=\tau_p$ for $p=1,\dots,n$. This implies that for all $s\in S$,
the fibers $\varphi^{-1}(s)$ and $\psi^{-1}(s)$ are isomorphic, by an isomorphism taking the marked points to the corresponding marked points. \med

Informally, a moduli space is an algebraic variety such that any point stands for a unique isomorphism class. We just defined families of stable curves $\CC\map S$, but, similarily, we can also define, for any variety $X$, ``points of $X$ over $S$'': these are just morphisms $S\map X$ in the category of varieties (for $S$ a point, the image of $S$ is just a "usual" point of $X$). This concept of generalised points fits nicely to the idea that varieties
are often considered as systems of algebraic equations whose solutions can be considered in varying domains. Typically, a point in $X$ over $\C$
is a solution with complex coordinates, a point in $X$ over $\Q$ is a solution with rational coordinates, and a point over the line $\A^1$ corresponds to a parametrized curve in $X$. \med

If $\MM$ is a moduli space for stable $n$-pointed curves of genus zero, then a point in $\MM$ over $S$, i.e., a morphism $S\map \MM$, ``stands for'' 
an isomorphism class of families of stable $n$-pointed curves over $S$. But what exactly does ``standing for'' mean? With the concept of a universal
family, we can make this very precise. A {\it universal family} is a family $\Phi:\wt\MM\map \MM$ of $n$-pointed stable curves over $\MM$, together with sections
$\sigma_1,\dots,\sigma_n:\MM\map \wt\MM$, with the following distinctive property: For any morphism $f:S\map \MM$, let $\Phi_f:S\times_\MM \wt\MM\map S$ be the pullback of $\Phi$ along $f$, and let $\tau_{p,f}:=\sigma_p\circ f$ for $p=1,\dots,n$. It is straightforward to see that $\Phi_f$ with sections
$\tau_{1,f},\dots,\tau_{n,f}$ is a family of $n$-pointed stable curves of genus zero. Universality means that for every family $\varphi:\CC\map S$ 
of $n$-pointed stable curves of genus zero, there is a unique morphism $f:S\map \MM$ such that $\varphi$ is isomorphic to $\Phi_f$. So the expression ``standing for'' has an intrinsic meaning. A moduli space that  possesses a universal family is called a {\it fine moduli space}. See e.g.~{\cite{\[Behr, Cav\]}} for more background on this.\med\gb

In some comparable situations, for instance for stable curves of genus $g>0$, it is not possible to construct a fine moduli space. One can
still obtain a so-called {\it course moduli space} satisfying a weaker condition. However, all moduli spaces in this paper are fine moduli spaces.
\med\gb



{\it Example.} \label{example:npointedstablecurve} 
Let $C \subset \P^2 \times \A^1$ be the subvariety defined by
$xy-yz+t(xy-xz)=0$, together with the projection $C\map\A^1$ sending $((x:y:z),t)$ to $t$. The fibers $C_t$ are irreducible conics for $t\neq 0,-1$.
In case $t=0$ we get two lines intersecting in $(x:y:z)=(1:0:1)$,
and in case $t=-1$ we get two different lines intersecting in $(1:1:0)$.
Then choose the constant sections $t\map (1:0:0)$, $t\map (0:1:0)$, $t\map (0:0:1)$, $t\map (1:1:1)$ to mark four points on the fibers $C_t$ .\med\med

\ignore

{\it Example.} \label{example:npointedstablecurve} [\ooo Still to be checked, see mail JS] The simplest example of a family of $n$-pointed stable curves is a family of hyperbolas in $\P^1\times \P^1$, marked with four points, deforming into the union of the two coordinate axes. More precisely, consider the variety $\CC$ in $\P^1\times \P^1\times \A^1$ given as the closure of the affine curve defined by the equation\med

\cl{$xy-t =0$}\med

in $\A^1\times \A^1\times \A^1$. [\ooo For $t=\infty$, it is not clear of how to define the curve, hence we take the base $S=\A^1$.] Map it to the affine line $\A^1$ by sending $(x,y,t)$ to $t$. The fibers are, for $t\neq 0$, smooth quadratic curves in $\P^1\times \P^1$, while, for $t=0$, we obtain the union $(\P^1\times \{0\})\cup (\{0\}\times \P^1)$ of the two axes. For the sections, one may choose $t\map (1,t,t),(t,1,t),(-1,-t,t),(-t,-1,t)$ [\ooo assuming characteristic $\neq 2$.]\med\med

\recognize


{\Bf 6. The phylogenetic tree of an $n$-pointed stable curve}\label{section:treeofcurve}\med

Phylogenetic trees appear naturally when studying stable curves. The {\it augmented dual graph} of an $n$-pointed stable genus zero curve $C$ is defined as follows: It is a finite undirected graph $\Gamma$ without loops and multiple edges whose nodes have either degree $1$ or degree $\geq 3$; the first are the {\it leaves} of the tree, corresponding to the $n$ marked points of $C$, the second are the {\it inner vertices}, they are in bijection with the irreducible components of $C$. An edge joins two inner vertices if the respective components of the curve intersect, while each leaf is connected by an edge to exactly one inner vertex, corresponding to the component on which the point sits. There are no edges between leaves. These conditions are precisely the axioms of a phylogenetic tree with $n$ leaves: a finite undirected graph with no vertices of degree $2$ and precisely $n$ (labelled) vertices of degree $1$, see Fig.~6. 


\cl{\phylogenetictreeexample}\vs-.3cm

\cl{{\it Figure 6.} The example of a phylogenetic trees with $n$ leaves.}\med\gb


\med


{\it Plan of the article.} Part II, Sections 7-16, provides the construction and properties of the main players (cross-ratios, strings, trees, partitions, and operations between them) and gives in Section 16 the precise statement of the theorem. Part III is devoted to the proof of the various assertions of the theorem. The three most exigent (and also most interesting) parts are the smoothness of $\XX_n$ (Section 17), the description of the boundary divisor $\BB_n=\XX_n\sm\UU_n$ (Section 19), and the construction of stable curves from string as the fibers of the projection map $\pi: \XX_{n+1}\map\XX_n$ (Section 20).\med

In Section 24 we construct the one-to-one correspondence between $\XX_n$ and $\ol{\MM_{0,n}}$. To show that this is in fact an isomorphism of algebraic varieties would go beyond the scope of this article, as it uses dualizing sheaves and higher image sheaves. We give, however, precise references. The problem here is to show that any $n$-pointed stable curve $\gamma:\CC\map S$ admits a (unique) morphism $S\map \XX_n$ such that $\gamma$ appears as the pullback to $S$ of the universal curve $\pi:\XX_{n+1}\map\XX_n$. The proof requires in particular to show that families of isomorphic curves having sufficiently many sections can be {\it trivialized}, i.e., seen as the projection to $S$ of a cartesian product $S\times \CC_s$ of $S$ with one special fiber $\CC_s$. This is a highly non-trivial result from algebraic geometry for which no elementary argument seems to be applicable.\med

The text is designed to be accessible and appealing also for people outside algebraic geometry. It is worth to sit down with pencil and paper to draw the geometric situations, and then, hopefully, the charm of the reasoning will become transparent.\med\big\gb


\cl{\Bf Part II: The space $\XX_n$ of strings of $n$-gons}\label{part:spacexxn}\med\med

The construction of the space $\XX_n\subseteq (\P^1)^{n{n\choose 3}}$ of strings of $n$-gons $\xx=(x^t)_{t\in{N\choose 3}}$ requires a few basic and mostly classical concepts from invariant theory, respectively, projective geometry. We give a brief summary.\med\gb

\med

{\Bf 7. Cross-ratios}\label{section:crossratios}\med

Let $N$ be a finite totally ordered set of cardinality $n\geq 4$. Typically one may take $N=\{1,...,n\}$ with the natural ordering $1<2<\ldots <n$. The elements of $N$ are called {\it labels}. A {\it triple} in $N$ is a three-point subset $t$ of $N$, i.e., an {\it unordered} three-tuple of pairwise distinct elements of $N$. We write it as $t=(ijk)$ where we have arranged the entries $i,j,k\in N$ such that $i<j<k$. A {\it quadruple} in $N$ is an {\it ordered} four-tuple $q=(ijk\ell)$ in $N^4$ with pairwise distinct entries. For quadruples, we do not assume that the entries are listed increasingly. The sets of triples and quadruples in $N$ are denoted by ${N\choose 3}$ and $N^{\ul 4}$, respectively.\med


For variables $\xi_1,...,\xi_n$ and a quadruple $q=(ijk\ell)$ in $N^{\ul 4}$, the {\it formal cross-ratio} $[ijk\ell]$ is defined as the element\med

\cl{$\ds [ijk\ell](\xi)={(\xi_i-\xi_k)(\xi_j-\xi_\ell)\over (\xi_i-\xi_\ell)(\xi_j-\xi_k)}\in \kk(\xi_1,...,\xi_n)$}\med

in the field $\kk(\xi_N):=\kk(\xi_1,...,\xi_n)$ of rational functions. We often just write $[ijk\ell]$ or even $[q]$ and drop the word ``formal''. As the formal cross-ratios are quotients of homogeneous polynomials of the same degree, they belong to the function field $K((\P^1)^n)$ of $(\P^1)^n$. The cross-ratios satisfy the relations\med

\cl{$[ijk\ell]=[ji \ell k] =[k\ell ij]= [\ell k ji]=1-[ikj\ell]$,}\med

\cl{ $[ijk\ell]=[jik\ell]\inv=[ij \ell k]\inv$,}\med

and, for five distinct labels $i,j,k,\ell,m$ in $N$, the {\it triple product formula} \med

\cl{$[ijk\ell]\cdot [ij\ell m]\cdot [ijmk]=1$.}\med

This last formula can be seen as a {\it cocycle condition}.%
%
\tex
{{\baselineskip 12pt \footnote{$\!\!{}^2$}{\smallTimes In \cite{\[GHP\]}, (1.4), p.~135, the same relations are used to define the space $B_n$ as a closed subvariety of a projective ambient variety.}}}
\footnote{In \cite{\[GHP\]}, (1.4), p.~135, the same relations are used to define the space $B_n$ as a closed subvariety of a projective ambient variety.}
%
It will be most often used in the form\med

\cl{$\ds [ijk\ell]=  [ijkm]\cdot  [ijm\ell]$.}\med\gb


If $r$ is the quadruple obtained by a permutation of the entries of $q=(ijk\ell)$, the  cross-ratio $[r]$ is one of the following rational functions in the cross-ratio $[q]$ of $q$,\med

\cl{$\ds [r]=[q], {1\over [q]}, 1-[q], {1\over {1-[q]}}, {[q]\over{[q]-1}}, {{[q]-1}\over[q]}$.}\med

These six functions form a group isomorphic to the permutation group $S_3$ on three elements. Cross-ratios are invariant under the action of $\PGL_2$ on $K((\P^1)^n)$: To see this, it suffices to consider the transformations $\xi_i\map c\xi_i$, $\xi_i\map\xi_i+c$, $c\in K$, and $\xi_i\map{\1\over \xi_i}$. In all cases, the cross-ratio does not change.\med

Note that when restricting the formal cross-ratios to $\xi_i=\xi_j$ or $\xi_k=\xi_\ell$ one gets the constant $1$, and restricting  to $\xi_i=\xi_k$ or $\xi_j=\xi_\ell$ one gets the constant $0$. The restriction to $\xi_i=\xi_\ell$ or $\xi_j=\xi_k$  would yield $\infty$ and is therefore not defined as an element of $K((\P^1)^n)$.\med\gb


Formal cross-ratios $[ijk\ell]$ define {\it cross-ratio functions}\med

\cl{$\cross_{(ijk\ell)}:(\P^1)^n\sm \nabla^3_n(ijk\ell)\map \P^1$,}\med

via 

\cl{$\ds\cross_q(x)=\cross_{(ijk\ell)}(x)={(x_i-x_k)(x_j-x_\ell)\over (x_i-x_\ell)(x_j-x_k)}$,}\med

where $\nabla^3_n(ijk\ell)$ denotes the set of $n$-gons $x$ for which (at least) three of the entries $x_i,x_j,x_k,x_\ell$ are equal. We call $\cross_q(x)$ the {\it evaluation of the cross-ratio} in $x$, or simply the {\it cross-ratio} of $x$. It is defined whenever no three of the entries $x_i$, $x_j$, $x_k$, $x_\ell$ are equal (with the obvious rules when the denominator becomes $0$ or when some entries are $\infty$). To be more precise, one may define the cross-ratio in terms of projective coordinates $x_i=(a_i:b_i)$ in $\P^1$, taking then the affine charts $a_i\over b_i$, respectively, $b_i\over a_i$. In the first chart, with $b_i,b_j,b_k,b_\ell\neq 0$, this reads as a ratio of products of determinants\med

\centerline{$\ds \cross_q(x)= {({a_i\over b_i}-{a_k\over b_k})\cdot({a_j\over b_j}-{a_\ell\over b_\ell})\over ({a_i\over b_i}-{a_\ell\over b_\ell})\cdot({a_j\over b_j}-{a_k\over b_k})}=
{({a_ib_k}-{a_kb_i})\cdot({a_jb_\ell}-{a_\ell b_j})\over ({a_ib_\ell}-{a_\ell b_i})\cdot({a_jb_k}-{a_kb_j})}$,}\med

and symmetrically for the second chart. Here, we have set $0=(0:1)$, $1=(1:1)$, $\infty=(1:0)$.\med

One has the following rules: If the involved entries $x_i, x_j,x_k,x_\ell$ are pairwise distinct, the cross-ratio $\cross_q(x)$ takes a value in $\P^1\sm\{0,1,\infty\}$, and if two entries or two pairs of two entries are equal, it takes a value in $\{0,1,\infty\}$. The special values $0,1,\infty$ thus govern the equality of entries. If three entries are equal, the cross-ratio $\cross_q(x)$ is not defined. The precise distribution of the values of $\cross_q(x)$ is as follows.\med

\hs .5cm If the four entries $x_i$, $x_j$, $x_k$ and $x_\ell$ are pairwise distinct, the cross-ratio is different from $0$, $1$, $\infty$;\med

\hs .5cm If the first two entries $x_i$, $x_j$ and/or the last two entries $x_k$, $x_\ell$ are equal, the cross-ratio is $1$;\med

\hs .5cm If the first and third entry $x_i$, $x_k$ and/or the second and fourth entry $x_j$, $x_\ell$ are equal, the cross-ratio is $0$;\med

\hs .5cm If the first and last entry $x_i$, $x_\ell$ and/or the second and third entry $x_j$, $x_k$ are equal, the cross-ratio is $\infty$;\med

\hs .5cm If three of the entries $x_i$, $x_j$, $x_k$ and $x_\ell$ are equal, the cross-ratio is not defined.\med

Two $n$-gons $x$ and $y$ in $(\P^1)^n$ may define different $\PGL_2$-orbits even though all their cross-ratios are equal: The $4$-gons $(0,1,\infty,0)$, $(0,1,1,\infty)$ have the same cross-ratio $\cross_{(1234)}(x)=\infty$ but are not $\PGL_2$-equivalent.\med

According to the context, cross-ratios will be considered as formal cross-ratios $[q]=[ijk\ell]$, i.e., elements of $\kk(\xi_1,...,\xi_n)$, or as a cross-ratio functions $\cross_q$ on open subsets of $(\P^1)^n$, with evaluations $\cross_q(x)$ at $x\in (\P^1)^n$ whenever these are defined. \med

{\it Example.} \label{example:sevengon} For $x=(0,1,1,\infty,\infty, \infty, a)\in (\P^1)^7$, with $a\neq 0,1,\infty$, the quadruples $q=(1247)$ and $(1237)$ give cross-ratios $a$ and $\infty$, whereas for $q=(1456)$ the cross-ratio is not defined. \med




As a matter of interest we state\med

{\bf Lemma.} \label{lemma:crossratiosgenerate} (Cross-ratios) (a) {\it The formal cross-ratios $[ijk\ell]$, for $(ijk\ell)\in N^{\ul 4}$ generate the subfield of $K(\xi_1,...,\xi_n)$ of rational $\PGL_2$-invariants.}

(b) {\it The algebraic relations between the formal cross-ratios are generated by the obvious ones given by the permutation of the indices and the triple product formula as indicated above.}%
%
\tex
{{\baselineskip 12pt \footnote{$\!\!{}^3$}{\smallTimes These relations can be interpreted as the {\it Pl\"ucker relations} between products of determinants of matrices.}}}
{\footnote{These relations can be interpreted as the {\it Pl\"ucker relations} between products of determinants of matrices.}}
\med


{\it Proof.} (a) Let $\varphi\in K(\xi_1,...,\xi_n)$ be a $\PGL_2$-invariant, and let $x=(x_1,...,x_n)$ be a generic $n$-gon in $(\P^1)^n\sm \Delta_n$. Then $\varphi$ is constant on the orbit $[x]$ of $x$. As $\PGL_2$ acts $3$-transitively, there is an $n$-gon in the orbit of the form $y=(0,1,\infty, y_4,...,y_n)$. The entries $x_i$ of $x$ are rational functions in $y_4,...,y_n$. Each $y_i$ equals the cross-ratio $\cross_{(i213)}(y)$, which, in turn, equals $\cross_{(i213)}(x)$ since $y\sim x$. This shows that $\varphi(x)=\varphi(y)$ is a rational function in cross-ratios evaluated in $x_1,...,x_n$.\qed\gb

(b) Rewrite the triple product formula\med

\cl{$[ijk\ell]\cdot [ij\ell m]\cdot [ijmk]=1$.}\med

as

\cl{$[ijk\ell]=[ijkm]\cdot [ijm\ell]$.}\med


Let now $R(z)\in K[z_{(ijk\ell)},\, (ijk\ell)\in N^{\ul 4}]$ be an arbitrary algebraic relation between formal cross-ratios, for new variables $z_{(ijk\ell)}$. By the displayed formula, we may replace all $z_{(ijk\ell)}$ with $1<i,j,k,\ell\leq n$ by variables whose indexing quadruple has one entry $i=1$. Up to a permutation of the entries of the quadruples, we may assume that $R$ only depends on variables $z_{(1jk\ell)}$, for $2\leq j,k,\ell\leq n$. Similarly, we may then replace all variables $z_{(1jk\ell)}$ with $2<j,k,\ell\leq n$ by variables whose indexing quadruple has entries $i=1$ and $j=2$. Repeating the trick, we reduce to variables $z_{(123\ell)}$ with $4\leq \ell\leq n$. The respective cross-ratios are rational functions in $\xi_1, \xi_2,\xi_3$ and $\xi_\ell$. As $\xi_\ell$ appears, for each $\ell$, in exactly one such cross-ratio, there is only the trivial relation $R=0$ left. This proves (b).\qed


\med\gb


{\Bf 8. Strings of $n$-gons}\label{section:strings}\med

An {\it $n$-gon} $x$ is a point $x=(x_1,...,x_n)\in (\P^1)^n$, where $\P^1=\P^1_\kk$ is the projective line over an arbitrary base field $\kk$, seen as $\P^1=K\cup\{\infty\}$. We call $x_i\in\P^1$ the {\it entries} of $x$. A {\it string} is a vector \med

\cl{$\xx=(x^t)_{t\in {N\choose 3}}\in (\P^1)^{n{n\choose 3}}$}\med

of $n$-gons $x^t=(x^t_1,...,x^t_n)$ in $(\P^1)^n$, where the superscript $t$ varies over all triples $t=(ijk)\in {N\choose 3}$. We will be mostly interested in strings whose $n$-gons have at least three different entries. Further on, we will put soon restrictions on the strings to be considered by requiring that their $n$-gons have equal cross-ratios $\cross_q(x^s)=\cross_q(x^t)$ for all quadruples $q$ in $N$ (whenever the cross-ratios on both sides of the equation are defined). \med \gb

Define a smooth closed subvariety \med

\cl{$\TT_n\subseteq(\P^1)^{n{n\choose 3}}$}\med

as the set of strings $\xx=(x^t)_{t\in {N\choose 3}}$ for which\med

\cl{$x^t_i=0$, $x^t_j=1$, $x^t_k=\infty$}\med

holds for all triples $t=(ijk)$ in $N\choose 3$. As such, every $n$-gon $x^t$ of a string $\xx$ in $\TT_n$ has  ab initio at least three different entries. Clearly, $\TT_n$ is a projective variety isomorphic to $(\P^1)^{(n-3){n\choose 3}}$. \med

For $t=(ijk)\in{N\choose 3}$, the further entries $x^t_\ell$ of $x^t$, for $\ell\neq i,j,k$, can be expressed as a rational function in the cross-ratio $\cross_{(ijk\ell)}(x^t)$ of $x^t$. More precisely, \med

\cl{$\ds \cross_{(ijk\ell)}(x^t)={(x^t_i-x^t_k)(x^t_j-x^t_\ell)\over 
(x^t_i-x^t_\ell)(x^t_j-x^t_k)}={(0-\infty)(1-x^t_\ell)\over (0-x^t_\ell)(1-\infty)}
={x^t_\ell -1\over x^t_\ell}$,}\med 

respectively, \med

\cl{$\ds x^t_\ell ={1\over 1-\cross_{ijk\ell}(x^t)}= \cross_{(ik\ell j)}(x^t)= \cross_{(\ell jik)}(x^t)$.}\med

Hence, if $x^t$ and $x^r$ have the same cross-ratios $\cross_{(ijk\ell)}(x^t)=\cross_{(ijk\ell)}(x^r)$, for triples $t=(ijk)$ and $r=(ij\ell)$, then $\cross_{(ij\ell k)}(x^r)=\cross_{(ijk\ell)}(x^r)\inv$ implies that $x^t_\ell+x^r_k=1$ holds, a formula to be used in Section 20.\med

Here is a typical string in $\TT_4$,\med

\cl{$\xx=(x^{123}, x^{124}, x^{134}, x^{234})=((\ul 0,\ul 1,\ul \infty,a_4), (\ul 0,\ul 1,a_3,\ul \infty), (\ul 0,a_2,\ul 1,\ul \infty),(a_1,\ul 0,\ul 1,\ul \infty))$,}\med

where the prescribed values are underlined and $a_1,...,a_4$ can take arbitrary values in $\P^1$. \med\med\gb


{\Bf 9. The varieties $\XX_n$ and $\YY_n$} \label{section:varietiesxxnyyx}\med

If $(\P^1)^n\sm \Delta_n\subseteq (\P^1)^n$ denotes the open subvariety of $n$-gons with pairwise distinct entries, for $\Delta_n$ the big diagonal, we get a natural morphism\med

\hs5cm{$\tau_n: (\P^1)^n\sm \Delta_n\hs.3cm \map \hs.5cm\TT_n$,}\med

\hs 5cm{$ x=(x_1,...,x_n)\hs.35cm \map \hs.5cm\xx=(x^t)_{t\in {N\choose 3}}$,}\med

sending an $n$-gon $x$ to the string $\xx=(x^t)_{t\in {N\choose 3}}\in\TT_n$ such that for every triple $t=(ijk)$, the $n$-gon $x^t$ is the unique element in the orbit of $x$ with $x^t_i=0$, $x^t_j=1$, $x^t_k=\infty$. In this way all $n$-gons $x^t$ of $\xx$ are different representatives of the same orbit, namely, the one of $x$. The string $\xx$ will be called the {\it symmetrization} of $x$ with respect to $\PGL_2$. For every quadruple $q$ in $N$, the cross-ratios $\cross_q(x^t)$, $t\in {N\choose 3}$, are equal whenever they are defined. This is the case, for instance, when $q$ involves the entries $i,j,k$ of $t$. This observation justifies to denote the common value of the cross-ratios by $\cross_q(\xx)$.\med


{\it Example.} \label{example:fourgon} The $4$-gon $x=(0,1,\infty,a)$ of $(\P^1)^4$ with $a\in \P^1\sm\{0,1,\infty\}$ is sent by $\sigma_4$ to the string\med

\cl{$\xx=((\ul 0,\ul 1,\ul \infty,a), (\ul 0,\ul 1,1-a,\ul \infty),(\ul 0,{1\over 1-a},\ul 1,\ul \infty), ({1\over a},\ul 0,\ul 1,\ul \infty))$,}\med

of $4$-gons with the same cross-ratio $\cross_{(1234)}(x^t)=1-{1\over a}$ for all $t$. If $a=0,1,\infty$, respectively, one obtains $\cross_{(1234)}(x^t)=\infty, 0,1$, respectively.\med


Let $\UU_n=((\P^1)^n\sm \Delta_n)/\PGL_2$ be the orbit space of {\it generic} $n$-gons under the action of $\PGL_2$, i.e., those $n$-gons with pairwise distinct entries. We get an embedding\med

\cl{$\Sigma_n:\UU_n\hookrightarrow \TT_n\subseteq (\P^1)^{n{n\choose 3}}$,}\med

\cl{$[x]\map\xx=\Sigma_n(x)$,}\med

sending the $\PGL_2$-orbit $[x]$ of $x\in (\P^1)^n\sm \Delta_n$ to the string $\xx=(x^t)_{t\in {N\choose 3}}$ as above. We identify $\UU_n$ with its image $\Sigma_n(\UU_n)$, say, with its symmetrization. Then define the closed subvariety $\XX_n=\XX_n(\P^1)$ as the Zariski closure of (the symmetrization of) $\UU_n$ in $\TT_n$ (or, equivalently, in $(\P^1)^{n{n\choose 3}}$),\med

\cl{$\XX_n= \ol{\Sigma_n(\UU_n)}\subseteq (\P^1)^{n{n\choose 3}}$.}\med

We often just write $\XX_n=\ol{\UU_n}$. This variety will be shown to be the appropriate compactification of $\UU_n$. It will eventually turn out to be isomorphic to the moduli space $\ol{\MM_{0,n}}$ of $n$-pointed stable curves of genus zero.\med

Its set of labels is $N=\{1,...,n\}$. The strings of the dense open subset $\UU_n$ correspond to {\it irreducible} $n$-pointed stable curves, that is, to $n$ pairwise distinct points on $\P^1$. The boundary divisor $\BB_n=\XX_n\sm \UU_n$ will consist of strings $\xx$ associated to $n$-gons $x\in (\P^1)^n$ where certain entries have come together and become equal. The embedding $\Sigma_n$ is a smart trick to control this coalescing systematically and to prove the required properties of $\XX_n$.%
\tex
{{\baselineskip 12pt \footnote{$\!\!{}^4$}{\smallTimes A similar construction as the one for $\XX_n$ via strings appears in {\cite{\[Br\]}}, Section 2.1, p.~381, and Section 2.8, p.~402: The orbit $[x]$ of a generic $n$-gon $x\in (\P^1)^n$ is mapped to the vector in $(\P^1)^{n\choose 4}$ of all its cross-ratios $\cross_q(x)$, $q$ a quadruple in $N^{\ul 4}$, thus getting as image of $\UU_n$ a locally closed subvariety ${\frak M}_{0,n}$ of $(\P^1)^{n\choose 4}$. Its Zariski-closure $\ol{\frak M}_{0,n}$ will then be the requested compactification. As the two embeddings $\Sigma_n:\UU_n\map (\P^1)^{n{n\choose 3}}$ and $\UU_n\map {\frak M}_{0,n}\subset (\P^1)^{n\choose 4}$ are compatible to each other, one has $\XX_n\isom \ol{\frak M}_{0,n}$. In \cite{\[GHP\]}, p.~135, the authors associate to $n$-pointed stable curves a vector of cross-ratios: The $n$ marked points on the curve are projected, for each triple $t$, to the {\smallit median} component of the curve selected by $t$, and then the cross-ratio is taken on this component (which is isomorphic to $\P^1$).}}}
{{\footnote{A similar construction as the one for $\XX_n$ via strings appears in \cite{\[Br\]}, Section 2.1, p.~381, and Section 2.8, p.~402: The orbit $[x]$ of a generic $n$-gon $x\in (\P^1)^n$ is mapped to the vector in $(\P^1)^{n\choose 4}$ of all its cross-ratios $\cross_q(x)$, $q$ a quadruple in $N^{\ul 4}$, thus getting as image of $\UU_n$ a locally closed subvariety ${\frak M}_{0,n}$ of $(\P^1)^{n\choose 4}$. Its Zariski-closure $\ol{\frak M}_{0,n}$ will then be the requested compactification. As the two embeddings $\Sigma_n:\UU_n\map (\P^1)^{n{n\choose 3}}$ and $\UU_n\map {\frak M}_{0,n}\subset (\P^1)^{n\choose 4}$ are compatible to each other, one has $\XX_n\isom \ol{\frak M}_{0,n}$. In \cite{\[GHP\]}, p.~135, the authors associate to $n$-pointed stable curves a vector of cross-ratios: The $n$ marked points on the curve are projected, for each triple $t$, to the {\smallit median} component of the curve selected by $t$, and then the cross-ratio is taken on this component (which is isomorphic to $\P^1$).}}}
%
\med

To define $\XX_n$ by equations, let $\YY_n\subseteq \TT_n\subseteq (\P^1)^{n{n\choose 3}}$ be the closed subvariety defined as the zero-set \med

\hs 1cm $\YY_n=\{\xx=(x^t)_t\in \TT_n,\, \cross_q(x^s)=\cross_q(x^t)$ for all $s,t\in{N\choose 3}$ and all $q\in N^{\ul 4}\}$.\med

The equations $\cross_q(x^s)=\cross_q(x^t)$ are understood here as the polynomial equations \med

\cl{$(x^s_i-x^s_k)(x^s_j-x^s_\ell) (x^t_i-x^t_\ell)(x^t_j-x^t_k) = 
(x^s_i-x^s_\ell)(x^s_j-x^s_k) (x^t_i-x^t_k)(x^t_j-x^t_\ell)$}\med

obtained from the equality of the cross-ratios after clearing denominators on both sides. As such, there occurs no problem when the equation is evaluated on $n$-gons $x^s$ or $x^t$ with three equal entries (for which the cross-ratio would not be defined) since it becomes the trivial equation $0=0$.  From what we have seen before, $\Sigma(\UU_n)$ is contained in $\YY_n$.  As $\YY_n$ is closed, we get $\XX_n\subseteq\YY_n$. It turns out that this inclusion is in fact an equality, $\XX_n=\YY_n$. This allows us to describe $\XX_n$ by equations. The equality of the two varieties will be proven in the course of the article.\med

\ignore

\comment{[\ooo Internal comment on Brown's paper in Annales ENS 2009 [Br]: 
 \vs.1cm\baselineskip 12pt {\smallTimes
Brown denotes the label set $S=\{1,...,n\}$ instead of our $N$. On page 381, he introduces a space ${\frak M}_{0,S}$ as the space of $\PGL_2$-orbits $[x]$ of generic $n$-gons $x=(x_1,...,x_n)$ (there is a slight inaccuracy in line 2 of Section 2.1 since his "curves of genus $0$" must be a $\P^1$, say, smooth and irreducible, and not just any rational curve). \vs.1cm
On page 381, line -2, he defines {\smallit simplicial coordinates} on ${\frak M}_{0,S}$ by sending an orbit $[x]$ to its distinguished representative $z:=x^t=(0,1,\infty, z_4,...,z_n)$ with triple $t=(123)$. The coordinates are then $z_4,...,z_n$. This defines an embedding of ${\frak M}_{0,S}$ into $(\P^1)^{n-3}$. {\smallit Cubical coordinates} $x_1,...,x_{n-3}$ are defined on page 382, and {\smallit dihedral coordinates} $u_{ij}$ are defined via selected cross-ratios on page 382, giving an embedding of ${\frak M}_{0,S}$ into affine space $\A^{n(n-3)/2}$, see line -1. The constructions are somewhat ad hoc and seem to lack symmetry.\vs.1cm
A variety ${\frak M}^\delta_{0,S}$ (called {\smallit dihedral extension}) is introduced in 2.11 on page 385, via a bit mysterious equations (2.10). Then ${\frak M}_{0,S}$ is embedded into ${\frak M}^\delta_{0,S}$ in Lemma 2.5 on page 385. See also Thm.~2.15, page 392, and Thm.~2.25, page 398, respectively, Cor.~2.32, page 403, for the smoothness of these spaces and the description of their boundary.\vs.1cm
The projection map $\ol {\frak M}_{0,S}\map \ol {\frak M}_{0,T}$ for $T\subset S$ is considered in Section 2.4, page 388.\vs.1cm
The compactification $\ol{\frak M}_{0,S}$ is introduced in Section 2.8, page 402, as the Zariski-closure of ${\frak M}_{0,S}$ via the embedding of ${\frak M}_{0,S}$ into $(\P^1)^{n\choose 4}$, taking the vector of all cross-ratios of an orbit $[x]$ of an $n$-gon $x\in (\P)^1$, see formula 2.41 on page 402.\vs.1cm
In Section 2.9, page 405, Brown talks about the isomorphism between his compactification $\ol{\frak M}_{0,S}$ and the Deligne-Mumford-Knudsen compactification $\ol \MM_{0,n}$. He claims that, by the universal property of the latter, there is a morphism $\ol{\frak M}_{0,S}\map \ol \MM_{0,n}$. But for having this, one would have to show that $\ol {\frak M}_{0,S}\map \ol {\frak M}_{0,T}$ is flat and that its fibers are $n$-pointed stable curves (with respect to some well chosen sections, which even seem not to have been defined in \cite{\[Br\]}). It also seems that the fibers of the projection map $\ol {\frak M}_{0,S}\map \ol {\frak M}_{0,T}$ are not discussed or properly described in \cite{\[Br\]}, only for the map ${\frak M}_{0,S}\map {\frak M}_{0,T}$ it is mentioned that the fibers are one-dimensional, see the last but one displayed formula on page 389. He then alludes to a universal property of $\ol{\frak M}_{0,S}$, but he does not say what the equivalence classes would be nor how a respective family of equivalence classes should be understood. Note here that the orbits $[x]$ of $n$-gons are only defined for generic $n$-gons, the elements of the boundary $\ol{\frak M}_{0,S}\sm {\frak M}_{0,S}$ are just vectors (of cross-ratios) and not equivalence classes. Eventually, he says "One could probably check that this is an isomorphism", but he then adds  "we have not done this".\ooo]}\med}

\recognize


\med\gb

{\Bf 10. Limits of orbits} \label{section:limitsorbits}\med

We pause for a moment to illustrate in an informal manner the preceding completion process -- it is given by passing to the Zariski-closure $\XX_n$ of $\Sigma(\UU_n)$ -- in terms of limits of a concrete $n$-gon $x=(x_1,...,x_n)$ together with its $\PGL_2$-orbit $[x]\in\UU_n$ as some entries of $x$ approach each other. Let us take $x=(0,1,\infty,a,b)\in(\P^1)^5$ with $a\neq b$ in $\P^1\sm\{0,1,\infty\}$. We represent this $5$-gon as a star-like graph with five leaves, where the central vertex represents the orbit of $x$ and the five different edges connecting it to the leaves express the fact that all five entries of $x$ are pairwise distinct (Fig.~10). 


\cl{\fivelabelstreegeneric}

\cl{{\it Figure 10.} The generic tree $T_*$ with one vertex and five leaves.}\med\med
 

In the limit, as $a$ and $b$ approach each other and coalesce, we obtain the $5$-gon $(0,1,\infty,a,a)$ with two equal entries. Take now the symmetrization of $x$: this will be a string $\xx$ in $(\P^1)^{5{5\choose 3}}=(\P^1)^{50}$, that is, consists of ten $5$-gons $x^t$, all having the same orbit as $x$. We have $x^{(123)}=x$, and, for instance,  $x^{(345)}= (a',b',0,1,\infty)$, with $a'$ and $b'$ rational functions in $a$ and $b$. As $a$ and $b$ come together, it follows that both $a'$ and $b'$ tend to $0$. So the limit of $x^{(345)}$ is $(0,0,0,1,\infty)$, now with three equal entries. Observe that the two limiting $5$-gons $(0,1,\infty,a,a)$ and $(0,0,0,1,\infty)$ define different orbits, but have, by continuity, the same cross-ratios. One could informally say that these two orbits are both {\it limits} of the orbit of $x$, that is, they lie in the boundary of the union of orbits $[x]$ as $x$ varies.\med

The symmetrization of this $x$ consists of ten $5$-gons defining precisely the two orbits. We represent this by the graph with two inner vertices (corresponding to the two orbits) and two, respectively three leaves attached to each of them by an edge (corresponding to the incidence of entries in the limit $5$-gons), see Fig.~100, left.\med
\med


Here are the exact formulas: Let $\xx=(x^t)_{t\in {\{1,2,3,4,5\}\choose 3}}\in \YY_5$ be a string of $5$-gons with equal cross-ratios. We order the triples in $\{1,2,3,4,5\}$ as follows\med

\centerline{$(123), (124), (125), (134), (135), (145), (234), (235), (245), (345)$.}\med
 
The ten $5$-gons of a generic string $\xx\in\YY_5$, say, with pairwise different entries of its $5$-gons, are listed below. We underline in $x^t$ the entries with indices $i,j,k$ if $t=(ijk)$. The constants $a\neq b$ may assume any value in $\P^1\sm\{0,1,\infty\}$. \med

{\parindent 3cm 
 
$x^{(123)}=(\ul 0,\ul 1,\ul \infty,a,b)$, \med
$x^{(124)}=(\ul 0,\ul 1,1-a,\ul \infty,{1-{1\over a}\over {1\over b}-{1\over a}})$, \med
$x^{(125)}=(\ul 0,\ul 1,1-b,{1-{1\over b}\over {1\over a}-{1\over b}},\ul \infty)$, \med
$x^{(134)}=(\ul 0,{1\over 1-a},\ul 1,\ul \infty,{1\over 1-{a\over b}})$,  \med
$x^{(135)}=(\ul 0,{1\over 1-b},\ul 1,{1\over 1-{b\over a}},\ul \infty)$, \med
$x^{(145)}=(\ul 0,{{1\over a}-{1\over b}\over 1-{1\over b}}, 1-{b\over a},\ul 1,\ul \infty)$\med
$x^{(234)}=({1\over a},\ul 0,\ul 1,\ul \infty,{1-{1\over b}\over 1-{a\over b}})$,\med
$x^{(235)}=({1\over b},\ul 0,\ul 1,{1-{1\over a}\over 1-{b\over a}},\ul \infty)$, \med
$x^{(245)}=({{1\over b}-{1\over a}\over 1-{1\over a}},\ul 0,{1-{b\over a}\over 1-{1\over a}},\ul 1,\ul \infty)$, \med
$x^{(345)}=(1-{a\over b},{1-{a\over b}\over 1-{1\over b}},\ul 0,\ul 1,\ul \infty)$.\med

}
Up to permutation of the entries of the quadruples, the five cross-ratios are\med

{\parindent 3cm 

$\cross_{(1234)}(\xx)= 1-{1\over a}$, \med
$\cross_{(1235)}(\xx)= 1-{1\over b}$,  \med
$\cross_{(1245)}(\xx)={1-{1\over b}\over 1-{1\over a}} $,  \med
$\cross_{(1345)}(\xx)= {a\over b}$,  \med
$\cross_{(2345)}(\xx)= {1-a\over 1-b}$.\med

}

We may now take limits of such strings $\xx$ as $a$ and/or $b$ tend to one of the values $0,1,\infty$ or become equal $a=b$. As a matter of illustration, let $a$ and $b$ both go to $\infty$. For $x^{(123)}=(\ul 0,\ul 1,\ul \infty,a,b)$, the limit will be the single $5$-gon $\ol x^{(123)}=(\ul 0,\ul 1,\ul \infty,\infty,\infty)$, regardless how fast $a$ and $b$ tend to $\infty$ (one may take $b=a-1$ to have distinct values of $a$ and $b$). Similarly, the limit $\ol x^{(124)}$ of $x^{(124)}$ also equals $(\ul 0,\ul 1,\infty,\ul \infty,\infty)$ as $a,b\map \infty$. For $x^{(134)}=(\ul 0,{1\over 1-a},\ul 1,\ul \infty,{1\over 1-{a\over b}})$, the situation is quite different: now the limit depends on the ratio of $a\over b$, and we get accordingly various limits $\ol x^{(134)}=(\ul 0,0,\ul 1,\ul \infty,c)$, for $c$ varying arbitrarily in $\P^1$. Indeed, to get the value $c=1$, take $a\map \infty$ and $b=a^2\map\infty$, and to get a value $c\neq 1$ take $a\map \infty$ and $b={c\over c-1}a$. Then $b\map \infty$ if $c\neq 0$; otherwise, for $c=0$, $b$ will be constant equal to $0$. In all cases we get $\ol x^{(134)}=(\ul 0,0,\ul 1,\ul \infty,c)$ as the limit of $x^{(134)}=(\ul 0,{1\over 1-a},\ul 1,\ul \infty,{1\over 1-{a\over b}})$.\med\gb

We summarize the situation for these first three $5$-gons $x^{(123)}$, $x^{(124)}$, and $x^{(134)}$ of $\xx$ as $a$ goes to $\infty$, letting for instance $b=a+1$ for $x^{(123)}$,  respectively, $b={1\over 2}a$ for $x^{(124)}$, and $b=a^2$ (case $c=1$), respectively, $b={c\over c-1}a$ (case $c\neq 1$) for $x^{(134)}$:\med

{\parindent 1.5cm 
 
$x^{(123)}=(\ul 0,\ul 1,\ul \infty,a,b)=(\ul 0,\ul 1,\ul \infty,a,a+1)\map \ol x^{(123)}=(\ul 0,\ul 1,\ul \infty,\infty,\infty)$, \med

$x^{(124)}=(\ul 0,\ul 1,1-a,\ul \infty,{1-{1\over a}\over {1\over b}-{1\over a}})=\ul 0,\ul 1,1-a,\ul \infty,a-1)\map \ol x^{(124)}=(\ul 0,\ul 1,\ul \infty,\infty,\infty)$, \med

$x^{(134)}=(\ul 0,{1\over 1-a},\ul 1,\ul \infty,{1\over 1-{a\over b}})=(\ul 0,{1\over 1-a},\ul 1,\ul \infty,{1\over 1-{1\over a}})\map \ol x^{(134)}=(\ul 0,0,\ul 1,\ul \infty,1)$,  \med

$x^{(134)}=(\ul 0,{1\over 1-a},\ul 1,\ul \infty,{1\over 1-{a\over b}})=(\ul 0,{1\over 1-a},\ul 1,\ul \infty,c)\map \ol x^{(134)}=(\ul 0,0,\ul 1,\ul \infty,c)$.  \med

}
Note here that $\ol x^{(123)}=\ol x^{(124)}$ and $\ol x^{(134)}$ still have equal cross-ratios (whenver they are defined), but their orbits are different. Analogous considerations can be applied for the limits of the remaining $5$-gons of $\xx$. \med

The different behaviour of the limits of $5$-gons reflects the fact that the limit of a string $\xx$ as above varies in a one dimensional subvariety of $\YY_5$, even though this is not apparent by taking solely the (unique) limit of the first $5$-gon $x^{(123)}$. Put differently, the boundary $\BB_5:=\YY_5\sm\UU_5$ is a projective curve. \med


Let us now describe all possible $5$-gons $x^{(123)}$ of strings $\xx$ in $\YY_5$,  up to a permutation of the components and of the values $0$, $1$, $\infty$.\med

{\parindent 3cm 

$(\ul 0,\ul 1,\ul \infty,a,b)$, $a,b\neq 0,1,\infty$, $a\neq b$ (all entries distinct, generic case),\med

$(\ul 0,\ul 1,\ul \infty,a,a)$, $a\neq 0,1,\infty$ (two equal entries, both outside $\{0,1,\infty\}$),\med

$(\ul 0,\ul 1,\ul \infty,a,\infty)$, $a\neq 0,1,\infty$ (two equal entries from $\{0,1,\infty\}$),\med

$(\ul 0,\ul 1,\ul \infty,1,\infty)$,  (two distinct pairs of equal entries in $\{0,1,\infty\}$),\med

$(\ul 0,\ul 1,\ul \infty,\infty,\infty)$, (three equal entries from$\{0,1,\infty\}$),\med

}

The five cross-ratios $\cross_{(1234)}$, $\cross_{(1235)}$, $\cross_{(1245)}$, $\cross_{(1345)}$, $\cross_{(2345)}$ are in each case \med\gb

{\parindent 3cm 

$(\ul 0,\ul 1,\ul \infty,a,b)$: $\cross_q=1-{1\over a}$, $1-{1\over b}$, ${1-{1\over b}\over 1-{1\over a}}$, ${a\over b}$, ${a\over b}$,\med

$(\ul 0,\ul 1,\ul \infty,a,a)$: $\cross_q=1-{1\over a}$, $1-{1\over a}$, $1$, $1$, $1$,\med

$(\ul 0,\ul 1,\ul \infty,a,\infty)$: $\cross_q=1-{1\over a}$, $1$, ${1\over 1-{1\over a}}$, $0$, $0$,\med

$(\ul 0,\ul 1,\ul \infty,1,\infty)$: $\cross_q=0$, $1$, $\infty$, $0$, $0$,\med

$(\ul 0,\ul 1,\ul \infty,\infty,\infty)$: $\cross_q= 1$, $1$, $1$, $0\over 0$, $0\over 0$.\med

}

To familiarize with these computations, let us determine all strings $\yy\in\YY_5$ whose $5$-gon $y^{(123)}$ is $(\ul 0,\ul 1,\ul \infty,\infty,\infty)$. The cross-ratios $\cross_{(1234)}$, $\cross_{(1235)}$, and $\cross_{(1245)}$ are defined and equal to $1$, the remaining two, $\cross_{(1345)}$ and $\cross_{(2345)}$, are not defined. This signifies that the first two or the last two entries of the $4$-gon selected by the quadruple $q=(ijk\ell)$ from $y^t$ must be equal for $q=(1234),(1235),(1245)$, and no conditions follow from  the other two quadruples $q=(1345), (2345)$. We get, for any $c\in \P^1$, the string $\yy=(y^t)_{t\in {\{1,...,5\}\choose 3}}$ given by\med

{\parindent 3cm 
 
$y^{(123)}=(\ul 0,\ul 1,\ul \infty,\infty,\infty)$,\med

$y^{(124)}=(\ul 0,\ul 1,\infty,\ul \infty,\infty)$,\med

$y^{(125)}=(\ul 0,\ul 1,\infty,\infty,\ul \infty)$,\med

$y^{(134)}=(\ul 0,0,\ul 1,\ul \infty,c)$,\med

$y^{(135)}=(\ul 0,0,\ul 1,1-c,\ul \infty)$, \med

$y^{(145)}=(\ul 0,0, {1\over 1-c},\ul 1,\ul \infty)$, \med

$y^{(234)}=(0,\ul 0,\ul 1,\ul\infty,c)$,\med

$y^{(235)}=(0,\ul 0,\ul 1,1-c,\ul \infty)$,\med

$y^{(245)}=(0,\ul 0,{1\over 1-c},\ul 1,\ul \infty)$, \med

$y^{(345)}=({1\over c},{1\over c},\ul 0,\ul 1,\ul \infty)$.\med

}

If $c\neq 0,1,\infty$, one obtains precisely two orbits, whereas for $c\in\{0,1,\infty\}$, three orbits occur. The situation is represented by the graphs depicted in Fig.~100 (see Section 13 for the details of this presentation).\med


\cl{\fivelabelstreesspecial}

\cl{{\it Figure 100.} The graphs associated to the string $\yy$, with two orbits (left) and three orbits (right).}\med
 

The vertex $u$ corresponds on both sides to the orbit $[y^{(123)}]$, $w$ corresponds to $[y^{(145)}]$ (for $c=0$), whereas $v$ corresponds to $[y^{(134)}]$ on the left hand side, but to $[y^{(234)}]$ on the right hand side (for $c=\infty$). 
\med\med\gb


{\Bf 11. The projection map $\pi_a: \XX_{n+1}\map \XX_n$} \label{section:projectionmap}\med

Denote by $N$ and $N^{+1}$ the sets of labels of strings $\xx=(x^t)\in (\P^1)^{n{n\choose 3}}$, respectively, of strings $\yy\in (\P^1)^{(n+1){n+1\choose 3}}$. We may assume that $N^{+1}=N\cup\{a\}$ where $a$ is the label of $N^{+1}$ which is not in $N$. In this way, triples $t$ in ${N\choose 3}$ are also triples in ${N^{+1}\choose 3}$, namely those which do not involve the label $a$. There is then a natural \projection map \med

\cl{$\Pi_a: (\P^1)^{(n+1){n+1\choose 3}}\map(\P^1)^{n{n\choose 3}}$}\med
 
given by ``forgetting'' entries involving the label $a$. More explicitly, if $\yy=(y^s)_{s\in{N^{+1}\choose 3}}$ is a string in $(\P^1)^{(n+1){n+1\choose 3}}$, define its image $\pi_a(\yy)=\xx$ in $(\P^1)^{n{n\choose 3}}$ as the string $\xx=(x^t)_{t\in{N\choose 3}}$ deleting first all $n$-gons $y^s$ of $\yy$ with $s\in {N^{+1}\choose 3}\sm {N\choose 3}$ and taking then, for $t\in{N\choose 3}$, the $n$-gon $x^t\in (\P^1)^n$ obtained from the $(n+1)$-gon $y^t\in(\P^1)^{n+1}$ by deleting from it the entry $y^t_a$ with index $a$.\med
 
The map $\Pi_a$ can be restricted to the subvariety $\TT_{n+1}$ of  $(\P^1)^{(n+1){n+1\choose 3}}$ of strings $\yy=(y^s)_{s\in {N^{+1}\choose 3}}$ with $y^s_i=0$, $y^s_j=1$, $y^s_k=\infty$ for $s=(ijk)$, and then induces a projection\med

\cl{$\Pi_a:\TT_{n+1}\map \TT_n$.}\med

We may restrict further to both $\XX_{n+1}$ and $\YY_{n+1}$ (we do not know yet that they are equal, only that $\XX_{n+1}\subseteq \YY_{n+1}$) and get well-defined projections (which, a posteriori, will be identical)\med

\cl{$\pi_a:\XX_{n+1}\map \XX_n$\hs.3cm and \hs.3cm $\pi_a:\YY_{n+1}\map \YY_n$.}\med

Indeed, observe that $\Pi_a$ sends $\UU_{n+1}$ into $\UU_n$ since the action of $\PGL_2$ on $(\P^1)^n$ is the restriction of the action of $\PGL_2$ on $(\P^1)^{n+1}$. Therefore the Zariski-closure $\XX_{n+1}$ of $\UU_{n+1}$ is mapped into $\XX_n$. Further, the for\-mal cross-ratios $[ijk\ell]$ in $K(\xi_N)$ are equal to the formal cross-ratios in $K(\xi_{N^{+1}})$ that do not involve the label $a$. Hence also $\YY_{n+1}$ is mapped into $\YY_n$.\med

It is not hard to see that if $\yy\in\XX_{n+1}$ maps to $\xx\in\XX_n$, the phylogenetic tree $\Gamma_\xx$ is obtained from $\Gamma_\yy$ by clipping off the leaf with label $a$ and contracting edges if required, see Section 12 for the details.\med

The map $\pi_a$ is a surjective proper morphism of varieties. It will be proven in part III, Section 20, that the fibers $\FF_\xx=\pi_a\inv(\xx)$ of $\pi_a$ have constant dimension $1$. As the source and target spaces are non-singular varieties (proven in part III, Section 17), one may apply the flatness criterion from {\cite{\[Mat\]}}, Thm.~23.1, to conclude that $\pi_a$ is in fact a flat morphism. In part III, Section 20, it is shown that the fibers are stable rational curves (unions of projective lines $\P^1$ meeting transversally and at most pairwise), and, in Section 21, $n$ sections $\sigma_p$ of $\pi_a$ will be constructed. The fibers thus become $n$-pointed stable curves, that is, such curves arise directly from the consideration of strings of $n$-gons. And it will then be no surprise that the augmented dual graph $\Gamma_{\FF_\xx}$ of $(\FF_\xx, \sigma_1(\xx),...,\sigma_n(\xx))$ equals the phylogenetic tree $\Gamma_\xx$ of the string $\xx$. In this way, the cycle of constructions closes up.\med

Moreover, the morphism $\pi_a:\XX_{n+1}\map \XX_n$ turns out to be {\it universal} for the moduli problem of $n$-pointed stable curves: For any morphism $X\map S$ of $n$-pointed stable curves (i.e., proper flat morphism together with $n$ disjoint sections), there exists a unique morphism $\varphi: S\map \XX_n$ such that the diagram\med

\hs 5cm $X\hs 1cm \map \hs .5cm \XX_{n+1}$\med

\hs 5cm $\downarrow$\hs 2.5cm $\downarrow$\med

\hs 5cm $S\hs 1.1cm \map \hs .7cm \XX_n$,\med

commutes and realizes $X$ as the fibre product $\XX_{n+1}\times_{\XX_n} S$. \med\med\gb


{\Bf 12. Phylogenetic trees} \label{section:phylogenetictrees}\med

A {\it (labelled) phylogenetic tree} with $n$ leaves is a finite non-directed planar graph $T=(V,E,N)$ without cycles or multiple edges and without vertices of valence (= degree) $2$: There are $n$ vertices of valence $1$, called the {\it leaves} of $T$; each of them is equipped with a different {\it label} from a totally ordered set $N$ of cardinality $n$, typically $N=\{1,...,n\}$. The vertices of valence $\geq 3$ will be called {\it inner vertices} of $T$, and the edges joining them, {\it inner edges}. The leaves are often identified with their labels. They are attached by an {\it outer edge} to a (unique) inner vertex. Phylogenetic trees with precisely one inner vertex are called {\it generic}, and those all whose inner vertices have minimal valence $3$ {\it extremal}. A {\it bamboo} is a phylogenetic tree whose inner vertices have only one or two inner edges. See Fig.~12 for various examples of such trees.


\cl{\phylogenetictrees}

\cl{{\it Figure 12.} Some phylogenetic trees. Inner vertices and edges in black, leaves and their edges in red.}\med\gb


{\bf Lemma 1.} \label{lemma:meetingpoints} (Meeting points) {\it For any three pairwise distinct leaves $i,j,k$ of a phylogenetic tree $T$ there is a unique vertex $v$ such that each of the three leaves belongs to a different connected component of $T\sm\{v\}$, see Fig.~120.}\med


We call the vertex $v=v(i,j,k)$ the {\it meeting point} of the three leaves $i,j,k$. It minimizes the sum of the lengths of the (simple) paths connecting $v$ with $i$, $j$, $k$.\med


\cl{\meetingpoints}

\cl{{\it Figure 120.} The meeting points (in blue) of triples of leaves (in yellow).}\med


{\it Proof.} Pick any vertex $v$ in $T$. There are unique simple paths from $v$ to $a$, $b$, $c$. If no two of the three paths share the first edge, $v$ is the meeting point. If not, let $w$ be the other endpoint of the shared edge. The sum of the lengths of the three paths from $w$ to $a$, $b$, $c$ is smaller than the respective sum for $v$. Iterating the process one arrives in finitely many steps at the meeting point. Its uniqueness is obvious. \qed


{\it Destination sets.} If $v$ is a vertex of $T$, consider an edge $e$ emanating from $v$ and the set of leaves $I^e_v\subseteq N$ which can be reached from $v$ by a simple path starting along $e$, see Fig.~121. We call these sets the {\it destination sets} of $v$ (think of a train station with railways lines going out in different directions). For each vertex $v$, we thus get the {\it destination partition} $\DP_v$ of $N$, i.e., $N=\coprod_{I\in \DP_v} I$. \med


\cl{\destinationsets}

\cl{{\it Figure 121.} The destination sets of vertices $v$ in phylogenetic trees.}\med\gb


{\bf Lemma 2.} \label{lemma:destinationsets} (Destination sets) {\it The destination partitions $\DP_v: N=\coprod_{I\in \DP_v} I$ associated to the vertices $v$ of a phylogenetic tree $T$ satisfy the following properties.

\hs 1cm (i) $\abs{\DP_v}\geq 3$ for all $v\in V$;

\hs 1cm (ii) For each $v\in V$ and $I\in \DP_v$ with $\abs I\geq 2$,  there exists a unique inner vertex $w$ in $T$   

\hs 1cm whose destination partition $\DP_w$ contains a set $J$ complementary to $I$, say, $I\sqcup J =N$;

\hs 1cm (iii) For each $i\in N$, there exists a unique inner vertex $v$ in $T$ with $\{i\}\in\DP_v$.}
\med

{\it Proof.} This is immediate by inspection, see Fig.~122 for assertion (ii).\qed\vs-.3cm


\cl{\destinationsetscomplementary}

\cl{{\it Figure 122.} Complementary destination sets $I^e_v$ and $I^e_w$ of adjacent vertices $v$ and $w$.}\med\med\gb


{\it Remark.} The analogous definition of destination sets can be made for any simple path $\gamma$ between two vertices $v$ and $w$ of $T$ by choosing the labels of all leaves which can be reached from $v$ by going along the path $\gamma$, respectively, conversely, starting at $w$. This will become relevant in Section 20.\med


We will apply various operations to phylogenetic trees $T$. Four of them are the following.\med


{\it Deleting edges.} If we take off an inner edge of $T$, but not its endpoints, we get two connected components (which, in general will not be phylogenetic trees). We also say that $e$ {\it separates} the leaves on the left hand side from the leaves on the right hand side of the edge. If we take off, from an extremal tree, an inner edge together with its endpoints, we get four connected components. These decompositions will be used in later sections. See Fig.~124.


\cl{\deletingedgeA}
\cl{\deletingedgeB}

\cl{{\it Figure 124.} Deleting an edge from a tree (without and with endpoints).}\med\gb


{\it Clipping leaves.} If $i$ is (the label) of a leaf of $T$, we may clip it off together with the edge joining it to an inner vertex $v$. This vertex may either have valence $\geq 3$ again, and then we are finished, or valence $2$, in which case it had just one more leaf $j$, the {\it sibling} of $i$. In this case we contract the edge of $j$ and the (unique) inner edge of $v$ to one edge so that $v$ becomes now the leaf $j$ of the new tree. See Fig.~123.\med


\cl{\clippingleavesA}

\cl{\clippingleavesB}

\cl{{\it Figure 123.} Clipping the leaf $a$ from a phylogenetic tree (left) with possible edge contraction (right).}\med\med\gb


{\it Contracting and inserting edges.} Any inner edge can be contracted, together with its endpoints, to a single vertex. The tree remains phylogenetic, and the valence of the resulting vertex is the sum of the valences of the endpoints of the contracted edge. Conversely, any inner vertex of valence $\geq 4$ can be replaced by a new edge, attaching the original edges of the vertex arbitrarily to the two new endpoints, but such that both maintain valence at least $3$. This process thus has several options, each yielding a phylogenetic tree. We also call it the {\it vertex splitting in} or the {\it extension of} the tree. Contracting the new edge produces the original tree, see Fig.~125.


\cl{\contractinginsertingedge}

\cl{{\it Figure 125.} Contracting (left to right) and inserting (right to left) an edge.}\med\gb


In the situation of part (ii) of Lemma 2, we will say that $\DP_v$ and $\DP_w$ {\it possess complementary destination sets} $I=I^e_v$ and $J=I^e_w$. A collection of partitions $\{\PP_v\}_{v\in V}$ of $N$, indexed by a finite set $V$, say, $N=\coprod_{I\in \PP_v} I$ for all $v\in V$, is called an {\it arboral covering} of $N$ if the three properties of the lemma are satisfied. It is not hard to see that the partitions $\DP_v$, $v\in V$, defined by the destination sets determine the phylogenetic tree. More precisely, phylogenetic trees and arboral coverings are the same thing.\med\gb


{\bf Proposition.} \label{proposition:arboralcovering} (Arboral covering) {\it Let $N$ be a finite set of labels. The following two constructions are inverse to each other.

(a) For every phylogenetic tree $T$ with labels $N$ the destination partitions $\DP_v=\{I^e_v,\, e$ an edge of $v\}$ of $N$ associated to the vertices $v$ of $T$ define an arboral covering $\AA_N$ of $N$.

(b) For every arboral covering $\AA_N=\{\PP_v,\, v\in V\}$ of $N$, define a graph $T=(V,E,N)$ as follows: the set of (inner) vertices is $V$, and $N$ is the set of (the labels of the) leaves; no two vertices from $N$ are connected by an edge $e\in E$; a vertex $i\in N$ and a vertex $v\in V$ are connected by an edge $e$ if and only if $\{i\}\in \PP_v$ as in (iii) of the lemma; two vertices $v$ and $w$ from $V$ are connected by an edge $e$ if and only if there exist $I\in \PP_v$ and $J\in\PP_w$ with $N=I\sqcup J$ as in (ii) of the lemma. Then $T$ is a phylogenetic tree with destination partitions $\DP_v=\PP_v$ for all $v\in V$.}\med

{\it Proof.} Part (a) is clear from Lemma 4. Conversely, property (i) ensures that all inner vertices of $T$ have valence $\geq 3$. By (iii), the leaves have valence $1$. Finally, (ii) implies that $T$ has no cycles and is  a tree.\qed


{\it Aside: The number of phylogenetic trees.} (E.~Schr\"oder's fourth problem)  Let $a(n,o)$ be the number of {\it rooted} phylogenetic trees with $n$ leaves and $o$ inner vertices, with $1\leq o\leq n-1$. (Rooted means here that there is one extra vertex (drawn at the bottom) which has valence $1$ but does not count as a leave. Deleting the root from the tree, its adjacent inner vertex may thus have valence $2$, see the pictures in {\cite{\[Fels\]}}, p.~28\&29.) In the definition of $a(n,o)$, the labels $1,2,...,n$ are ordered, so that each permutation of the labels is counted separately, except if the permutation only permutes the labels within the sets of leaves attached to each vertex, but stabilizes the sets. The recursion is, for $o\geq 2$, $n\geq 2$,\med

\cl{$a(n,o)=o\cdot a(n-1,o)+ (o+n-2)\cdot a(n-1,o-1)$,}\med

with $a(n,1)=1$. The total number of pylogenetic trees with $n$ leaves equals \med

\cl{$\ds a(n)=\sum_{o=1}^{n-1} a(n,o)$.}\med\gb

The recursion is \med

\cl{$a(n) = -(n-1) a(n-1) + \sum_{k=1}^{n-1} a(k)\cdot a(n-k)\cdot {n\choose k}$,}\med

with first values \med

\cl{$a(n)=0, 1, 1, 4, 26, 236, 2752, 39208, 660032, 12818912, 282137824, 6939897856, \ldots$}\med

The associated exponential generating function \med

\cl{$A(x)= x + x^2/2! + 4x^3/3! + 26x^4/4! + 236x^5/5! + 2752x^6/6! + \ldots$}\med

satisfies the functional equation $\exp A(x) = 2A(x) - x + 1$ and the non-linear differential equation \med

\cl{$(1 + x - 2 y) \cdot y' = 1$.}\med

For further information, one may consult {\cite{\[Fels\]}}, p.~29, {\cite{\[Schr, Stan\]}}, Ex.~5.2.5, Equ.~(5.27), Fig.~5-3, p.~14, and notes p.~66, {\cite{\[Comt\]}}, p.~224. \med



\med\gb

{\Bf 13. The phylogenetic tree $\Gamma_\xx$ of a string $\xx$} \label{section:treeofstring}\med

To every string $\xx\in \YY_n$ we will now associate a finite graph, its phylogenetic tree $\Gamma_\xx$ with $n$ leaves. This object will be the key combinatorial tool of the whole story: It captures succinctly essential information about $\xx$ and then serves as a touristic guide when one is held to find proofs of statements about $\XX_n=\YY_n$. Recall that, by definition, strings in $\YY_n$ satisfy the equalities of cross-ratios $\cross_q(x^s)=\cross_q(x^t)$ for all quadruples $q\in N^{\ul 4}$ and all triples $s,t\in {N\choose 3}$. As we will prove later that $\XX_n=\YY_n$, we may as well take the string $\xx$ in $\XX_n$.\med

The construction of $\Gamma_\xx$ goes as follows. \med

If $x=(x_1,...,x_n)\in(\P^1)^n$ is an arbitrary $n$-gon, the set  $N=\{1,...,n\}$ of labels decomposes naturally into disjoint subsets $I_j$ collecting the indices $i\in N$ with equal $x_i$. One obtains the {\it incidence partition} \med

\cl{$\IP_x: N=\bigsqcup_j I_j$}\med

of $x$ with {\it incidence sets} $I_j$. If $\abs{I_j}=1$, we call $I_j=\{i_j\}$ a {\it singleton} of $x$. Clearly, $n$-gons which are $\PGL_2$-equivalent have equal incidence partition, so that we can write $\IP_x=\IP_{[x]}$. We will show soon that, conversely, $n$-gons with the same cross-ratios and incidence partition are $\PGL_2$-equivalent (see the proposition below).\med


{\it Example.} \label{example:incidencepartition} For $x=(0,1,1,\infty,\infty, \infty, a)$ with $a\in\P^1\sm \{0,1,\infty\}$ one has \med

\cl{$\IP_x: \{1,2,3,4,5,6,7\}= \{1\}\sqcup \{2,3\}\sqcup \{4,5,6\} \sqcup \{7\}$.}\med\gb


For later use we state a first and somewhat surprising fact.\med

{\bf Lemma 1.} \label{lemma:incidenceequivalence} (Incidence and $\PGL_2$-equivalence) (a) {\it Two $n$-gons $x$ and $y$ with at least three different entries and with the same cross-ratios are $\PGL_2$-equivalent if and only if the have a common incidence set $I$.}

(b) {\it Two $n$-gons $x$ and $y$ with at least three different entries and with the same cross-ratios are $\PGL_2$-equivalent if and only if their incidence partitions $\IP_{x}$ and $\IP_{x}$ coincide.}\med


{\it Proof.} We show first that if $x$ and $y$ have one common incidence set then they already have the same incidence partitions $\IP_x=\IP_y$. This claim reduces the proof to part (b) of the lemma. The argument is a bit tricky.\med

Let $J\neq K$ be two further incidence sets of $x$. These exist since $x$ has at least three different entries. Choose $i\in I$, $j\in J$, $k\in K$ and some $\ell\in N\sm I$. Then $x_j,x_k,x_\ell$ are different from $x_i$, and the analogous condition holds for the components of $y$. Set $q=(ijk\ell)$. We know that $\cross_q(x)=\cross_q(y)$ holds by assumption, provided that no three of the involved components of $x$ and $y$ are equal. Three cases: If $x_\ell=x_k \neq x_j$, then $\cross_q(x)=1$ (see Section 7 on cross-ratios) and hence $y_\ell=y_k$ (we cannot exclude here that $y_\ell=y_k= y_j$ with undefined cross-ratio $\cross_q(y)$, but this does not matter). If $x_\ell=x_j \neq x_k$ then $\cross_q(x)=0$ and hence $y_\ell=y_j$ (and possibly even $y_\ell=y_j=y_k$). Finally, if $x_\ell\neq x_k, x_j$, then $\cross_q(x)\in \P^1\sm\{0,1,\infty\}$, and hence $y_\ell\neq y_k,y_j$. Varying now $j,k,\ell$ we conclude that $y$ has two incidence sets $J'$ and $K'$ containing $J$ and $K$ respectively. Inverting the role of $x$ and $y$ in the above, we end up by symmetry with equalities $J'=J$ and $K'=K$. This shows that $y$ has the same incidence sets as $x$. This proves our claim.\med\gb

So assume that $x$ and $y$ have the same incidence partitions $\IP_x=\IP_y$. We wish to show that $x$ and $y$ are $\PGL_2$-equivalent, the converse implication being obvious. Without loss of generality we may assume that the labels $i,j,k$ lie in different incidence sets $I$, $J$ and $K$ of $x$ and $y$. Up to $\PGL_2$-equivalence, we may then further assume that $x_i=y_i=1$, $x_j=y_j=0$ and $x_k=y_k=\infty$. Choose any $\ell\in N\sm (I\cup J\cup K)$ and set $q=(ijk\ell)\in N^{\ul 4}$. From $\cross_q(x)=\cross_q(y)$ now follows that $x_\ell=y_\ell$ (see again Section 7). As this holds for all $\ell$, we are done, and assertion (b) is proven.\qed\gb


{\it Definition of phylogenetic tree $\Gamma_\xx$.} Let $\xx=(x^t)_{t\in {N\choose 3}}$ be a string in $\YY_n$ with $n$-gons $x^t$. Associate to it a finite (labelled) graph \med

\cl{$\Gamma_\xx=(V,E,N)$}\med

with vertices $V=V(\Gamma_\xx)$, edges $E=E(\Gamma_\xx)$ and labels $N$ as follows: 

{\it \parindent 1.5cm

\litem{}{--} The set of (inner) vertices $P$ consists of $\PGL_2$-orbits $v=[x^t]$ of the $n$-gons $x^t$ of $\xx$; 

\litem{}{--} The set of leaves (or labels) is the set $N$; 

\litem{}{--} Thus $V=P\cup N$; 

\litem{}{--} Two vertices $v=[x^s]$ and $w=[x^t]$ are connected by an edge $e\in E$ if $x^s$ and $x^t$ have 

\hs .3cm an incidence set $I$, respectively, $J$, which are complementary to each other, $I\,\sqcup\, J=N$; 

\litem{}{--} A leaf $i\in N$ is attached to an inner vertex $[x^s]$ by an edge $e\in E$ if $I=\{i\}$ is a  

\hs .3cm singleton of $x^s$, i.e., if $x^s_i$ is different from all other entries $x^s_j$ of $x^s$; 

\litem{}{--} Leaves $i$, $j$ in $N$ are not connected by an edge. \med

}
The graph $\Gamma_\xx$ will be called the {\it incidence graph} of the string $\xx$. By definition, it is completely determined by the incidence partitions $\IP_{x^t}$ of the $n$-gons $x^t$ of $\xx$. We will prove below that it is, as expected, indeed a phylogenetic tree (and we will hence speak of $\Gamma_\xx$ as {\it the} phylogenetic tree of $\xx$).\med


{\it Example.} \label{example:fivegon} The case $n=4$ of strings in $\XX_4$ being left to the reader, let $n=5$. Consider the string $\xx\in\XX_5$ from Section 10 given by\med

\cl{$x^{(123)}=(\ul 0,\ul 1,\ul \infty,a,b),...,x^{(345)}=(1-{a\over b},{1-{a\over b}\over 1-{1\over b}},\ul 0,\ul 1,\ul \infty)$,}\med

with $a,b\in\P^1$. If $a\neq b$ and $a,b\neq 0,1,\infty$, all incidence sets are singletons and the tree $\Gamma_\xx$ is the generic tree $T_*$ with one inner vertex and five leaves, see Fig.~10 in Section 10. If $a=b\neq 0,1,\infty$, the $5$-gons of $\xx$ define two different orbits, for instance $u=[x^{(123)}]=[(0,1,\infty,\infty,\infty)]$ and $v=[x^{(134)}]=[(0,0,1,\infty,c)]$. All others are equal to either one of these two. The tree $\Gamma_\xx$ has two inner vertices, one inner edge, and five leaves, of which two are attached to one vertex, and three to the other, see Fig.~13, left. Finally, if $a=b=\infty$ (and similarly for $a=b=0$ or $a=b=1$), there are three orbits, namely $u=[x^{(123)}]=[(0,1,\infty,\infty,\infty)]$, $v=[x^{(134)}]=[(0,0,1,\infty,\infty)]$, $w=[x^{(145)}]=[(0,0,0,1,\infty)]$. The tree $\Gamma_\xx$ is symmetric, with two edges $e$ and $f$ connecting $u$ with $v$, respectively, $v$ with $w$, and the leaves are distributed in two pairs of two and one singleton, see Fig.~13, right.\med


\cl{\fivelabelstreesspecial}

\cl{{\it Figure 13.} The two incidence graphs associated to the string $\xx$.}\med\gb
 

As a first warm-up, let us recover the incidence relations between the entries of $\xx$ from the geometry of $\Gamma_\xx$.\med



{\bf Lemma 2.} \label{lemma:incidencedestination} (Incidence and destination sets) {\it For every string $\xx=(x^t)_{t\in {N\choose 3}}$ in $\YY_n$, the incidence partitions $\IP_{x^t}$ of the $n$-gons $x^t$ of $\xx$ coincide with the destination partitions $\DP_v$ of the inner vertices $v=[x^t]$ of $\Gamma_\xx$.}\med

{\it Proof.} This is immediate from the definition of $\Gamma_\xx$.\qed


The next result will be the crucial step for the proof that $\Gamma_\xx$ is a phylogenetic tree. It is a first instance where one can nicely observe the beneficial interplay between the $\PGL_2$-action on $n$-gons, the use of cross-ratios, and the combinatorial geometry of $\Gamma_\xx$.\med


{\bf Lemma 3.} \label{lemma:complementaryincidence} (Complementary incidence and singletons) {\it Let $\xx=(x^t)_{t\in{N\choose 3}}$ be a string in $\YY_n$ and let $\Gamma_\xx$ be its incidence graph with label set $N$.}

(a) {\it Let $I\subseteq N$, with $\abs I\geq 2$, be an incidence set of the orbit $v=[x^t]$ of an $n$-gon $x^t$ of $\xx$. There exists a unique orbit $w=[x^s]$ of $\xx$ with complementary incidence set $J=N\sm I$.}

(b) {\it Each singleton label $i\in N$ is attached to precisely one inner vertex $v=[x^t]$ of $\Gamma_\xx$.}
\med

{\it Proof.} (a) Uniqueness follows from the lemma at the beginning of this section. 
For proofing the existence, we can draw the idea from the geometry of $\Gamma_\xx$, see Fig.~130, assuming to know that $\Gamma_\xx$ is a phylogenetic tree. This makes no harm as long as we do not use this information in the proof.\med

Here is the strategy: Let $v$ be a vertex of $\Gamma_\xx$ with some incidence set $I$ with $\abs I\geq 2$. In the picture, $I$ is drawn as the destination set to the left of $v$. We wish to determine the neighbouring vertex $w$ to the right of $v$, since it will have the complementary destination set $N\sm I$. This will also be, by the last lemma, the incidence set of any $n$-gon $x^s$ defining $w$. To find $w$, we first choose {\it some arbitrary} vertex $u$ to the right of $v$ (i.e., in the connected component of $\Gamma_\xx\sm\{v\}$ containing $w$) and then approach $w$ from there stepwise along the path from $u$ to $v$, determining successively vertices $u'$, $u''$,..., until we reach $w$. To find $u$ is easy: choose $i,k\in I$ and $j\in N\sm I$ as in Fig.~130, take $r=(ijk)$ and set $u=[x^r]$. It will have a destination set $J$ containing $N\sm I$. If $J=N\sm I$, we are done, since $J$ is then also an incidence set of $x^s$, say, $w=[x^s]$. Otherwise, replace $k\in I$ suitably by a label $k'\in I\cap J$ such that $u'=[x^{r'}]$ with $r'=(ijk')$ is closer to $w$ than $u$ as indicated in the picture. Eventually, one will arrive at $w$.\med


\cl{\lemmacomplementaryincidence}

\cl{{\it Figure 130.} Vertices $u,u',w$  in $\Gamma_\xx$ with incidence sets $J\supsetneq J'\supsetneq N\sm I$.}\med\med\gb


Let us make this more precise (some patience will be needed). Let $x^r$ with $r=(ijk)$ be the chosen $n$-gon of $\xx$ as described before. By definition of $\YY_n$ as a subset of $\TT_n$, it has (pairwise) different components $x^r_i$, $x^r_j$ and $x^r_k$. We show that, for $\ell\in N\sm I$, all entries $x^r_\ell$ are equal, namely, equal to $x^r_j$. This signifies that  $x^r$ has an incidence set $J$ containing $N\sm I$.\med

Set $q=(ijk\ell)\in N^{\ul 4}$ with arbitrary $\ell\in N\sm I$. As $r=(ijk)$, the cross-ratio $\cross_q(x^r)$ is defined. On the other hand, also the cross-ratio $\cross_q(x^t)$ is defined because $i,k\in I$ and $j,\ell\in N\sm I$ belong to different incidence sets of $x^t$; and $\cross_q(x^t)$ equals $0$ by the formulas in Section 7, using that $x^t_i=x^t_k$ because of $i,k\in I$. As $\xx\in\YY_n$ it follows that also $\cross_q(x^r)=0$ which, in turn, implies $x^r_\ell=x^r_j$ because $x^r_i\neq x^r_k$ are different due to $r=(ijk)$. This shows that the entries $x^r_\ell$ of $x^r$ are equal to $x^r_j$ for all $\ell\in N\sm I$. Thus $x^r$ has an incidence set $J$ containing $N\sm I$. Further, we know that $i,k\not\in J$ because $x^r_i,x^r_k\neq x^r_j$. If $J=N\sm I$ we are done. \med

So assume that $N\sm I \subsetneq J$, i.e., that $I$ and $J$ intersect. In this case, replace $k$ by a label $k'$ in $J\sm(N\sm I)=I\cap J$ and repeat the preceding procedure with $r'=(ijk')$. One obtains an orbit $u'=[x^{r'}]$ with incidence set $J'$ containing again $N\sm I$, see Fig.~130, and with $i, k'\not \in J'$.  We claim that $J'\subseteq J$. If this is the case, we get $J'\subsetneq J$ since $k'\in J$ but $k'\not\in J'$. Hence $N\sm I\subseteq J'\subsetneq J$ will follow and we are done by induction.\med

To see that $J'\subseteq J$, choose the quadruple $q'=(ijk'\ell)$ with $i\in I$, $j\in N\sm I\subseteq J$, $k'\in J\cap I$ as before, and with an arbitrary $\ell\in N\sm J$. Recall that $r=(ijk)$ and $r'=(ijk')$. The cross-ratio $\cross_{q'}(x^r)$ is defined since $j,k'\in J$ and $i,\ell'\not\in J$, which implies that $x^r_i=x^r_{k'}$ is different from $x^r_j$ and $x^r_\ell$; it equals $\infty$ because of $x^r_j=x^r_{k'}$, see Section 7. The cross-ratio $\cross_{q'}(x^{r'})$ is defined because $r'=(ijk')$. It then follows that also $\cross_{q'}(x^{r'})=\infty$, which, in turn, implies $x^{r'}_\ell=x^{r'}_i$ since $x^{r'}_j\neq x^{r'}_{k'}$. As $\ell\in N\sm J$ was arbitrary, this shows that $x^{r'}$ has an incidence set $K$ which contains $N\sm J$. This $K$ cannot be $J'$ because $i\in N\sm J$ but $i\not \in J'$. Therefore $K\subseteq N\sm J'$. It follows that $J'\subseteq N\sm K\subseteq J$. This proves the claim and assertion (a).\med


(b) We now come to singletons. Uniqueness follows from Lemma 1 at the beginning of the section. As for the existence, let $x^t$ be an arbitrary $n$-gon of $\xx$. There is a unique incidence set $I$ of $x^t$ containing $i$. If $I=\{i\}$ we are done. Otherwise, $I$ has at least cardinality $2$ and part (a) applies: There exists a unique $n$-gon $x^s$ of $\xx$ with complementary incidence set $J=N\sm I$. The remaining incidence sets of $x^s$ form a partition of $I$. One of them, say $I'\subseteq I$, contains $i$. But $I'\subsetneq I$ since $x^s$ has at least one more incidence set. The claim now follows by induction on $\abs I$.\qed

\ignore
[\ooo optional]\med

(c) We show the three properties of arboral coverings. We have already seen that the incidence partitions of $\xx$ coincide with the destination partitions of $\Gamma_\xx$. But as each $n$-gon $x^t$ of $\xx$ has at least three different entries, the vertex $v=[x^t]$ of $\Gamma_\xx$ has at least three destination sets. This proves property (i) of arboral coverings. Property (ii) follows from assertion (a) of the proposition, and (iii) from (b).\qed
\recognize
\med

{\bf Proposition.} \label{proposition:incidencegraphisphylogenetictree} (Incidence graph is phylogenetic tree) (a) {\it The incidence graph $\Gamma_\xx=(V,E,N)$ attached to a string $\xx\in\YY_n$ is a phylogenetic tree with $n$ leaves.}

(b) {\it Every (labelled) phylogenetic tree with $n$ leaves arises as the incidence graph $\Gamma_\xx$ of a string $\xx\in\YY_n$.}\med


It is now manifest to call the incidence graph $\Gamma_\xx$ of $\xx$ {\it the phylogenetic tree} of $\xx$. Assertion (b) is interesting per se, but will not be used in the rest of the paper.\med\gb


{\it Proof.} Assertion (a) follows from the proposition together with the 
considerations preceding it. For assertion (b), we apply induction on $n$. The case $n=3$ is obvious since then the tree must be generic with one vertex and three leaves $1,2,3$. It is the tree of the unique string $\xx=(x^{123})\in\TT_3$, where $x^{123}=(0,1,\infty)$.\med

So assume that $n\geq 4$. Pick an extremal inner vertex $v$ of $T$, i.e., one which is attached to a single other inner vertex. It must have at least two leaves, and without loss of generality we may assume that one of it has label $n$. Clipping off the leaf from $T$ produces a phylogenetic tree $T'$ with $n-1$ leaves and labels $N'=N\sm\{n\}$, see Section 12 for the precise construction. By induction on $n$, let $\yy\in \YY_{n-1}$ be a string of $(n-1)$-gons $y^s$, $s\in {N'\choose 3}$, with $\Gamma_\yy=T'$. See Fig.~131 for the case where an edge of $T$ has to be contracted since $v$ has only two leaves.\med


\cl{\clippingleaf}

\cl{{\it Figure 131.} Clipping off the leaf with label $n$ from the tree $T$ (left) producing the tree $T'=\Gamma_\yy$ (right).}\med\med


We now define from $\yy$ a string $\xx=(x^t)_{t\in {N\choose 3}}\in\YY_n$ such that $\Gamma_\xx=T$. Choose for every vertex $w$ of $T$ a triple $s$ defining it as its meeting point, in the following manner: (i) If $w=v$, take $s=(i\ell n)$, where $i\in N\sm L$ and $\ell\in L$ are chosen arbitrarily. (ii) If $w\neq v$, take any $s\in {N'\choose 3}$ defining it (as $n$ has at least one sibling $\ell$ in $L$, it is not necessary to use $n$ for defining $w$ as meeting point). Observe that in this latter case all entries $y^s_\ell$ of $y^s$ are equal, for $\ell\in L$. So assume that a triple $s$ has been chosen for each vertex $w$. We then define for each of these $s$ an $n$-gon $x^s$ by setting, according to the two cases,\med

(i) \hs .5cm $x^s_j=0$ for all $j\in N\sm L$, $x^s_\ell=1$, $x^s_n=\infty$, $x^s_k\neq 1,\infty$ pairwise different for all $k\in L\sm\{\ell,n\}$;\med

(ii) \hs .4cm $x^s=(y^s,x^s_n)$ with $x^s_n=y^s_\ell$ for $\ell\in L$.\med

For the remaining triples $t$, different from one of the triples $s$, there is a always a unique triple $s$ from the list above defining the same vertex as $t$ in $T$. Then choose $x^t$ as the unique $n$-gon $\PGL_2$-equivalent to $x^s$ with prescribed values $0,1,\infty$ at the entries given by the labels of $t$.\med

This defines a string $\xx=(x^t)_{t\in{N\choose 3}}$ of $n$-gons in $\TT_n$. We are left to show that $\xx\in \YY_n$, i.e., that the $n$-gons $x^t$ have equal cross-ratios. So let $x^r$ and $x^t$ be two $n$-gons defining two different vertices $u=[x^r]$ and $w=[x^t]$ of $T$. If both $u$ and $w$ are different from the chosen vertex $v$, we are in case (ii). Let $q$ be a quadruple in $N^{\ul 4}$. If $q$ does not involve $n$, we have $\cross_q(x^r)=\cross_q(y^r)$ and $\cross_q(x^t)=\cross_q(y^t)$ and the assertion follows from $\yy\in \YY_{n-1}$. If $q$ does involve $n$, but no other label from the set $L$ of labels of $v$, one may replace $n$ by any $\ell\in L$ and get the same cross-ratios, by definition (ii) above. If, finally, $q$ involves $n$ and a further label $\ell$ of $L$, say, without loss of generality, in the last two entries of $q$, both cross-ratios $\cross_q(x^r)$ and $\cross_q(x^t)$ are equal to $1$.\med

We are left with the case where $u\neq v$ and $w= v$, say case (i). If $q$ involves more than two labels from $N\sm L$, the cross-ratio $\cross_q(x^t)$ is not defined and nothing is to show. If $q$ involves exactly two labels from $N\sm L$, say, at the first two entries, the cross-ratio $\cross_q(x^t)$ equals $1$. But as $q$ has in this case two labels from $L$, namely, at the last two entries, also $\cross_q(x^r)=1$. Finally, if $q$ involves more than two labels from $L$, the cross-ratio $\cross_q(x^r)$ is not defined and nothing is to show. This proves assertion (b) of the proposition.\qed 


{\bf Corollary.} \label{corollary:equivalentngons} {\it Two $n$-gons $x^s$ and $x^t$ of a string $\xx\in\YY_n$ are $\PGL_2$-equivalent if and only if there is a triple $r=(ijk)\in {N\choose 3}$ such that  $x^s_i, x^s_j, x^s_k$, and, respectively, $x^t_i,x^t_j,x^t_k$, are pairwise different.}\med

{\it Proof.} By the definition of meeting points in a phylogenetic tree, the assertion is immediate from part (a) of the proposition since both $[x^s]$ and $[x^t]$ are the meeting point of the same triple $r$.\qed


{\bf Lemma 4.} \label{lemma:characterizationtree} (Characterization of tree) {\it For any phylogenetic tree $T$, let $Q_T^1$ be the set of quadruples $q=(ijk\ell)\in N^{\ul 4}$ such that there is an inner edge $e$ in $T$ separating $i,j$ from $k,\ell$. Then, for a string $\xx\in\XX_n$, one has $\Gamma_\xx=T$ if and only if $\cross_q(\xx)=1$ for $q\in Q_T^1$ and $\cross_q(\xx)\neq  1$ for $q\not \in Q_T^1$.}\med

{\it Proof.} Recall that if $q=(ijk\ell)$ then $\cross_q(\xx)=1$ if and only if, for all triples $x^t$ of $\xx$ for which $\cross_q(x^t)$ is defined, one has either $x^t_i=x^t_j$ or $x^t_k=x^t_\ell$. Assume first that $\Gamma_\xx=T$. Let $e$ be any inner edge of $T$. Then, for all quadruples $q=(ijk\ell)$ such that $i,j$ are separated by $e$ from $k,\ell$ in $T$ one has $\cross_q(x^t)=1$, for all triples $t$ such that $v=[x^t]$ is one of the endpoints of $e$. Thus,  $\cross_q(\xx)=1$ for all $q\in Q_T^1$. For all other quadruples $q=(ijk\ell)$ one has, accordingly, $\cross_q(\xx)\neq 1$. \med




Conversely, let $\xx$ be a string in $\XX_n$ such that $\cross_{(ijkl)}(\xx)=1$ if and only if $T$ has an edge seperating $\{i,j\}$ from $\{k,l\}$. But if $\cross_{(ijkl)}(\xx)=1$, then, by its very definition, $\Gamma_\xx$ has an edge seperating $\{i,j\}$ from $\{k,l\}$. If $\cross_{(ijkl)}(\xx)\ne 1$, then $\Gamma_\xx$ has no such edge. So for any edge of $T$, and any  $\{i,j\}$ and $\{k,l\}$ separated by that edge, there is an edge of $\Gamma_\xx$ separating $\{i,j\}$ from $\{k,l\}$, and conversely. Therefore $\Gamma_\xx=T$. \qed\med\med\gb



{\Bf 14. The phylogenetic tree $\Gamma_C$ of an $n$-pointed stable curve $C$} \label{section:treeofcurve}\med

We specify the assertions made in Part I, Section 6. The {\it dual graph} of an $n$-pointed stable curve $C$ has as vertices the components $C_i$ of $C$. Two vertices are connected by an edge if and only if the respective components intersect. The dual graph neglects the position of the marked points. It is a tree, that is, a connected finite graph without loops.\med

The {\it augmented dual graph} $\Gamma_C$ of an $n$-pointed stable curve $C$ is obtained from the dual graph by attaching to each vertex as many edges as there are marked points on the corresponding irreducible component and by labelling the endpoint of each of these edges by the respective marked point $p_i$. We call $N=\{p_1,...,p_n\}$ the set of labels of $\Gamma_C$. It is just considered as an abstract set, and ignores any information about the precise location of the marked points on the components, see Fig.~14.\med


\cl{\stablecurvesandtree}\vs-.3cm

\cl{{\it Figure 14.} An $n$-pointed stable curve $C$ (left) and its phylogenetic tree $\Gamma_C$ (right).}\med\gb

 
{\bf Proposition 1.} \label{proposition:dualgraph} (Dual graph) {\it The augmented dual graph $\Gamma_C$ of an $n$-pointed stable curve $C$ is a phylogenetic tree with labels given by the marked points.}\qed


We will therefore call $\Gamma_C$ the phylogenetic tree of the $n$-pointed stable curve $C$.\med


Meeting points (see Section 12) allow us to make the following observation on the geometric configuration of stable curves with respect to their marked points. Define a {\it path} in $C$ between a marked point $p\in C_i$ and an irreducible component $C_j$ of $C$ as the sequence of irreducible components of $C$ given by the vertices of the path in $T_C$ connecting the two vertices corresponding to $C_i$ and $C_j$. \med


{\bf Proposition 2.} \label{proposition:projectionmediancomponent} (Projection to median component) {\it Let $T=\Gamma_C$ be the phylogenetic tree of an $n$-pointed stable curve with label set $N$. For every triple $t=(ijk)\in N$ let $v^t$ be the vertex of $T$ where the marked points $p_i$, $p_j$, $p_k$ meet, and denote by $C^t$ the irreducible component of $C$ specified by $v^t$. Let $D_i$, $D_j$, $D_k$ be the connected curves of $C\sm C^t$ containing $p_i$, $p_j$, $p_k$, and denote by $p'_i$, $p'_j$, $p'_k$ their intersection points with $C^t$. Then the three points $p'_i$, $p'_j$, $p'_k$ on $C^t$  are pairwise distinct.}\med

We call $p'$ the {\it projection} of a marked point $p$ of $C$ onto $C^t$. Obviously, $p'=p$ if and only if $p\in C^t$. The irreducible component $C^t$ is known as the {\it median component} associated to the triple $t$, see for instance {\cite{\[GHP\]}}, p.~133.\med


{\it Proof.} See Fig.~140.\qed\vs-.2cm


\cl{\medianprojection}\med

\cl{{\it Figure 140.} The projection points (green) of marked points (red) onto the median component $C^{(ijk)}$.}\med

\med\med\gb


{\Bf 15.~Constructing a string $\xx$ from a stable curve $C$} \label{section:stringofcurve}\med

The last proposition allows us to associate, to every $n$-pointed stable curve $C$ with marked points $p_\ell$, $\ell\in N$, and every triple $t=(ijk)$ in $N$, the median component $C^t$ of $C$ and $n$ points $p'_\ell$ on $C^t$ such that $p'_i$, $p'_j$, $p'_k$ are pairwise distinct. There is then a unique isomorphism $C^t\map \P^1$ sending $p'_i$, $p'_j$, $p'_k$ to $0$, $1$, $\infty$. In this way we get an $n$-gon $x^t$ in $(\P^1)^n$ whose entries are the images of the points $p'_\ell$ under this isomorphism. This $n$-gon clearly only depends on the isomorphism class of $C$. Varying over all triples $t$ we end up with a string $\xx=(x^t)_{t\in {N\choose 3}}$ in $\TT_n$, that is, every $n$-gon $x^t$ satisfies $x^t_i=0$, $x^t_j=1$, $x^t_k=\infty$, for $t=(ijk)$. \med


{\bf Proposition.} \label{proposition:stringstablecurve} (String of stable curve) {\it The map from the set of isomorphism classes of $n$-pointed stable curves to the projective variety $\TT_n$ sending $C$ to the string $\xx=(x^t)_{t\in {N\choose 3}}$ as defined above is injective with image included in $\YY_n$.}\med

{\it Proof.} Let $C$ be a stable curve with marked points $p_\ell$, $\ell\in N$. Let $t=(ijk)$ be a triple in ${N\choose 3}$. Let $C^t$ be the median component of $p_i, p_j, p_k$, and let $p_1^t,...,p_n^t$ be the projections of the points $p_1,...,p_n$ to $C^t$. By Proposition 2 in Section 14, the three projections $p_i^t,p_j^t,p_k^t$ are pairwise different. Hence there is a unique isomorphism $\varphi^t:C^t\map\P^1$ such that $\varphi^t(p_i^t)=0,\varphi^t(p_j^t)=1,\varphi^t(p_k^t)=\infty$. Set\med 

\cl{$x^t=(\varphi^t(p_1^t),\dots,\varphi^t(p_n^t))\in(\P^1)^n$.} \med

This gives the string $\xx=\xx_C=(x^t)_{t\in {N\choose 3}}$. It remains to show that $\xx\in\YY_n$. For this, it suffices the show that the $n$-gons of $\xx$ have equal cross-ratios. This goes as follows.\med

If $s$ and $t$ are two triples with the same median components $C^s=C^t$, then the two $n$-gons $x^s$ and $x^t$ are $\PGL_2$-equivalent, and therefore
all their cross-ratios are equal. Assume now that $C^s\neq C^t$. Let $q=(ijk\ell)$ be a quadruple such that $\cross_q(x^s)$ and $\cross_q(x^t)$ are defined.  No three of the four points in $C$ selected by $i,j,k,\ell$ can have pairwise distinct projections in both $C^s$ and $C^t$, because of the uniqueness of the median component. If there are two among the four points which have equal projections in both $C^s$ and $C^t$, then it follows again that $\cross_q(x^s)=\cross_q(x^t)$. Up to renaming the labels in $q$, there is only one case left:\med

\cl{$x^s_\ell\neq x^s_i=x^s_j\neq x^s_k\neq x^s_\ell$ and $x^t_j\neq x^t_k=x^t_\ell\neq x^t_i\neq x^t_j$.}\med

As $\cross_q(x^s)=\cross_q(x^t)=1$, the string $\xx$ belongs to $\YY_n$.\qed

\med\med\gb




{\Bf 16. The main theorem about $\XX_n$} \label{section:maintheorem}\med

Here is the phylogenetic or ``stringy'' version of the result of Deligne, Mumford and Knudsen on the moduli of stable curves of genus $0$. It will imply that $\XX_n$ is in one-to-one correspondence with the compactification $\ol{\MM_{0,n}}$ and that it has the asserted properties. Recall that \med

\hs 2cm $\UU_n=((\P^1)^n\sm \Delta_n)/\PGL_2$, orbit space of generic $n$-gons,\med

\hs 2cm $\Sigma_n:\UU_n\map (\P^1)^{n{n\choose 3}}$, $[x]\map \xx=(x^t)_{t\in {N\choose 3}}$, symmetrization map,\med

\hs 2cm $\TT_n\subseteq(\P^1)^{n{n\choose 3}}=\{\xx=(x^t)_{t\in{N\choose 3}},\, x^t_i=0, x^t_j=1, x^t_k=\infty$ for $t=(ijk)\in{N\choose 3}\}$,\med

\hs2cm $\YY_n=\{\xx\in\TT_n,\,\cross_q(x^s)=\cross_q(x^t),\,$ for all $s,t\in {N\choose 3},\, q\in N^{\ul 4}\}\subseteq \TT_n$,\med

\hs 2cm $\XX_n=\ol{\Sigma(\UU_n)}\subseteq \YY_n\subseteq\TT_n$, Zariski-closure of the set of generic strings,\med

\hs 2cm $N=\{1,...,n\}$, $N^{+1} =\{1,...,n,a\}$, sets of labels,\med

\hs2cm $\pi_a: \XX_{n+1} \map \XX_n$, projection map, forgetting entries involving $a$,\med

\hs2cm $\sigma_p: \XX_n \map \XX_{n+1}$, $p=1,...,n$, disjoint sections of $\pi_a$.\med\med\gb


{\bf Main Theorem.} \label{theorem:maintheorem} (1) {\it The Zariski-closure $\XX_n$ of $\Sigma(\UU_n)$ in $\TT_n$ equals $\YY_n$; it is a smooth, closed and irreducible subvariety of $(\P^1)^{n{n\choose 3}}$ of dimension $n-3$.}  

(2) {\it As a variety, $\XX_n$ is naturally stratified by smooth locally closed strata $\SS_{n,T}$, each consisting of strings $\xx$ with the same phylogenetic tree $\Gamma_\xx=T$, the generic tree $T_*$ corresponding to the dense open stratum $\UU_n\subset\XX_n$. The smallest strata consist of a single string $\xx$, and their trees are extremal (every inner vertex of valence $3$). Adjacent strata have phylogenetic trees related to each other by contraction of edges, respectively, extension of vertices.}

(3) {\it The boundary $\BB_n=\XX_n\sm\UU_n$ is a divisor in $\XX_n$ with normal crossings. Its components are smooth hypersurfaces $\DD_{I,J}$ indexed by pairs $(I,J)$ of complementary subsets $I, J$ of $N$ of cardinality $\geq 2$.}

(4) {\it The projection $\pi_a: \XX_{n+1} \map \XX_n$ is a flat projective morphism of algebraic varieties with one-dimensional reduced fibers.}

(5) {\it The strings $\xx=(x^t)_{t\in{N\choose 3}}$ in $\XX_n$ are in one-to-one correspondence with isomorphism classes of $n$-pointed stable curves $C$ of genus zero:}

\hs.6cm (i) {\it The projections of the $n$ marked points on $C$ to the median components $C^t$,  $t\in{N\choose 3}$, define   

\hs.6cm  the $n$-gons $x^t$ of a string $\xx\in\XX_n$ whose phylogenetic tree $\Gamma_\xx$ equals the augmented dual graph

\hs.6cm  $\Gamma_C$ of $C$.}

\hs.6cm (ii) {\it The fibers $\FF_\xx=\pi_a\inv(\xx)$ of  the \projection $\pi_a: \XX_{n+1} \map \XX_n$, marked with the $n$ images of  

\hs.6cm suitably chosen sections $\sigma_p$ of $\pi_a$, are $n$-pointed stable curves with augmented dual graph $\Gamma_{\FF_\xx}$  

\hs.6cm equal to the phylogenetic tree $\Gamma_\xx$ of $\xx$.}

\hs.6cm (iii) {\it The two operations are inverse to each other.}

(6) {\it For any flat and proper family $f:X\map S$ of $n$-pointed stable rational curves over a base $S$ there exists a unique morphism $g: S\map \XX_n$ such that $f$ is the pull-back $S\times_{\XX_n}\XX_{n+1}\map S$ of $\pi_a:\XX_{n+1}\map\XX_n$ under $g$.}

(7) {\it The variety $\XX_n$ is a fine moduli space for isomorphism classes of $n$-pointed stable curves of genus zero, with universal family $\pi_a: \XX_{n+1} \map \XX_n$. It is hence isomorphic to the Deligne-Mumford-Knudsen compactification $\ol\MM_{0,n}$ of $\MM_{0,n}$.}


\med\med

\med\gb

\cl{\Bf Part III: Proofs}\label{part:proofs}\med

As indicated in the introduction, the proofs of assertions (6) and (7) of the theorem require more advanced techniques and have to be omitted. The other assertions will be proven in all details in Part III together with the construction of the sections $\sigma_p$. Assertion (1) is proven in Section 17 and 22, assertion (2) in Section 18, assertion (3) in Section 19, assertion (5) in Sections 19 - 21. In Section 23 and 24 the proof of assertions (6) and (7) is briefly outlined. The most interesting (and also most challenging) proofs concern the smoothness of $\XX_n$ and the fact that the fibers of $\pi_a:\XX_{n+1}\map\XX_n$ are stable curves. \med

As a general principle, we tend to design the proofs pictorially by investigating the geometry of the associated phylogenetic trees rather than doing blindly algebraic or combinatorial computations with cross-ratios. This will allow the reader to capture much better the flavour of the arguments. It will then be a straightforward task to turn the reasoning into a rigorous  formal proof.\med

\med\med\gb



{\Bf 17. The smoothness of $\XX_n$} \label{section:smoothness}\med



The proof goes in several stages and will cover the whole section. All constructions rely on applying specific operations to the phylogenetic trees and to exploit then the resulting combinatorial structure.%
\tex
{{\baselineskip 12pt \footnote{$\!\!{}^5$}{\smallTimes  For the proof of the smoothness of the compactification $\ol{\frak M}_{0,n}$ considered by Brown and based on different combinatorial constructions, see {\cite{\[Br\]}}, Section 2.5 and 2.8.}}}
%
{{\footnote{For the proof of the smoothness of the compactification $\ol{\frak M}_{0,n}$ considered by Brown and based on different combinatorial constructions, see {\cite{\[Br\]}}, Section 2.5 and 2.8.}}}
%
\med


(a) {\it The known facts.} We assume that we have already defined $\XX_n$ as the Zariski closure of the set $\UU_n=((\P^1)^n\sm \Delta_n)/\PGL_2$, embedded into the set of strings $\TT_n\subseteq (\P^1)^{n{n\choose 3}}$ by symmetrization. Further, $\YY_n$ was defined as the subvariety of $\TT_n$ given by the equality of cross-ratios $\cross_q(x^s)=\cross_q(x^t)$ of the $n$-gons of strings $\xx$. It is clear that $\XX_n$ is a subvariety of $\YY_n$, and we will prove in Section 22 that $\XX_n$ equals $\YY_n$. This proof of the equality does not rely on the arguments in this section, so we will assume for convenience and without loss of generality that $\XX_n=\YY_n$ (more accurately, we will prove in this section the smoothness of $\YY_n$, and via Section 22 and $\XX_n=\YY_n$ the smoothness of $\XX_n$ will follow).\med 

For every quadruple $q$, the cross-ratio function $\cross_q:\YY_n\map \P^1$ is defined on whole $\YY_n$: there is always an $n$-gon $x^t$ in each string $\xx$ where $\cross_q(x^t)$ is defined, and one then sets $\cross_q(\xx):=\cross_q(x^t)\in\P^1$.\med

Also, we already introduced the concept of an extremal string $\xx$ in $\YY_n$ in terms of its phylogenetic tree $\Gamma_\xx$. It can be characterized by one of the following properties (being equivalent with each other):\med

{\parindent 1.1cm

\litem {\hs .5cm (i)} Every inner vertex of $\Gamma_\xx$ is of degree 3, say, has $3$ emanating edges;

\litem {\hs .5cm (ii)} $\Gamma_\xx$ has $2n-3$ edges and $2n-2$ vertices (i.e., $n$ leaves and $n-2$ inner vertices);

\litem {\hs .5cm (iii)} Every $n$-gon of $\xx$ has exactly three distinct entries (with values $0$, $1$, and $\infty$);

\litem {\hs .5cm (iv)} For all quadruples $q$, the cross-ratio $\cross_q(\xx)$ is {\it special}, i.e., equal to $0$, $1$, or $\infty$.

}\med

The associated phylogenetic trees will also be called {\it extremal}. Extremal strings will be used to define an open covering $\{O_\xx,\, \xx$ extremal$\}$ of $\YY_n$ and, eventually, also the chart maps $\alpha_\xx:O_\xx\map (\kk^*)^{n-3}$. We will see later that extremal strings correspond precisely to the zero-dimensional strata of the stratification of $\YY_n$ whose strata are given by the constancy of the phylogenetic tree. In particular, extremal strings are uniquely determined by their phylogenetic tree: If $\Gamma_\xx=\Gamma_\yy$ is extremal, then $\xx=\yy$.\med

It will be convenient to distinguish in the sequel, for a given quadruple $q=(ijk\ell)$,  between the {\it formal} cross-ratio as an (abstract) element \med

\cl{$\ds [q]:= [ijk\ell]={(\xi_i-\xi_k)(\xi_j-\xi_\ell)\over (\xi_i-\xi_\ell)(\xi_j-\xi_k)}$}\med

of the field $\kk(\xi):=\kk(\xi_1,...,\xi_n)$ of rational functions in variables $\xi_i$, and the associated cross-ratio {\it function} \med

\cl{$\cross_q=\cross_{ijk\ell}: (\P^1)^n\sm\nabla^3(ijk\ell)\map \P^1$,}\med

\cl{$\ds (x_1,...,x_n)\map {(x_i-x_k)(x_j-x_\ell)\over (x_i-x_\ell)(x_j-x_k)}$,}\med

where $\nabla^3(ijk\ell)\subseteq (\P^1)^n$ denotes the set of $n$-gons $x=(x_1,...,x_n)$ with at least three equal entries among $x_i, x_j,x_k,x_\ell$. When we briefly say ``cross-ratio'' it should always be clear from the context whether we mean the formal cross-ratio or the cross-ratio function.\med



(b) {\it Charts on $\XX_n$.} For every extremal string $\xx\in\YY_n$, we define an open neighborhood $O_\xx\subset\YY_n$ of $\xx$ and a regular map $\alpha_\xx:O_\xx\map (\kk^*)^{n-3}$. Both the open neighborhood and the map are prescribed in terms of cross-ratios. We will show that the sets $O_\xx$ form an open covering of $\YY_n$ and that each map is a (biregular) isomorphism onto some Zariski-open subset $V_\xx$ in $(\kk^*)^{n-3}$. \med


The open neighborhood $O_\xx$ of an extremal $\xx$ will consist of all strings $\yy\in\YY_n$ whose cross-ratio $\cross_q(\yy)$ is either non-special (i.e., $\neq 0,1,\infty$) or, if special, then equal to $\cross_q(\xx)$ (recall that $\cross_q(\xx)\in\{0,1,\infty\}$ since $\xx$ is extremal). Clearly, $O_\xx$ is an open subset of $\YY_n$ containing $\xx$.\med



{\bf Lemma 1.} \label{lemma:opencovering} (Open covering) {\it The sets $O_\xx$, with $\xx$ ranging over all extremal strings, form an open covering of $\YY_n$.}\med


{\it Proof.} Let $\yy\in\YY_n$ be a string with tree $\Gamma_\yy$. If it is not extremal, then there is an inner vertex $v$ of degree at least four in $\Gamma_\yy$. Then we can perform a vertex splitting, replacing this inner vertex by a new edge $(v_1,v_2)$ and connecting each edge of $\Gamma_\yy$ adjacent to the replaced vertex $v$ by an edge adjacent to either one of the two vertices $v_1$ or $v_2$ of the replacing edge, see Fig.~17 as well as Section 12. Repeating this process, we construct (in a non unique way) a phylogenetic tree $T$ with $n$ leaves whose inner vertices all have degree three, i.e., which is an extremal tree. By part (b) of the proposition in Section 13 about incidence graphs there exists a (non-empty) stratum $\SS_T$ in $\YY_n$ consisting of strings with phylogenetic tree $T$. But as $T$ is extremal, the stratum is reduced to a single string, call it $\xx$. The proof of this fact is an exercise (or see the proposition about stratification in Section 19). We claim that $\yy\in O_{\xx}$.


\cl{\vertexsplitting}

\cl{{\it Figure 17.} Vertex splitting at the inner vertex $v$ (left) replacing it by the edge $e$ (right).}\med


Let $q$ be a quadruple, and assume that $\cross_q(\yy)$ is special and $\neq \cross_q(\xx)$. We will derive a contradiction. Write $q=(ijkl)$. We may assume, up to a permutation of the labels, that $\cross_q(\xx)=1$, and, then, that $\cross_q(\yy)=0$ or $\cross_q(\yy)=\infty$. Up to symmetry, let $\cross_q(\yy)=0$. This signifies, by the formulas for the cross-ratios, that $\Gamma_{\yy}$ has an inner edge $e$ separating $\{i,k\}$ from $\{j,l\}$ (i.e., if you remove the edge $e$, then $i$ and $k$ lie in one connected component and $j$ and $k$ lie in the other). The vertex splitting procedures preserves this property. Therefore also $\cross_q(\xx)=0$ holds, and we get a contradiction. This shows that $\yy\in O_\xx$.\qed


We now come to the definition of the chart map $\alpha_\xx$. It depends on the choice of a total ordering of the set $N$ of leaves. Recall that the tree
$\Gamma_\xx$ is extremal and thus has $n-3$ inner edges. For each such edge $e$, we define a quadruple $q_e$ in $N^{\ul 4}$ as follows: Remove $e$ and its two endpoints from the tree. One obtains a forest (finite disjoint union of trees) with four connected components. Let $\{i,j,k,l\}$ be the smallest leaves of $\Gamma_\xx$ in each of these components. Then set $q_e:=(ijkl)$, where the numeration of the four labels $i,j,k,\ell$ is chosen such that the associated cross-ratio satisfies $\cross_{q_e}(\xx)=1$ (this choice is just for notational convenience, also $0$ or $\infty$ would work). There are eight possible ways to numerate $i,j,k,\ell$ in $q_e$, any choice is fine, so just pick one for each edge $e$. \med

The $n-3$ quadruples $q_e=(ijk\ell) $ selected by these choices, with $e$ ranging over all inner edges of $\Gamma_\xx$, respectively, the associated (formal) cross-ratios $[q_e]=[ijk\ell]$, will be called the {\it edge quadruples}, respectively, the {\it edge cross-ratios} of $\xx$. Recall that their definition depends on the choice of a total ordering of the set of labels $N$, and that we may permute the four entries as long as the equality $\cross_{q_e}(\xx)=1$ is preserved.\med

Denote by $\cross_{q_e}:O_\xx\map \kk^*$ the associated cross-ratio map. By our choices, it is well-defined on the whole neighborhood $O_\xx$ of $\xx$. The chart map $\alpha_\xx:O_\xx\map (\kk^*)^{n-3}$ is then defined by sending a string $\yy$ in $O_\xx$ to the vector of evaluations $\cross_{q_e}(\yy)$ of the $\xx$-edge cross-ratios $\cross_{q_e}$,\med

\cl{$\alpha_\xx: O_\xx\map (\kk^*)^{n-3}, \, \yy\map (\cross_{q_e}(\yy))_{q_e\,  {\rm edge\, quadruple}}$.}\med

Changing the numeration of the labels $i,j,k,\ell$ of $q_e$, but subject to the equality $\cross_{q_e}(\xx)=1$,  may result in having one of the components of $\alpha_\xx$ replaced by its multiplicative inverse. The function $\alpha_\xx$ maps $O_\xx$ into $(\kk^*)^{n-3}$, so the replacement of any component by its multiplicative inverse is just a composition with the isomorphism that maps one of the coordinates in $(\kk^*)^{n-3}$ to its multiplicative inverse.\med\gb



(c) {\it Preview on cross-ratio identities.} To prove that $\alpha_\xx$ is injective, we will show that, for each $\yy\in O_\xx$, the values $\cross_q(\yy)$ of the cross-ratios of $\yy$ with respect to arbitrary quadruples $q$ can be expressed as rational functions in the values $\cross_{q_e}(\yy)$ of the edge cross-ratios $\cross_{q_e}$ of $\xx$. This, in turn, implies that the entries of the $n$-gons of $\yy$ - which are actually themselves values of cross-ratios - are determined by all $\cross_{q_e}(\yy)$, say, by $\alpha_\xx(\yy)$, thus establishing the injectivity of $\alpha_\xx$. \med


Denote by \med

\cl{$\CR_n(\xi):=\kk([q],\, q\in N^{\ul 4}$ a quadruple$)$}\med

the subfield of $\kk(\xi_1,...,\xi_n)$ generated by all formal cross-ratios $[q]=[ijk\ell]$. We will first show that $\CR_n(\xi)$ is generated, for all extremal strings $\xx$, or, equivalently, for all extremal trees $T$ with $n$ leaves, by the edge cross-ratios $[q_e]$ of $\xx$, respectively, $T$, for $e$ an inner edge of the tree $\Gamma_\xx$, respectively, of $T$,\med

\cl{$\CR_n(\xi)=\kk([q_e],\, e$ an edge of $\Gamma_\xx)$.}\med

Note here that the presentation holds for any extremal string $\xx$. It follows that every formal cross-ratio $[q]$ can be written as a rational expression \med

\cl{$\ds [q] = {{P_q([q_1],...,[q_{n-3}])}\over {Q_q([q_1],\dots,[q_{n-3}])}}$}\med

of the edge cross-ratios $[q_1],...,[q_{n-3}]$ of $\xx$. Here it suffices to consider quadruples $q$ with $\cross_q(\xx)=1$, since the remaining ones can be obtained by a permutation of the entries of $q$, with the already discussed effect on the cross-ratios (being one of the six rational functions $c, 1/c, 1-c, {1\over 1-c}, {c\over c-1},{c-1\over c}$). In a second step, we will show that this formal identity carries over to the associated cross-ratio maps, when restricted to the open neighborhoods $O_\xx$ of $\xx$. In this way it will follow that any string $\yy$ in $O_\xx$ can be reconstructed from its image $\alpha_\xx(\yy)$ under $\alpha_\xx$, i.e., that $\alpha_\xx$ is injective. Moreover, one can determine the image $V_\xx\subseteq (\kk^*)^{n-3}$ of $O_\xx$ under $\alpha_\xx$, using the fact that for all $\yy\in O_\xx$ the evaluation $\cross_q(\yy)$ is either non-special or equal to $\cross_q(\xx)$. Assuming without loss of generality that $\cross_q(\xx)=1$ for all $q$ (and hence $\cross_q(\yy)\neq 0,\infty$ for all $\yy\in O_\xx$, by definition of $\O_\xx$), we will have to show more specifically the inequalities  \med

\cl{$P_q(\cross_{q_1}(\yy),...,\cross_{q_{n-3}}(\yy))\neq 0$}\med

and \med

\cl{$Q_q(\cross_{q_1}(\yy),...,\cross_{q_{n-3}}(\yy))\neq 0$}\med

for all $\yy\in O_\xx$ (both are a priori $\neq \infty$ since all edge cross-ratio values $\cross_{q_e}(\yy)$ are $\neq \infty$, again by the definition of $O_\xx$). If this holds, the quotient\med

\cl{$\ds {P_q(\cross_{q_1}(\yy),...,\cross_{q_{n-3}}(\yy))\over Q_q(\cross_{q_1}(\yy),...,\cross_{q_{n-3}}(\yy))}$}\med

is defined for all $\yy\in O_\xx$ and equal to $\cross_q(\yy)$ as desired.\med


(d) {\it The case $n=5$.} This case represents already some of the main ideas of the argument, so let us start with it. Choose an extremal string $\xx$ in $\YY_5$ (there are $15$ of them, and all their phylogenetic trees have the same shape with two inner edges - they differ only in the distribution of the five labels on the leaves). To fix one, let $\Gamma_\xx$ have an inner vertex $u$ with leaves $i<j$ attached to it, a second inner vertex $v$ with leaf $k$, and a third inner vertex $w$ with leaves $l<m$ attached to it, see Fig.~170.


\cl{\fivelabelstreeextremal} 

\cl{{\it Figure 170.} The unique extremal tree with five leaves; it has three vertices.}\med
 

Up to a permutation of the entries, the edge quadruples of $\xx$ are $q_1=(ijlk)$ and $q_2=(lmik)$, associated to the inner edges $e_1=(uv)$ and $e_2=(vw)$, and the numeration of the entries is chosen such that $\cross_{q_1}(\xx)=\cross_{q_2}(\xx)=1$. Another choice of the numeration but subject to these two equalities would at most yield cross-ratio maps $\cross_{q_1}$ and $\cross_{q_2}$ which are multiplicative inverses of the original ones. So there is no loss of generality in our choice.\med

Consider now other quadruples $q$ and their respective formal cross-ratio $[q]$. We claim that $[q]$ can be expressed as a rational function in $[q_1]$ and $[q_2]$. After this, it will be shown that the same identity holds for the associated cross-ratio maps $\cross_q$ on the whole open set $O_\xx$.\med

Up to a permutation of the entries and up to the symmetry in $\Gamma_\xx$ given by swapping $u$ with $w$, it suffices to consider the two quadruples \med

\cl{$q=(ij\ell m), q'=(ijmk)$.}\med

By the triple product formula $[ijk\ell][ij\ell m][ijmk]=1$ and the transformation rules for formal cross-ratios under permutation of the entries one has\med

\cl{$1-[q] = (1-[q_1])(1-[q_2])$, }\med

say,\med

(1) \hs 5.1cm{$[q] = [q_1]+[q_2]-[q_1][q_2]$.}\med

Similarly, one gets\med

\cl{$\ds [q_1]=[q'][q] $,}\med

say,\med

(2) \hs 5.1cm{$\ds [q'] = {[q_1]\over [q_1]+[q_2]-[q_1][q_2]}$.}\med

We have shown\med


{\bf Lemma 2.} \label{lemma:crossratiosnfive} (Cross-ratios for $n=5$) {\it The field $\CR_5(\xi)\subset\kk(\xi_1,...,\xi_5)$ of formal cross-ratios in $n=5$ variables is generated by the two edge-cross ratios $[q_{e_1}]$ and $[q_{e_2}]$ of the (unique) extremal tree $T$ with five leaves, denoting by $e_1$ and $e_2$ the two inner edges of $T$.}\qed


We claim that equalities (1) and (2) also hold on whole $O_\xx$ for the respective cross-ratio functions. For the first, this follows from the definition of $O_\xx$ since the values of $\cross_{q_1}$, $\cross_{q_2}$ and $\cross_q$ in $\xx$ are $1$, hence they are never $\infty$ nor $0$ in any $\yy\in O_\xx$ and (1) holds there. But then also the second equality (2) holds on $O_\xx$, since the evaluations of the numerator and denominator on the right hand side never become $0$ or $\infty$ on $O_\xx$. This proves the claim.\med

Finally, to determine the image $V_\xx=\alpha_\xx(O_\xx)$, let $c_1,c_2$ be coordinates on $(\kk^*)^2$. As all cross-ratio maps considered before are never $0$ or $\infty$ on $O_\xx$, it follows that $V_\xx$ is contained in $(\kk^*)^2$. The preceding formulas then show that $V_\xx$ is actually defined in $(\kk^*)^2$ by the inequality\med

\cl{$V_\xx: c_1+c_2-c_1c_2\neq 0$,}\med

since $\alpha_\xx$ can be inverted there. This finishes the discussion of cross-ratio identities in the case of trees with $n=5$ leaves.\med




(e) {\it The bridge of a quadruple.} For the general case with $n$ labels, one has to construct, for the given extremal string $\xx\in\YY_n$, and each quadruple $q$, a specific subtree $\H_{\xx,q}$ of $\Gamma_\xx$, the {\it $H$-tree}, that will be used to express the associated formal cross-ratio $[q]$ in terms of certain edge cross-ratios defined by $\H_{\xx,q}$. \med

For any quadruple $q=(ijk\ell)$, the tree $\Gamma_\xx$ has a minimal subtree $\H_{\xx,q}$ containing the leaves $i,j,k,\ell$: it looks like a letter `$H$' with endpoints $i,j,k,\ell$, and has two inner vertices of degree $3$; call them $v$ and $w$. The path $\beta_{\xx,q}$ connecting $v$ and $w$ is called the {\it bridge} of $\H_{\xx,q}$ or of $q$, see Fig.~171. If $q=q_e$ is an edge quadruple of $\xx$, defined by the inner edge $e$ of $\Gamma_\xx$, the bridge consists of the single edge $e$.\med


\cl{\Htree}

\cl{{\it Figure 171.} The $H$-tree in blue (dotted) with respect to the leaves labelled by $i,j,k,\ell$,}

\cl{with bridge (bold) between the vertices $v$ and $w$.}\med\med


If we remove the bridge together with its endpoints $v$ and $w$ from $\Gamma_\xx$, we get a forest with several connected components. Each of the labels in $i,j,k,\ell$ is sitting in a different component; there may be more than these four components, but the others do not play any role in the sequel. If each of the four labels $i,j,k,\ell$ is minimal among all labels of their respective component (with respect of the chosen ordering of the labels), then we say that $q$ is {\it minimal} for $\Gamma_\xx$. With this definition, we can also say that the edge quadruples are exactly the minimal quadruples having a bridge of length one.\med\gb

In a first step we will now show that the (formal and function) cross-ratios associated to {\it minimal} quadruples $q$ can be expressed as rational functions in the cross-ratios of the edge quadruples $q_e$ of $\xx$. In the next subsection, we will extend this to arbitrary quadruples.\med\gb\gb


{\bf Lemma 3.} \label{lemma:edgecrossratiosbridge} (Edge cross-ratios of bridge) {\it For an extremal string $\xx\in\YY_n$, let $q$ be a minimal quadruple for $\Gamma_\xx$, and let $e_1,\dots,e_b$ be the edges of the bridge $\beta_{\xx,q}$ of its $H$-tree $\H_{\xx,q}$. There exist polynomials $P$ and $Q$ in $b$ variables expressing the formal cross-ratio $[q]$ as a rational function in the formal cross-ratios of $e_1,\dots,e_b$,\med

\cl{$\ds [q]={P([q_{e_1}],..., [q_{e_b}])\over Q([q_{e_1}],..., [q_{e_b}])}$.}\med

If $\cross_q(\xx)=1$, the same formula is valid on whole $O_\xx$ for the cross-ratio maps: the equality\med

\cl{$\ds \cross_q(\yy)={P(\cross_{q_{e_1}}(\yy),..., \cross_{q_{e_b}}(\yy))\over Q(\cross_{q_{e_1}}(\yy),..., \cross_{q_{e_b}}(\yy))}$}\med

holds for all $\yy\in O_\xx$. In particular, the numerator and denominator of the quotient are non-zero.}\med


{\it Remark.} Considering the case $\cross_q(\xx)=1$ is no restriction as it can always be achieved by a permutation of the entries of $q$. The advantage is then that the numerator and denominator of the quaotient will be non-zero (and, obviously, also $\neq \infty$), so the right hand side is defined.\med


{\it Proof.} If $b=1$ the bridge of $q$ consists of the unique edge $e_1=e$, and the statement is trivial as $q=q_e$. It had been proven already for trees with five leaves in the case $n=5$. So assume that $n\geq 6$. 
%
%
%
We proceed by induction on $b$. Decompose the bridge into a path of length $b-1$ and an edge $e_b$. Let $u$ be their common vertex. Assume that $q=(ij\ell m)$, with $i<j$ and $\ell<m$, and let $k$ be the smallest label of the leaves in the (additional third) component that is attached to $u$, see Fig.~172. \med


\cl{\Htreeinduction}

\cl{{\it Figure 172.} The component at the vertex $u$ of $\Gamma_\xx$ contains the leaf with minimal label $k$.}\med


The quadruple $q_1:=(ijk\ell)$ is minimal for $\Gamma_\xx$ since $\ell<m$ (and  also for the tree $\Gamma^1_\xx$ obtained from $\Gamma_\xx$ by contracting the connected component of $\Gamma_\xx\sm\{u\}$ containing $e_b$ to a single edge with one leaf of label $j$). Its bridge in $\Gamma^1_\xx$ consists of the edges $e_1,...,e_{b-1}$. By induction on $b$, there is an equality expressing the formal cross-ratio \med

\cl{$\ds [q_1]={P_1([q_{e_1}],..., [q_{e_{b-1}}])\over Q_1([q_{e_1}],..., [q_{e_{b-1}}])}$}\med

as a rational function in the formal cross-ratios $[q_{e_1}],...,[q_{e_{b-1}}]$ defined by the edges $e_1,\dots,e_{b-1}$ and their edge quadruples $q_{e_1},...,q_{e_{b-1}}$. And, again by induction, this equality also holds for the evaluation of the respective cross-ratio map $\cross_{q_1}$ at strings $\yy\in O_\xx$ (after applying, if necessary, a suitable permutation of the entries of the quadruples such that $\cross_q(\xx)=\cross_{q_1}(\xx)=1$). In particular, the numerator and denominator of the quotient on the right hand side are non-zero.\med

Similarly, the quadruple $q_2:=(ijkm)$ is minimal for $\Gamma_\xx$ (and also for the tree $\Gamma^2_\xx$ obtained from $\Gamma_\xx$ by contracting the connected component of $\Gamma_\xx\sm\{u\}$ containing $e_b$ to a single edge with one leaf of label $m$). Its bridge consists of the unique edge $e_b$. The same argument as before applies to the associated cross-ratios $[q_2]={P_2([q_{e_b}])\over Q_2([q_{e_b}])}$ and $\cross_{q_2}$ and yields the respective rational expressions for them in terms of $[q_{e_b}]$ and $\cross_{q_{e_b}}$.\med

We will now express the formal cross-ratio $[q]$ in terms of $[q_1]$ and $[q_2]$, and, similarly, the evaluation of the cross-ratio $\cross_q$ at strings $\yy\in O_\xx$.
This uses again the triple product formula,\med

\cl{$[ijk\ell][ij\ell m][ijmk]=1$.}\med

It gives 

\cl{$\ds [q]={[q_2]\over [q_1]}$,}\med

and then, by substitution, the expression of $[q]$ as a rational function in $[q_{e_1}],...,[q_{e_{b-1}}]$,\med

\cl{$\ds [q]={P([q_{e_1}],..., [q_{e_b}])\over Q([q_{e_1}],..., [q_{e_b}])}$.}\med

To prove the same identity for the evaluations of the cross-ratio maps on strings $\yy\in O_\xx$, it is convenient to permute the entries of the quadruples such that $\cross_q(\xx)=\cross_{q_1}(\xx)=\cross_{q_2}(\xx)=1$, i.e., to take $q=(i\ell jm)$, $q_1=(i\ell jk)$, and $q_2=(ikjm)$. With this modification, the evaluations of the associated cross-ratio functions $\cross_{q_1}$ and $\cross_{q_2}$ at strings $\yy\in O_\xx$ are non-zero, and yields for all $\yy\in O_\xx$ the equality \med

\cl{$\ds \cross_q(\yy)={P(\cross_{q_{e_1}}(\yy),..., \cross_{q_{e_b}}(\yy))\over Q(\cross_{q_{e_1}}(\yy),..., \cross_{q_{e_b}}(\yy))}$.}\med


(f) {\it The induction argument.} We now treat the case of arbitary quadruples $q$. Let $a\in N$, and $N':=N\setminus\{a\}$, and consider the projection map $\pi_a:\YY_n\map\YY_{n-1}$ forgetting all entries involving the label $a$. Recall the transformation of the phylogenetic trees defined by $\pi_a$:  if $\yy$ is a string in $\YY_n$ with image $\yy'$ in $\YY_{n-1}$, the tree $\Gamma_{\yy'}$ of $\yy'$ is obtained from the tree $\Gamma_\yy$ of $\yy$ by clipping off from $\Gamma_\yy$ the leaf with label $a$ together with its edge connecting it to an inner vertex of $\Gamma_\yy$, say, $v$, and, in case that $v$ had degree $3$ in $\Gamma_\yy$, by contracting the two remaining edges of $v$ to a single edge, thus eliminating the vertex $v$.\med

This shows in particular that extremal strings $\xx\in\YY_n$ are mapped to extremal strings $\xx'$ in $\YY_{n-1}$. Moreover, $\pi_a$ sends the open neighborhood $O_\xx$ of $\xx$ surjectively onto $O_{\xx'}$ because the inequalities defining $O_\xx$ are still valid after projection for $O_{\xx'}$. A quadruple $q=(ijk\ell)$ with $a\not\in\{i,j,k,\ell\}$ satisfies $\cross_q(\xx)=1$ if and only if $\cross_q(\xx')=1$. The label set $N'$ inherits a total ordering from $N$, hence the edge quadruples of $\xx'$ are again defined.  The edge quadruples of $\xx$ that do not contain $a$ are, in general, not edge quadruples of $\xx'$ (namely, this happens whenever $a$ is attached to an inner vertex $v$ to which a second leaf is attached: clipping off $a$ contracts the inner edge of $v$ and the second outer edge of $v$ to a new outer edge attached to the unique inner neighbour $w$ of $v$). Hence we have to distinguish quadruples in $N'$ that are edge quadruples for $\xx$ and from those that are edge quadruples for $\xx'$.\med

We will embed the field $\CR_{n-1}(\xi')$ of formal cross-ratios in $n-1$ variables $\xi_i'$ into the field $\CR_n(\xi)$ via $N=N'\cup\{a\}$ and the identification $\CR_n(\xi)\isom\CR_{n-1}(\xi')(\xi_a)$.\med\gb


{\bf Lemma 4.} \label{lemma:edgecrossratiosprojection} (Edge cross-ratios and projection) {\it Let $\xx\in\YY_n$ be an extremal string and let $\xx':=\pi_a(\xx)\in\YY_{n-1}$ be its image under the projection $\pi_a:\YY_n\map\YY_{n-1}$, as described above. Let  $q_{e_1},...,q_{e_{n-3}}$ be the edge quadruples of $\xx$ in $N$. For every edge quadruple $q$ of $\xx'$ in $N'$, there exist polynomials $P$ and $Q$ in $n-3$ variables such that the formal cross-ratio $[q]$ of $q$ can be expressed as the quotient\med

\cl{$\ds [q]={P([q_{e_1}],..., [q_{e_{n-3}}])\over Q([q_{e_1}],..., [q_{e_{n-3}}])}$.}\med\gb

The same formula is valid on whole $O_{\xx}$ for the cross-ratio maps: the equality\med

\cl{$\ds \cross_q(\pi_a(\yy))={P(\cross_{q_{e_1}}(\yy),..., \cross_{q_{e_{n-3}}}(\yy))\over Q(\cross_{q_{e_1}}(\yy),..., \cross_{q_{e_{n-3}}}(\yy))}$,}\med

holds for all $\yy\in O_{\xx}$. In particular, the numerator and denominator of the quotient are non-zero.}\med




{\it Proof.} Let $q=(ijk\ell)$ be an edge quadruple for $\xx'$ in $N'$, corresponding to some edge $e=(vw)\in\Gamma_{\xx'}$. The vertices $v$ and $w$ of $\Gamma_{\xx'}$ will also be treated as the respective vertices in $\Gamma_\xx$ (since they are not eliminated by the projection). Then the path from $v$ to $w$ in $\Gamma_\xx$ is either again the edge $e$, or it is a path of length $2$ with a midpoint $u$ whose (unique) leaf has label $a$. To establish the identity for the formal cross-ratios, we have to distinguish two cases.\med

(i) Assume first that the path in $\Gamma_\xx$ from $v$ to $w$ is the edge $e=(vw)$. Remove this edge from $\Gamma_\xx$ together with $v$ and $w$: if then the leaf with label $a$ is not minimal in any of the four connected components, it follows that $q$ is also an edge quadruple for $\xx$, again for the edge $e$; in this case, the proof is finished. Otherwise, we may assume that $a$ is minimal, say, in the component of $\ell$, see Fig.~173. Then, as $q=(ijk\ell)$ is the edge quadruple of $\xx'$ associated to the edge $e$ of $\Gamma_{\xx'}$, the quadruple $q_1:=(ijka)$ is an edge quadruple for $\xx$, associated to the edge $e$ of $\Gamma_\xx$ (up to a permutation of the entries we may assume that $\cross_{q_1}(\xx)=1$). Consider further the quadruple $q_2:=(ij\ell a)$ in $N$, again with $\cross_{q_2}(\xx)=1$ up to a permutation of the entries. Its $\H$-tree in $\Gamma_\xx$ has bridge $\beta_{\xx,q_2}$ ending at the inner vertex $u$ as depicted in Fig.~173. Therefore, $q_2$ is minimal in $\Gamma_\xx$. By Lemma 3, the formal cross-ratio $[q_2]$ is a rational function in the formal cross-ratios $[q_{e_1}],..., [q_{e_{n-3}}]$. By the triple product formula we can express $[q]$ as the quotient $[q]={[q_1]\over [q_2]}$ in $[q_1]=[q_e]$ and $[q_2]$. Substituting for $[q_2]$ yields the desired identity for $[q]$ in terms of $[q_{e_1}],..., [q_{e_{n-3}}]$.\med


\cl{\clippingleafinduction} 

\cl{{\it Figure 173.} The tree $\Gamma_\xx$ (left) and its contraction $\Gamma_{\xx'}$ (right) clipping off the leaf $a$.}\med\gb


(ii) If the path in $\Gamma_\xx$ from $v$ to $w$ has length $2$, then $a$ must have been, by the definition of $\pi_a$, the label of the (unique) leaf attached to the midpoint $u$ of this path, see Fig.~174. But then the bridge $\beta_{\xx,q}$ of $q=(ijk\ell)$ in $\Gamma_\xx$ equals the path from $v$ to $w$. It follows that $q$, whose entries are minimal in the four components resulting from the deletion of $e$ and its endpoints in $\Gamma_{\xx'}$, is is also an edge quadruple in $\Gamma_\xx$, since deleting in $\Gamma_\xx$ the path from $v$ to $w$ produces the same four components. Now the statement follows from Lemma 3. This proves the identity for the formal cross-ratios in all cases.


\cl{\clippingleafinductionii}

\cl{{\it Figure 174.} The tree $\Gamma_\xx$ (left) and its contraction $\Gamma_{\xx'}$ (right) clipping off the leaf $a$.}\med\gb


We are left with the proof of the identity for the evaluations of the cross-ratio functions at strings $\yy\in O_\xx$. By Lemma 3, the numerator and the denominator of the rational function expressing $\cross_{q_2}(\yy)$ as before, for $q_2=(ij\ell a)$,  are non-zero for all $\yy\in O_\xx$. Also, $\cross_{q_1}(\yy)$ and $\cross_{q_2}(\yy)$ are non-zero by definition of $O_\xx$ and since, without loss of generality, $\cross_{q_1}(\xx)=\cross_{q_2}(\xx)=1$. Hence the expression of $\cross_q(\yy)$ as a rational function in $\cross_{q_{e_1}}(\yy), ..., \cross_{q_{e_{n-3}}}(\yy)$ holds for all $\yy$ in $O_\xx$. This proves also the second assertion of the lemma.\qed


Combining the preceding arguments and lemmata we can now prove the injectivity of the chart map $\alpha_\xx$.\med


{\bf Proposition.} \label{proposition:injectivitychartmap} (Injectivity of chart map) {\it Let $n\geq 5$, and let $\xx$ be an extremal string in $\YY_n$ with edge quadruples $q_{e_1},...,q_{e_{n-3}}$ associated to the inner edges $e_1,...,e_{n-3}$ of $\Gamma_\xx$. Let $O_\xx\subset\YY_n$ be the open neighborhood of $\xx$ defined earlier.}

(a) {\it The field of formal cross-ratios $\CR_n(\xi)$ is generated by the edge cross-ratios $[q_{e_1}],...,[q_{e_{n-3}}]$,}\med

\cl{$\CR_n(\xi)=\kk([q_{e_1}],...,[q_{e_{n-3}}])$.}\med

(b) {\it Assume that the entries of $q_{e_1},...,q_{e_{n-3}}$ are numerated such that $\cross_{q_{e_i}}(\xx)=1$ holds for all $i=1,...,n-3$. For every quadruple $q$ with $\cross_q(\xx)=1$ there are polynomials $P_q$ and $Q_q$ in $n-3$ variables such that \med

\cl{$\ds \cross_q(\yy)={P_q(\cross_{q_{e_1}}(\yy),..., \cross_{q_{e_{n-3}}}(\yy))\over Q_q(\cross_{q_{e_1}}(\yy),..., \cross_{q_{e_{n-3}}}(\yy))}$}\med

holds for all $\yy\in O_\xx$. In particular, the numerator and the denominator  are non-zero.}

(c) {\it For every extremal string $\xx\in \YY_n$, the chart map $\alpha_\xx:O_\xx\map (\kk^*)^{n-3}$ sending $\yy\in O_\xx$ to the vector of cross-ratios $\cross_{q_{e_1}}(\yy),..., \cross_{q_{e_{n-3}}}(\yy)$ is injective.}

(d) {\it The image $V_\xx$ of the map $\a_\xx:O_\xx\map (K^*)^{n-3}$ equals the complement of the closed subvariety defined by the vanishing of all polynomials $P_q$ and $Q_q$ occurring in the numerators and denominators of the quotients displayed in assertion (b).}\med


{\it Proof.} Let $q=(ijk\ell)$ be a quadruple in $N$. Let $a\in N$ be any index different from $i,j,k,\ell$, and set $N':=N\setminus\{a\}$. Consider the projection $\pi_a:\YY_n\map\YY_{n-1}$, and let $\xx':=\pi_a(\xx)\in \YY_{n-1}$ be the image of $\xx$. By induction on $n$ and since $q$ is also a quadruple in $N'$, there is a cross-ratio identity expressing $[q]$ in the formal edge cross-ratios of $\xx'$. By Lemma 4, any formal edge cross-ratio of $\xx'$ can be expressed as a rational function in the formal edge cross-ratios of $\xx$. Composing these rational functions, we get the desired cross-ratio identity expressing $[q]$ in the formal edge cross-ratios of $\xx$. This proves (a).\med

As for assertion (b), the same reasoning applies, using at each step the respective identity for the evaluations of the cross-ratio functions in strings $\yy$ of $O_\xx$ and taking into account that the involved numerators and denominators never vanish.\med

Assertion (c) is an immediate corollary of (b), recalling that the entries of a string $\yy\in \YY_n$ can be expressed as rational functions in the cross-ratios $\cross_q(\yy)$ of $\yy$. This establishes the assertions of the proposition and also concludes the proof of the smoothness of $\XX_n$.\med

Finally, assertion (d) holds because the inverse map of $\a_\xx$ is defined and maps $(n-3)$-tuples into the open set $O_\xx$ of strings $\yy$ whose cross-ratios are either non-special or equal to the cross-ratio of $\xx$ if and only if all  polynomials $P_q$ and $Q_q$ do not vanish in $(K^*)^{n-3}$. This implies (d) and establishes the proposition.\qed\med\gb



{\Bf 18. The stratification of $\XX_n$} \label{section:stratification}\med

The strata of $\XX_n$ are given by strings $\xx$ with the same phylogenetic tree $\Gamma_\xx =T$. The dense open stratum $\UU_n$ corresponds to the generic tree $T_*$ with one inner vertex and $n$ leaves. The remaining strata fill up the boundary divisor $\BB_n=\XX_n\sm\UU_n$.   \med


{\bf Proposition.} \label{proposition:stratification} (Stratification) {\it For $T$ a phylogenetic tree with $n$ leaves, denote by $\SS_T$ the stratum of strings $\xx\in \XX_n$ with $\Gamma_\xx=T$.}

(a) {\it The strata $\SS_T$ are non-empty locally closed smooth subvarieties of $\XX_n$.}

(b) {\it A string $\yy\in\XX_n$ lies in the Zariski closure of $\SS_T$ if and only if $T$ can be obtained from $\Gamma_\yy$ by edge contractions.}\med

{\it Proof.} By the characterization of the trees of strings by cross-ratios in Lemma 3 of Section 13,  there exists a set $Q_T^1\subset N^{\ul 4}$ of quadruples $q$ such that, for a string $\xx\in\XX_n$, one has $\Gamma_\xx=T$ if and only if $\cross_q(\xx)=1$ for $q\in Q_T^1$ and $\cross_q(\xx)\neq  1$ for $q\not\in Q_T^1$. This shows that the strata are locally closed, and that their closure is given by the equations $\cross_q(\xx)=1$ for all $q\in Q_T^1$. Furthermore, we have already seen in the proposition about incidence graphs in Section 13 that every phylogenetic tree arises as the tree $\Gamma_\xx$ of a string. Hence the strata $\SS_T$ are non-empty. The smoothness of the strata follows from the explicit chart maps in the proof of the smoothness of $\XX_n$ in Section 17. Each stratum is an intersection of components of the boundary $\BB_n$ of $\XX_n$; as the boundary is a divisor with normal crossings, the strata are smooth. This shows (a).\med

For (b), assume first, for a given string $\yy\in\XX_n$, that $T$ is obtained from $\Gamma_\yy$ by contracting some edges in $\Gamma_\yy$. We have to show that $\yy$ belongs to the closure of $\SS_T$. By the proof of (a) it suffices to show that, if $\xx\in \XX_n$ has tree $\Gamma_\xx=T$ and $q$ is a quadruple with $\cross_q(\xx)=1$, then this implies $\cross_q(\yy)=1$ (so that $\yy$ lies in the closure of $\SS_T$). Set $q=(ijk\ell)$. Then $\Gamma_\xx$ has an edge $e$ that separates $i,j$ from $k,\ell$. As $T=\Gamma_\xx$ is obtained from $\Gamma_\yy$ by edge contratctions, also $\Gamma_\yy$ has an edge that separates $i,j$ from $k,\ell$. This shows that $\cross_q(\yy)=1$ as desired.\med

Conversely, let $\yy$ belong to the closure of $\SS_T$. We have to show that $T$ is induced from $\Gamma_\yy$ by edge contractions. Again, we may use the characterization of strata via cross-ratios from the proof of (a). But replacing an inequality $\cross_q(\yy)\neq 1$ by an equality $\cross_q(\yy)= 1$ corresponds precisely to an edge contraction in the tree $\Gamma_\yy$. This proves (b).\qed
\med


{\Bf 19. The boundary divisor $\BB_n=\XX_n\sm \UU_n$} \label{section:boundary}\med



The variety $\XX_n$ was defined as the Zariski-closure of the image of $\UU_n= ((\P^1)^n\sm\Delta_n)/\PGL_2$ in $(\P^1)^{n{n\choose 3}}$ under the symmetrization map $\Sigma_n$. Its {\it boundary} $\BB_n=\XX_n\sm\UU_n$ consists of the strings $\xx$ which have been added as limits of generic strings. \med\gb


{\bf Proposition.} \label{proposition:boundary} (Boundary of $\XX_n$) {\it The boundary $\BB_n= \XX_n\setminus{\cal U}_n$ of $\XX_n$ is a divisor with normal crossings. It is a union of smooth hypersurfaces $\DD_{I,J}$, indexed by pairs $(I,J)$ of complementary subsets $I,J$ of cardinality $\geq 2$ of the label set $N$, $I\sqcup J=N$, and these hypersurfaces intersect transversally and at most pairwise. Each $\DD_{I,J}$ is the Zariski-closure in $\BB_n$ of the set of strings whose tree has exactly two inner vertices $v$ and $w$ with the two destination sets $I_v=I$ and $I_w=J$ of cardinality $\geq 2$ (see Fig.~19). } \med


{\it Proof.} Recall that the open dense subset $\UU_n$ of $\XX_n$ consists of strings $\xx$ whose $n$-gons have pairwise distinct entries. Therefore, the phylogenetic tree $\Gamma_\xx$ of a string $\xx$ in $\UU_n$ is {\it generic}, i.e., it has one inner vertex to which $n$ leaves are attached, as in the left most picture of Fig.~12 in Section 12. Consequently, the boundary $\BB_n$ consists of strings $\xx$ whose $n$-gons have at least two equal entries (we may assume $n\geq 4$, since for $n=3$ one has $\BB_n=\emptyset$ and nothing is to prove). Hence, each $n$-gon of a string $\xx\in\BB_n$ must have at least one incidence set $I$ of cardinality $\geq 2$. This implies that the tree $\Gamma_\xx$ has at least two inner vertices $v$ and $w$. We claim that the set $\BB'_n$ of strings $\xx$ whose tree $\Gamma_\xx$ has exactly two inner vertices is dense in $\BB_n$.
These strings have two orbits $v=[x^s]$ and $w=[x^t]$ corresponding to the two vertices of $\Gamma_\xx$, and the respective $n$-gons $x^s$ and $x^t$ have unique complementary incidence sets $I_v$ and $I_w$ of cardinality $\geq 2$, say, $I_v\sqcup I_w=N$, while all other incidence sets are singletons, see Fig.~19.\med


\cl{\twoverticestree}

\cl{{\it Figure 19.} A tree $\Gamma_\xx$ with two vertices and two complementary destination sets $I_v$ and $I_w$.}\med\gb


Let $N=I\sqcup J$ be an arbitrary partition of $N$, with $\abs I, \abs J\geq 2$, and let $\DD_{I,J}\subset\BB_n$ be defined as in the proposition.%
\tex
{{\baselineskip 12pt \footnote{$\!\!{}^6$}{\smallTimes This notation is also used in {\cite{\[Keel\]}} for the corresponding divisors in $\overline{{\cal M}_{0,n}}$.}}}
{{\footnote{\smallTimes This notation is also used in \cite{\[Keel\]} for the corresponding divisors in $\overline{{\cal M}_{0,n}}$.}}}
%
Further, let $\xx\in\XX_n$ be an extremal string, and consider the inner edges of its tree $\Gamma_\xx$. Pick one such edge $e$. Deleting $e$ from $\Gamma_\xx$ while keeping its endpoints cuts $\Gamma_\xx$ into two connected components. Denote by $I_e$ and $J_e$ the sets of labels of the leaves in each component. Contracting all other inner edges $f$ of $\Gamma_\xx$ produces a tree with exactly two inner vertices, the two endpoints of $e$, and incidence sets $I_e$ and $J_e$ of cardinality $\geq 2$. By the proposition in Section 19 on the closure of strata of $\XX_n$, we get $\xx\in \DD_{I_e,J_e}$. Conversely, taking any partition $N=I\sqcup J$ with $\abs I, \abs J\geq 2$, but distinct to all pairs $(I_e, J_e)$, for $e$ an inner edge of $\Gamma_\xx$, the induced tree $T_{I,J}$ with two inner vertices and complementary destination sets $I$ and $J$ cannot be obtained from $\Gamma_\xx$ by edge contractions. This shows that the subvarieties $\DD_{I,J}$ of $\BB_n$ which contain $\xx$ are precisely the subvarieties $\DD_{I_e,J_e}$ associated to edges $e$ of $\Gamma_\xx$ as above. We are left to prove\med



{\bf Lemma.} \label{lemma:boundarycomponents} (Boundary components) {\it Let $\xx$ be an extremal string in $\XX_n$, with chart map $\alpha_\xx:O_\xx\map (\kk^*)^{n-3}, \, \yy\map (\cross_{q_e}(\yy)),$ as in Section 17 on the smoothness of $\XX_n$. Let  $\DD_{I_e,J_e}$ be the subvarieties of $\BB_n$ associated to the $n-3$ edges $e$ of $\Gamma_\xx$ as before, and denote by $c_e$, for $e$ an edge quadruple, the induced coordinates in $\kk^{n-3}$. Then \med

\cl{$\alpha_\xx(\DD_{I_e,J_e}) =V(1-c_e)\subseteq (\kk^*)^{n-3}$.} \med

In particular, the union $\DD_\xx$ of the varieties $\DD_{I_e,J_e}$ is a normal crossings divisor in $(\P^1)^{n{n\choose 3}}$ defined in the open neighborhood $O_\xx$ of $\xx$ by the equation\med

\cl{$\prod_{e\, edge\, of\, \Gamma_\xx}(1-[q_{e}])=0$,}\med

where $[q_{e}]$ denotes the formal cross-ratio of the edge quadruple $q_e$ associated to the edges $e$ of $\Gamma_\xx$.}\med


\comment{[\ooo Internal comment: This proof is not as rigorous as one could wish. Should be reworked in a second version. \ooo]\med}


{\it Proof.} Fix an edge $e$ of $\Gamma_\xx$. Recall that the cross-ratio $\cross_{q_e}$ associated to the edge quadruple $q_e$ of $e$ satisfies $\cross_{q_e}(\xx)=1$. For any $\yy\in O_\xx\cap \DD_{I_e,J_e}$, we therefore have $\cross_{q_e}(\yy)\ne 0$ or $\infty$. Moreover, in any $n$-gon of $\yy$, at least two of the four entries corresponding to the indices in $q_e$ are equal. So $\cross_{q_e}$ is special. Hence $\cross_{q_e}(\yy)=1$. This shows that $\alpha_\xx$ maps $O_\xx\cap \DD_{I_i,J_i}$ into $V(1-c_e)$. \med

Conversely, let $\yy$ be a string in $O_\xx$ which maps under $\a_\xx$ into $V(1-c_e)$, for some inner edge $e$ of $\Gamma_\xx$. Then $\cross_{q_e}(\yy)=1$, by definition of $\a_\xx$. Denote by $q_e=(ijk\ell)$, with $i,j\in I_e$ and $k,\ell\in J_e$, the edge quadruple of $e$. Let $m\in N$ be a label distinct from $i,j,k,\ell$. Because of the cross-ratio relation $1-[ijk\ell]=(1-[ijkm])(1-[imk\ell])$ (see Section 7) we have either $\cross_{(ijkm)}(\yy)=1$ or $\cross_{(imk\ell)}(\yy)=1$. Assume first that $\cross_{(ijkm)}(\yy)=1$. Then $\cross_{(ijkm)}(\xx)=1$ as well. It follows that $m\in I_e$. In the second case, when $\cross_{(imk\ell)}(\yy)=1$, we will have $m\in J_e$. If $o,p\in I_e$, then we also have $\cross_{(opk\ell)}(\yy)=1$, by the triple product formula. So, in any $n$-gon $y^t$ of $\yy$ with $y^t_i\ne y^t_j$, all entries with labels in $J_e$ are equal. Then there is a contraction of edges of $\Gamma_\yy$ such that all labels in $J_e$ are in a single destination set of the contracted tree. Likewise, we can contract more edges such that also $I_e$ becomes a single destination set in the resulting tree. This then implies, using the proposition on stratifications of the last section, that $\yy$ lies $\DD_{I_e,J_e}$. \qed \med\gb

{\Bf 20.~Constructing a stable curve $\CC$ from a string $\xx$.} \label{section:constructingcurvefromstring}\med



In this section, we wish to describe the fibers $\FF_\xx=\pi_a\inv(\xx)$ of the projection $\pi_a:\YY_{n+1}\map \YY_n$ over a string $\xx\in\YY_n$, see assertion (5) of the main theorem and the proposition below. Recall that we will prove later that $\XX_n=\YY_n$, so $\pi_a$ can also be written $\pi_a:\XX_{n+1}\map\XX_n$.  Our goal ist to prove that the fibers are {\it stable curves} in the sense of Deligne-Mumford-Knudsen, i.e., unions of smooth curves meeting transversally, and that the dual graph of $\FF_\xx$ coincides with the skeleton $\sk(\Gamma_\xx)$ of $\Gamma_\xx$, i.e., the tree obtained form the phylogenetic tree $\Gamma_\xx$ of $\xx$ by deleting all leaves together with their edges connecting them to inner vertices. The proof goes in several stages and uses the repeated interplay between the geometry of the tree and the equations defining the fiber in a suitable ambient space.\med

After having finished the proof, we will construct $n$ disjoint sections $\sigma_p$ of $\pi_a$, thus equipping the fibers of $\pi_a$ with $n$ distinct points. As a result, we can interpret the fibers $\FF_\xx$ as {\it $n$-pointed} stable curves, with augmented dual graph equal to $\Gamma_\xx$.\med

Before stating and proving the main result of this section (proposition label{fibers}), we need a couple of auxiliary results collected in Lemmata 1 to 5.\med


Let $N^{+1}=N\cup\{a\}$ and $N$ be the label sets of $\YY_{n+1}$ and $\YY_n$, respectively. The projection map $\pi_a:\YY_{n+1}\map \YY_n$ is given by forgetting all entries $y^t_i$ of strings $\yy=(y^t)_{t\in {N^{+1}\choose 3}}\in \YY_{n+1}$ which involve the label $a$, namely, either as a label of the triple $t$ or because of $i=a$. More precisely, consider the projection\med

\cl{$\theta=\theta_a:(\P^1)^{(n+1){n+1\choose 3}}\map (\P^1)^{(n+1){n\choose 3}}$}\med

forgetting in each string $\yy\in(\P^1)^{(n+1){n+1\choose 3}}$ the $(n+1)$-gons  $y^t\in (\P^1)^{n+1}$ whose triple $t$ involves the label $a$, as well as the projection map\med

\cl{$\rho=\rho_a :(\P^1)^{(n+1){n\choose 3}}\map (\P^1)^{n{n\choose 3}}$}\med

forgetting in each $(n+1)$-gon $y^t$ of a string $\yy$ of $(\P^1)^{(n+1){n\choose 3}}$ the entry $y^t_a$ with index $a$. Then $\pi_a$ is defined as the restriction to $\YY_{n+1}$ of the composition $\rho_a\circ \theta_a:(\P^1)^{(n+1){n+1\choose 3}}\map  (\P^1)^{n{n\choose 3}}$, \med

\cl{$\pi_a={\rho_a\circ \theta_a}_{\vert \YY_{n+1}}: \YY_{n+1}\subseteq (\P^1)^{(n+1){n+1\choose 3}}\map  \YY_n\subset(\P^1)^{n{n\choose 3}}$.}\med


In formulas, let $\yy\in \YY_{n+1}$ be a string, write $\yy=(y^t)_{t\in {N^{+1}\choose 3}}$ as a vector of $(n+1)$-gons $y^t\in (\P^1)^{n+1}$, for triples $t\in {N^{+1}\choose 3}$. Then \med

\cl{$\pi_a(\yy) = \pi((y^t)_{t\in {N^{+1}\choose 3}})=:\zz=(z^t)_{t\in {N\choose 3}}$,}\med

where, for $t\in{N\choose 3}$, each $n$-gon $z^t$ of the image $\zz$ of $\yy$ is obtained from the $(n+1)$-gon $y^t$ of $\yy$ by deleting the entry $y^t_a$, say, such that $y^t=(z^t, y^t_a)$ for $t\in {N\choose 3}$.\med\gb


{\bf Lemma 1.} \label{lemma:thetamap} (Projection $\theta_a$) {\it Let $\YY_{n+1}\subseteq (\P^1)^{(n+1){n+1\choose 3}}$ have label set $N^{+1}=N\cup\{a\}$. The projection \med

\cl{$\theta_a: (\P^1)^{(n+1){n+1\choose 3}} \map (\P^1)^{(n+1){n\choose 3}}$}\med

\cl{$\yy=(y^t)_{t\in {N^{+1}\choose 3}}\map (y^t)_{t\in {N\choose 3}}$}\med

forgetting the $(n+1)$-gons $y^t$ of $\yy$ whose triples $t$ involve $a$ induces by restriction to $\YY_{n+1}$ an isomorphism of $\YY_{n+1}$ onto its image $\theta_a(\YY_{n+1})$ in $(\P^1)^{(n+1){n\choose 3}}$. 

Moreover, $\theta_a(\YY_{n+1})$ equals the subvariety of vectors $(y^t)_{t\in {N\choose 3}}$ of $(n+1)$-gons $y^t$ indexed by triples $t$ in $N\choose 3$ with equal cross-ratios $\cross_q(x^s)=\cross_q(x^t)$ for all quadruples $q$ in $(N^{+1})^{\ul 4}$. In particular, the image of $\theta_a(\YY_{n+1})$ under the projection $(\P^1)^{(n+1){n\choose 3}}\map (\P^1)^{n{n\choose 3}}$ forgetting the last entry in each $(n+1)$-gon is contained in $\YY_n$.}\med\gb



Said differently, the $(n+1)$-gons $y^t$ of a string $\yy\in\YY_{n+1}$ whose triple $t$ does not involve $a$ determine the $(n+1)$-gons $y^t$ whose triple does involve $a$.\med


{\it Proof.} Let $i,j\in N$ be distinct labels, and consider the $(n+1)$-gon $y^{(ija)}$ of $\yy$. We wish to express, for $k\in N\sm\{i,j\}$, the $k$-th entry $y^{(ija)}_k$ of $y^{(ija)}$ in terms of the $(n+1)$-gons $y^t$ of $\yy$ with triples $t\in {N\choose 3}$ not involving $a$. For $k\in\{i,j,a\}$, the entries $y^{(ija)}_k$ are already prescribed as one of the values in $\{0,1,\infty\}$, so these need not be taken into account. Now, the cross-ratio $\cross_{(ijka)}(y^{(ija)})$ determines the entry $y^{(ija)}_k$ of $y^{(ija)}$, since $y^{(ija)}_i=0$, $y^{(ija)}_j=1$, $y^{(ija)}_a=\infty$ are fixed and pairwise distinct. Moreover, both sides of the equality $\cross_{(ijka)}(y^{(ija)}) = \cross_{(ijka)}(y^{(ijk)})$ are defined. It follows that the entries $y^{(ija)}_k$ are determined by $\cross_{(ijka)}(y^{(ijk)})$, and hence by $y^{(ijk)}$, for all $k\in N\sm\{i,j\}$. This is what had to be shown.\qed


The next result, which will not be used further on, describes how the tree $\Gamma_\yy$ of a string $\yy\in\YY_{n+1}$ has to be transformed by clipping off leaves to get the tree $\Gamma_\xx$ of the image $\xx$ of $\yy$ under $\pi_a$,  see Fig.~123 in Section 12.\med 


{\bf Lemma 2.} \label{lemma:treeprojection} (Tree of projection) {\it Let $N^{+1}=N\cup\{a\}$ and $N$ be the label sets of $\YY_{n+1}$ and $\YY_n$. Let $\yy\in\YY_{n+1}$ be mapped to $\xx\in\YY_n$ under the projection $\pi_a:\YY_{n+1}\map \YY_n$. The  phylogenetic tree $\Gamma_\xx$ of $\xx$ is obtained from $\Gamma_\yy$ by deleting from $\Gamma_\yy$ the leaf with label $a$ together with the edge connecting it to an inner vertex, and by contracting, in case that this inner vertex $v$ of $\Gamma_\yy$ had degree three, the two other edges of $v$ to one edge.}\med


{\it Proof.} As the phylogenetic trees are determined by the incidence partitions of the string, it suffices to compute the incidence sets of all $n$-gons $x^t$ of $\xx$ from the incidence sets of the $(n+1)$-gons $y^t$ of $\yy$ for triples $t\in {N\choose 3}$. This is a case by case check which we omit. See Fig.~123 for the various cases which can occur.\qed


For what follows we will rely on a neat result from combinatorics. Let a finite set $N=\{1,...,n\}$ be partitioned into $p$ disjoint non-empty subsets $N_i$, say, $N=N_1\sqcup\ldots\sqcup N_p$.  Define a {\it subjacent $k$-set} of $N$ with respect to the chosen partition as a subset $B$ of $N$ with $k$ elements, each from a different set $N_i$. Of course, such a set only exists if $k\leq p$; it is unique if and only if $k=p=n$ and hence all $N_i$ are one-element sets.\med

{\bf Lemma 3.} \label{lemma:countinglemma} (Counting lemma) {\it For any finite partitioned set $N=N_1\cup\ldots\cup N_p$ there exists, for any subjacent $k$-set $B$ of $N$, a counting (= linear ordering) $B_1=B<B_2<\ldots$ of all subjacent $k$-sets of $N$ starting with $B$ and such that any two consecutive $k$-sets $B_i$ and $B_{i+1}$ share exactly $k-1$ elements.}\med



{\it Proof.} Let a subjacent $k$-set $B$ of $N$ be given. Let $b$ be any element of $B$, write $B=C\cup \{b\}$ with $b\not\in C$. Without loss of generality, $b\in N_p$. Set $M=N\sm N_p$ be equipped with the partition given by the sets $N_1,...,N_{p-1}$. It is clear that $C$ is a subjacent $(k-1)$-set of $M$. By induction on the cardinality $\abs N$, there exists a counting $C_1=C < C_2<\ldots$ of all subjacent $(k-1)$-sets of $M$ starting with $C$ such that any two consecutive sets share $k-2$ elements. The sets $B_i:=C_i\cup\{b\}$ then count all subjacent $k$-sets of $N$ containing $b$; by construction, any two consecutive sets $B_i$ and $B_{i+1}$ share $k-1$ elements. Let $B_m=C_m\cup\{b\}$ with $b\not\in C_m$ be the last set in this counting. If $B_m=N$, we are done. Otherwise, choose any $b'\in N\sm B_m$, and set $B_{m+1}:=C_m\cup\{b'\}$. This is a subjacent $k$-subset of $N$ which does not contain $b$ and which is hence a subjacent $k$-subset of $N'=N\sm\{b\}$ with respect to the partition $N_1,..., N_{p-1}, N_p\sm\{b\}$ (omit the last set if $N_p=\{b\}$). It shares $k-1$ elements with $B_m$, namely the set $C_m$.  By induction on $\abs N$, there exists a counting $B_{m+1}<B_{m+2}<\ldots$ of all subjacent $k$-sets in $N'$ starting with $B_{m+1}$ such that any two consecutive sets share $k-1$ elements. Then the counting \med

\cl{$B_1=B=C\cup\{b\}<B_2<\ldots <B_m=C_m\cup\{b\}<B_{m+1}=C_m\cup\{b'\}<B_{m+2}<\ldots$}\med

provides the required counting of all subjacent $k$-sets of $N$: Indeed, any two consecutive sets will share $k-1$ elements.\qed


{\bf Corollary.} \label{corollary:countingtriples} (Counting of triples) {\it Let $\Gamma$ be a phylogenetic tree with $n$ leaves. For an inner vertex $v$ of $\Gamma$, let $\triples(v)$ be the set of triples $t=(ijk)$ in $N=\{1,...,n\}$ defining $v$ as their meeting point. Then there exists a total ordering of $\triples(v)$ such that any two consecutive triples share two entries.}\qed



{\bf Lemma 4.} \label{lemma:variationquadruples} (Variation of quadruples) {\it Let $\xx=(x^t)\in\YY_n$, ${t\in {N\choose 3}}$, be a string. For each $t\in{N\choose 3}$, let $y^t_a$ be a variable and consider the $(n+1)$-gon $y^t:=(x^t,y^t_a)$. For any triples $s,t\in {N\choose 3}$ and any quadruples $p=(ijka)$ and $q=(ij\ell a)$ in $(N^{+1})^{\ul 4}\sm N^{\ul 4}$ involving the label $a\in N^{+1}\sm N$, one has, setting $o=(ijk\ell)$ and $c=\cross_o(\xx)\in\kk\cup\{\infty\}$,} \med

\cl{$\cross_p(y^s)= c\cdot\cross_q(y^s)$ \hs.5cm and \hs.5cm $\cross_p(y^t)= c\cdot\cross_q(y^t)$.}\med

{\it Therefore, if $o\neq 0,\infty$, the equation $\cross_p(y^s)=\cross_p(y^t)$ is equivalent to $\cross_q(y^s)=\cross_q(y^t)$.}\med


{\it Remarks.} (a) Permuting suitably the entries of $p$ and $q$ one can always achieve that $o\neq 0,\infty$. Therefore, by the transformation rules for cross-ratios under permutation of the entries of the involved quadruple, the equation $\cross_p(y^s)= \cross_p(y^t)$ is equivalent to $\cross_q(y^s)= \cross_q(y^t)$ for all $p$ and $q$ as in the lemma.

(b) By the Counting Lemma 3, any two quadruples $p$ and $q$ in $(N^{+1})^{\ul 4}\sm N^{\ul 4}$ can be connected by a sequence $p_1=p, p_2,..., p_{m-1}, p_m=q$ of quadruples for which any two consecutive ones have three equal entries. By transitivity, this implies that the equations given by the equality of cross-ratios $\cross_q(y^s)=\cross_q(y^t)$ are independent of the choice of the quadruple $q$ in $(N^{+1})^{\ul 4}\sm N^{\ul 4}$.\med


{\it Proof.} We use the triple product formula $[ijk\ell][ij\ell a][ijak] = 1$, say $[ijka]=[ijk\ell][ij\ell a]$. It gives the equalities\med

\cl{$\cross_p(y^s)= \cross_o(y^s)\cdot\cross_q(y^s)$ \hs.5cm and \hs.5cm $\cross_p(y^t)= \cross_o(y^t)\cdot\cross_q(y^t)$,}\med

where $o$ denotes the quadruple $o=(ijk\ell)$. As it does not involve the label $a$, we get that\med

\cl{$\cross_o(y^s)=\cross_o(x^s)$ \hs.5cm and \hs.5cm $\cross_o(y^t)=\cross_o(x^t)$.}\med

Substitution gives\med

\cl{$\cross_p(y^s)= \cross_o(x^s)\cdot\cross_q(y^s)$ \hs.5cm and \hs.5cm $\cross_p(y^t)= \cross_o(x^t)\cdot\cross_q(y^t)$.}\med

But $\cross_o(x^s)=\cross_o(x^t)=\cross_o(\xx)$ are equal, by definition of $\YY_n$ and since $x^s$ and $x^t$ are $n$-gons of $\xx\in\YY_n$. With $c=\cross_o(\xx)\in\kk\cup\{\infty\}$ gives\med

\cl{$\cross_p(y^s)= c\cdot\cross_q(y^s)$ \hs.5cm and \hs.5cm $\cross_p(y^t)= c\cdot\cross_q(y^t)$.}\med

This proves the lemma.\qed\gb



{\bf Lemma 5.} \label{lemma:birationaltransformations} (Birational transformations) {\it For a given string $\xx\in\YY_n$, let $s$ and $s'$ be two (not necessarily increasingly ordered) triples in ${N\choose 3}$ with equivalent $n$-gons $x^s$ and $x^{s'}$. Let further be given a string $\yy$ in the fiber $\FF_\xx=\pi_a\inv(\xx)$, and write $y^s=(x^s,y^s_a)$ and $y^{s'}=(x^{s'},y^{s'}_a)$ for the two $(n+1)$-gons $y^s$ and $y^{s'}$ of $\yy$ associated to $s$ and $s'$. Then $y^s_a$ and $y^{s'}_a$ depend birationally on each other. More explicitly (and for later use), one has, for $s=(ijk)$ and $m\neq i,k$, the M\"obius transformations\med

\hs 3.5cm $\ds y^{s'}_a=1-y^s_a$, \hs 2cm for ${s'}=(jik)$,\med

\hs 3.5cm $\ds y^{s'}_a={1\over 1-y^s_a}$, \hs 1.95cm for ${s'}=(kij)$,\med

\hs 3.5cm $\ds y^{s'}_a={1\over x^s_m} y^s_a$, 
\hs 2.1cm with $x^s_m\neq 0,\infty $, for ${s'}=(imk)$.\med

The remaining transformations follow by composition from the listed ones.}\med\gb


{\it Proof.} By the Counting Lemma 3 above and by transitivity, it is sufficient to consider triples $s$ and $t$ which share two or all entries. Therefore one only has to consider the three cases listed in the lemma. Notice also that $y^s$ and $y^{s'}$ need not be $\PGL_2$-equivalent, so we have to resort for the proof to the equality of cross-ratios of the $(n+1)$-gons of $\yy$ in $\YY_{n+1}\subset(\P^1)^{(n+1){n+1\choose 3}}$.\med


Assume that $s=(ijk)$ and ${s'}=(jik)$. Choose $q=(ijka)$ and get from 
$\cross_q(x^s,y^s_a)=\cross_q(x^t,y^t_a)$ that \med

\cl{$(x^s_i-x^s_k)(x^s_j-y^s_a)(x^{s'}_i-y^{s'}_a)(x^{s'}_j-x^{s'}_k)=
(x^s_i-y^s_a)(x^s_j-x^s_k)(x^{s'}_i-x^{s'}_k)(x^{s'}_j-y^{s'}_a)$.}\med

Now use that $x^s_i=x^{s'}_j=0$, $x^s_j=x^{s'}_i=1$, $x^s_k=x^{s'}_j=\infty$, and get\med

\cl{$(0-\infty)(1-y^s_a)(1-y^{s'}_a)(0-\infty)=
(0-y^s_a)(1-\infty)(1-\infty)(0-y^{s'}_a)$,}\med

say,

\cl{$\ds y^{s'}_a=1-y^s_a$.}\med


Assume that $s=(ijk)$ and ${s'}=(kij)$. Choose $q=(ijka)$ and get from 
$\cross_q(x^s,y^s_a)=\cross_q(x^{s'},y^{s'}_a)$ that \med

\cl{$(x^s_i-x^s_k)(x^s_j-y^s_a)(x^{s'}_i-y^{s'}_a)(x^{s'}_j-x^{s'}_k)=
(x^s_i-y^s_a)(x^s_j-x^s_k)(x^{s'}_i-x^{s'}_k)(x^{s'}_j-y^{s'}_a)$.}\med

Now use that $x^s_i=x^{s'}_k=0$, $x^s_j=1$, $x^s_k=\infty$, $x^{s'}_i=1$, $x^{s'}_j=\infty$ and get\med

\cl{$(0-\infty)(1-y^s_a)(0-y^{s'}_a)(\infty-1)=
(0-y^s_a)(1-\infty)(0-1)(\infty-y^{s'}_a)$,}\med

say,

\cl{$\ds y^{s'}_a={1\over 1-y^s_a}$.}\med


Assume that $s=(ijk)$ and ${s'}=(imk)$. Choose $q=(imka)$ and get from $\cross_q(x^s,y^s_a)=\cross_q(x^{s'},y^{s'}_a)$ that\med

\cl{$(x^s_i-x^s_k)(x^{s}_m-y^s_a)(x^{s'}_i-y^{s'}_a)(x^{s'}_m-x^{s'}_k)=
(x^s_i-y^s_a)(x^{s}_m-x^s_k)(x^{s'}_i-x^{s'}_k)(x^{s'}_m-y^{s'}_a)$.}\med

Now use that $x^s_i=x^{s'}_i=0$, $x^s_j=x^{s'}_m=1$, $x^s_k=x^{s'}_k=\infty$  and get\med

\cl{$(0-\infty)(x^s_m-y^s_a)(0-y^{s'}_a)(1-\infty)=
(0-y^s_a)(x^s_m-\infty)(0-\infty)(1-y^{s'}_a)$,}\med

say,

\cl{$\ds y^{s'}_a={1\over x^s_m}y^s_a$.}\med

We show that $x^s_\ell\neq 0,1$. By definition, $x^{s'}_i=0$, $x^{s'}_m=1$ and $x^{s'}_k=\infty$ are pairwise different, hence, as $[x^{s'}]=[x^s]$, also 
$x^s_i=0$, $x^s_m=0$, $x^s_k=\infty$ are pairwise different. It follows that $x^s_m\neq 0,\infty$.\qed



After these preparations, we can proceed to the description of the fibers $\FF_\xx$ of the projection map $\pi_a$ as stable curves.  \med


{\bf Proposition.} \label{proposition:fibersprojection} (Fibers of projection) {\it  Let $N^{+1}=N\cup\{a\}$ and $N=\{1,...,n\}$ be the label sets of $\YY_{n+1}$ and $\YY_n$, respectively, and let $\pi_a:\YY_{n+1}\map \YY_n$ be the associated projection map forgetting the entries of strings involving $a$. The fibers $\FF_\xx=\pi_a\inv(\xx)\subseteq\YY_{n+1}$ of strings $\xx\in\YY_n$ under $\pi_a$ are connected unions of smooth irreducible rational curves meeting transversally, with dual graph equal to the skeleton $\sk(\Gamma_\xx)$ of the phylogenetic tree $\Gamma_\xx$ of $\xx$.}\med


{\it Proof.} As the argument requires a series of reduction steps, we will first describe the overall strategy and organization. In most stages of the proof it will be crucial to exploit the information encoded in the geometry of the phylogenetic tree $\Gamma_\xx$.\med

The fiber $\FF_\xx$ lives in the high-dimensional projective variety $(\P^1)^{(n+1){n+1\choose 3}}$ and is defined there by an abundant number of equations, each given by an equality of cross-ratios. The idea is then to reduce the ambient dimension $(n+1){n+1\choose 3}$ in two steps by projecting $\FF_\xx$ isomorphically to smaller dimensional projective ambient varieties until one arrives at a closed subvariety $\HH_\xx$ of $(\P^1)^d$ whose ambient dimension is exactly the number $d=d(\Gamma_\xx)$ of inner vertices of $\Gamma_\xx$. The final step is then to show that $\HH_\xx$ is actually defined in $(\P^1)^d$ by $d-1$ equations forming a complete intersection, and that these equations define a curve with smooth rational components intersecting in the way as predicted by the proposition.\med

Recall at that point that every string $\yy$ in $\FF_\xx\subset(\P^1)^{(n+1){n+1\choose 3}}$ has  $(n+1){n+1\choose 3}$ many entries, and we will have to select precisely $d$ of them to get the image of $\yy$ in $\HH_\xx$. This means that we have to choose for each vertex $v$ of $\Gamma_\xx$ a triple $t\in {N\choose 3}$ such that $v=[x^t]$ and then pick the respective entry $y^t_a$ of $\yy$, where $a$ is the label in $N^{+1}\sm N$. To see that $\HH_\xx$ is a curve, we will need (at least) $d-1$ equations defining it in $(\P^1)^d$. But note that there are exactly $d-1$ inner edges in $\Gamma_\xx$. This suggests to define for each such edge $e=(vw)$ an equation $E^e=E^{vw}$ in $(\P^1)^d$. If $v=[x^s]$ and $w=[x^t]$ are the endpoints of $e$ with selected triples $s$ and $t$, the equation $E^e$ will be given by the equality of cross-ratios $\cross_q(y^s)= \cross_q(y^t)$, where $q$ is an (arbitrarily) chosen quadruple in $(N^{+1})^{\ul 4}\sm N^{\ul 4}$. As $y^s=(x^s,y^s_a)$ and $y^t=(x^t,y^t_a)$ have only the last entry unspecified, this will be equations in our coordinates $y^s_a$ and $y^t_a$ on $(\P^1)^d$.\med

There is a small nuisance in this procedure concerning the choice of the triples $t\in {N\choose 3}$ defining the vertices $v=[x^t]$ of $\Gamma_\xx$: The shape of the equations $E^{vw}$ depends on the chosen triples for $v$ and $w$, according to the transformation rules for the variables $y^t_a$ described in Lemma 5. There is no global choice of triples such that all equations $E^{vw}$ assume a systematic form. In view of this, we will work for each edge $e=(vw)$ with specifically chosen triples $s$ and $t$ for $v$ and $w$, keeping in mind that, eventually, the choice does not matter. \med

The first ambient reduction maps the fiber $\FF_\xx=\pi_a\inv(\xx)\subseteq (\P^1)^{(n+1){n+1\choose 3}}$ isomorphically onto a closed subvariety $\GG_\xx$ of $(\P^1)^{{n\choose 3}}$. The coordinates $y^t_a$ in $(\P^1)^{{n\choose 3}}$ will be indexed by triples $t=(ijk)$ in $N\choose 3$. In the next step we construct an isomorphism of $\GG_\xx$ onto the closed subvariety $\HH_\xx$ of $(\P^1)^d$ described before. The coordinates $y^v_a$ in $(\P^1)^d$ will be indexed by the vertices $v$ of $\Gamma_\xx$. It then suffices to prove the asserted properties for $\HH_\xx$. This will go in two steps: first, we select $d-1$ polynomials in the ideal of $K[x^v,\, v\in V(\Gamma_\xx)]$ defining $\HH_\xx$ in $(\P^1)^d$. They will define a closed subvariety $\wt\HH_\xx$ of $(\P^1)^d$ containing $\HH_\xx$. It will be shown that $\wt\HH_\xx$ is a union of smooth curves with the properties listed in the proposition. In a second step, one shows that actually $\wt\HH_\xx=\HH_\xx$: in fact, it will be shown that the chosen polynomials generate the whole ideal defining $\HH_\xx$ in $(\P^1)^d$. This will prove the first part of the proposition via the isomorphisms $\FF_\xx\isom\GG_\xx\isom\HH_\xx$.  \med\gb


(a) {\it Construction of} $\GG_\xx\subseteq (\P^1)^{n\choose 3}$. We start with the isomorphism $\FF_\xx\isom\GG_\xx$. Each string $\yy\in\pi_a\inv(\xx)$ is of the form $\yy=(y^t)_{t\in{N^{+1}\choose 3}}$ with $(n+1)$-gons $y^t\in (\P^1)^{n+1}$. Denote by ${N_a\choose 3}= {N^{+1}\choose 3}\sm{N\choose 3}$ the set of triples $(ija)$ in $N^{+1}\choose 3$ involving the label $a$. Decompose $\yy$ accordingly  into\med

\cl{$\yy=((y^t)_{t\in{N\choose 3}}, (y^t)_{t\in {N_a\choose 3}})\in (\P^1)^{(n+1){n\choose 3}}\times (\P^1)^{(n+1)({n+1\choose 3}-{n\choose 3})}$.}\med




From Lemma 1 we know that the $(n+1)$-gons $(y^t)_{t\in {N_a\choose 3}}$ are completely determined by the $(n+1)$-gons $(y^t)_{t\in{N\choose 3}}$ of $\yy$. Upon replacing $\FF_\xx=\pi_a\inv(\xx)$ with its isomorphic image $\wt\FF_\xx$ under the projection \med

\cl{$\theta_a: (\P^1)^{(n+1){n+1\choose 3}} \map (\P^1)^{(n+1){n\choose 3}}$}\med

forgetting the $(n+1)$-gons $y^t$ indexed by triples $t$ in $N_a$ we may ignore the $(n+1)$-gons $y^t$ of $\yy\in\pi_a\inv(\xx)$ whose triple $t$ involves $a$. We will thus only be concerned with $(n+1)$-gons $y^t$ of $\yy$ for triples $t\in{N\choose 3}$.\med

Whenever $t\in{N\choose 3}$ and $i\in N$, the entry $y^t_i$ of $y^t$ equals $x^t_i$, since $\pi_a(\yy)=\xx$. So only $y^t_a$ is not determined yet - the projection $\pi_a$ forgets this entry. This entry is, however, subject to fulfill the equations imposed by the membership $\yy\in \YY_{n+1}$, namely, the equality of cross-ratios\med

\cl{$E^{st}_q:  \cross_q(x^s,y^s_a)= \cross_q(x^t,y^t_a)$,}\med

for all triples $s$ and $t$ in $N$ and all quadruples $q$ in $N\cup\{a\}$ (with the standard convention that denominators are cleared in order to have polynomial equations). Here, $x^s$ and $x^t$ are considered as constants, whereas $y^s_a$ and $y^t_a$ are treated as variables. These equations define a closed subvariety $\GG_\xx$ of $(\P^1)^{{n\choose 3}}$ isomorphic to $\FF_\xx$, \med

\cl{$\GG_\xx=\{\yy_a:=(y^t_a)_{t\in{N\choose 3}}\in(\P^1)^{n\choose 3},\,  \cross_q(x^s,y^s_a)= \cross_q(x^t,y^t_a)$ for all $s,t\in {N\choose 3}$ and all $q\in N^{\ul 4}\}$.}\med


We have shown that $\FF_\xx\subseteq (\P^1)^{(n+1){n+1\choose 3}}$ is isomorphic to $\GG_\xx\subseteq (\P^1)^{{n\choose 3}}$.\med\gb



(b) {\it Construction of} $\HH_\xx\subseteq (\P^1)^d$. Next we show that $\GG_\xx \subseteq (\P^1)^{{n\choose 3}}$ is isomorphic to a subvariety $\HH_\xx\subset(\P^1)^d$, with $d$ the number of vertices of $\Gamma_\xx$. The construction depends on the choice of a selected triple $t$ for each vertex $v$ of $\Gamma_\xx$ such that $v=[x^t]$ equals the orbit of the $n$-gon $x^t$ defined by $t$. There is some freedom to do so (for instance, one may take for $t$ the triple with lexicographically smallest entries defining $v$). Denote the chosen triple by $t_v$, and let \med

\cl{$\HH_\xx=\{(y^{t_v}_a)_{v\in V(\Gamma_\xx)}\in(\P^1)^d,\,  \cross_q(x^{t_v},y^{t_v}_a)= \cross_q(x^{t_w},y^{t_w}_a)$ for all $v$ and $w$ and all $q\in N^{\ul 4}\}$.}\med

By definition, $\HH_\xx$ is the projection of $\GG_\xx$ on the components indexed by $t_v$, for $v$ a vertex of $\Gamma_\xx$. The birational correspondences of Lemma 5 show that this is actually an isomorphism, so $\GG_\xx\isom\HH_\xx$.
To simplify the notation, and after fixing a choice of triples defining the vertices, we will use coordinates $y^v_a$  on $(\P^1)^d$, indexed directly by the vertices $v\in V(\Gamma_\xx)$ and replacing $y^{t_v}_a$. Setting $y^v=(x^{t_v},y^{t_v}_a)$ we can then write\med

\cl{$\HH_\xx=\{(y^v_a)_{v\in V(\Gamma_\xx)}\in(\P^1)^d,\,  \cross_q(y^v)= \cross_q(y^w)$ for all $v,w\in V(\Gamma_\xx)$ and all $q\in N^{\ul 4}\}$.}\med



(c) {\it Equations for $\HH_\xx$.} We will now calculate the equations defining $\HH_\xx$ in $(\P^1)^d$ in terms of the chosen variables $y^v_a$, for $v\in V(\Gamma_\xx)$, subject to the choice of triples $t$ associated to the vertices $v$. Let $v=[x^s]$ and $w=[x^t]$ be two distinct vertices of $\Gamma_\xx$ and consider the path $\gamma$ between $v$ and $w$. Denote by $I$ and $J$ the two {\it destination sets} of $\gamma$: $I$ is the set of labels of leaves which can be reached from $v$ and going through $\gamma$, i.e., the leaves which from the perspective of $v$ lie {\it behind} $w$. Symmetrically, $J$ is the set of labels of leaves which can be reached from $w$ and going through $\gamma$, i.e., the leaves which  from the perspective of $w$ lie behind $v$, see Fig.~20.


\cl{\destinationsimplepaths}

\cl{{\it Figure 20.} The destination sets $I$ and $J$ of the path between $v=[x^s]$ and $w=[x^t]$.}\med\med\gb


Pick $i,\ell$ in $I$ and $j,k$ in $J$ such that $s=(ijk)$ defines $v$ and $t=(ji\ell)$ defines $w$, see Fig.~20. It does not matter for the moment which labels $i$ and $j$ in $I$ and $J$ are chosen since we show that other choices yield (up to birational coordinate changes) the same equations. But notice that the choice of $s$ (and, symmetrically, that of $t$) depends on both $v$ and $w$, not just on $v$ alone (respectively, $w$ alone). To memorize this preferred choice of triples, we say that the first entries $i$ and $j$ of $s$ and $t$ {\it lie opposite} to $v$ and $w$, as the leaf with label $i$ can be reached from $v$ only passing through $w$, and, conversely, the leaf with label $j$ can be reached from $w$ only passing through $v$. \med

As mentioned earlier, the choice of such distinguished triples $s$ and $t$ for vertices $v$ and $w$ as in Fig.~20 cannot be performed uniformly on $\Gamma_\xx$: For another vertex $w'\neq w$, the distinguished triple defining $v$ with respect to $w'$ may be different from the triple defining $v$ with respect to $w$. The simplest example for this obstruction is a path of length $2$ from $u$ to $w$ with midpoint $v$: For the edge from $u$ to $v$ the distinguished triple $s$ defining $v$ is not the same as the one for the edge from $v$ to $w$, see Fig.~170 in Section 17: For the edge $e_1=(uv)$ the triples $r=(kij)$ and $s=(ik\ell)$ defining $u$ and $v$ are distinguished, whereas for the edge $e_2=(vw)$ the triples $s=(\ell ik)$ and $t=(i\ell m)$ defining $u$ and $v$ are distinguished. There is no choice of $s$ which works for both $r$ and $s$ (except if allowing permutations of the entries, e.g., $s=(i\ell k)$, $r=(\ell ij)$, respectively $s=(lki)$ and $t=(klm)$.) For others situation where this is no longer feasible, see Figs.~200 and  201. \med

 
\cl{\ribisltree}

\cl{{\it Figure 200.} The tree {\it Ribisl} with three edges, four vertices, and six leaves.}\med\gb


{\bf Lemma 6.} \label{lemma:equationsEst}(Equations $E^{st}_{q_{st}}$) {\it  Let $\yy=(y^t)_{t\in {N\choose 3}}\in \FF_\xx$ be a string over $\xx\in\YY_n$ mapping to $(y^v_a)_{v\in V(\Gamma_\xx)}\in \HH_\xx\subset (\P^1)^d$ under the isomorphisms $\FF_\xx\isom\GG_\xx\isom\HH_\xx$. Let $v=[x^s]$ and $w=[x^t]$ be two vertices of $\Gamma_\xx$ with chosen triples $s=(ijk)$ and $t=(ji\ell)$ as just described, say, $i,\ell\in I$, $j,k\in J$, for the destination sets $I$ and $J$ of the path from $v$ to $w$. Set $q_{st}=(ijka)$ and write $y^s_a$ and $y^t_a$ for $y^v_a$, respectively, $y^w_a$. The equation $\cross_{q_{st}}(y^s)=\cross_{q_{st}}(y^t)$ given by the equality of cross-ratios of the $n$-gons of $\yy$ is of the form}\med

\cl{$E^{st}_{q_{st}}:y^s_ay^t_a=0$.}\med

{\it The symmetric choice $q_{ts} = (ji\ell a)$ yields the same equation.}\med


{\it Remark.} We have seen in Lemma 4 that the choice of the quadruple $q$ does not matter at all.\med


{\it Proof.} Recall first that $y^s=(x^s,y^s_a)$ and $y^t=(x^t,y^t_a)$. Hence $\cross_q(y^s)=\cross_q(y^t)$ is equivalent to $\cross_q(x^s,y^s_a)= \cross_q(x^t,y^t_a)$. We thus get, in more explicit form, the equation\med

\cl{$E^{st}_{q_{st}}:
(x^{s}_i-x^{s}_\ell)(x^{s}_j-y^{s}_a)(x^t_i-y^t_a)(x^t_j-x^t_\ell)=
(x^{s}_i-y^{s}_a)(x^{s}_j-x^{s}_\ell)(x^t_i-x^t_\ell)(x^t_j-y^t_a)$,}\med

with variables $y^s_a, y^t_a$. By our convention on strings in $\TT_n$, we know that $x^s_i=x^t_j=0$, $x^s_j=x^t_i=1$,  $x^s_k=x^t_\ell=\infty$. The values of $x^s_\ell=x^s_i=0$ and $x^t_k=x^t_j=0$ follow from the position of the leaves with labels $\ell\in I$ and $k\in J$ in the tree $\Gamma_\xx$. Substitution gives\med

\cl{$E^{st}_{q_{st}}:(0-\infty)(1-y^s_a)(1-y^t_a)(0-0)=
(0-y^s_a)(1-\infty)(1-0)(0-y^t_a)$.}\med
 
The multiplication rules in $\P^1=K\cup\{\infty\}$ show that one may divide on both sides by $0-\infty =-\infty$ and $1-\infty=-\infty$, resulting in the equation\med

\cl{$E^{st}_{q_{st}}: y^s_ay^t_a=0$.}\med

This proves the lemma.\qed\gb



{\it Remark.} The choice of the triples $s=(ijk)$ and $t=(ji\ell)$ for $v$ and $w$ is rather special. Other choices modify the equation by birational correspondences as described in Lemma 6, see Lemma 7 below for the respective formulas. \med\gb

\ignore

(c) [REDUNDANT SINCE $\cross_q(y^s)=\cross_q(y^{s'})$ was used to eliminate $y^{s'}_a$] {\it Variation of quadruples, same vertex.} We show that the equalities of cross-ratios $\cross_q(y^s)= \cross_q(y^t)$, for fixed triples $s,t\in {N\choose 3}$ defining the {\it same} vertex $v=[x^s]=[x^t]$ of $\Gamma_\xx$ do not depend on the quadruple $q$: different quadruples yield equivalent equations in the sense that each equation $\cross_r(y^s)=\cross_r(y^t)$ for some other quadruple $r$ is a constant non-zero and non-infinite multiple of $\cross_q(y^s)=\cross_q(y^t)$. Therefore, in order to describe $\HH_\xx$, it suffices to select one quadruple $q=q_{st}$ for each pair $s,t$ of triples defining $v$ and to then take the respective equation $E^{st}_{q_{st}}$.\med

By the Counting Lemma 3 and transitivity, it suffices to compare equations  $\cross_q(y^s)=\cross_q(y^t)$ and $\cross_r(y^s)=\cross_r(y^t)$ where the quadruples $q$ and $r$ share, aside from the label $a$, exactly two other entries. 
We will show that in this case the two equations are constant multiples of each other.
\med

Assume that $q=(ijka)$ and $r=(ij\ell a)$. We resort to the triple product formula $[ijk\ell][ij\ell a][ijak] = 1$, say $[ijka]=[ijk\ell][ij\ell a]$, for abstract cross-ratios. It gives here the equalities\med

\cl{$\cross_r(y^s)= \cross_o(x^s)\cdot\cross_q(y^s)$,\hs 1cm $\cross_r(y^t)= \cross_o(x^t)\cdot\cross_q(y^t)$,}\med

where $o$ denotes the quadruple $o=(ijk\ell)$ which does not involve the label $a$. If both cross-ratios $\cross_o(x^s)$ and $\cross_o(x^t)$ are defined then $\cross_o(x^s)=\cross_o(x^t)=\cross_o(\xx)$ is a constant in $\P^1$, i.e., independent of $s$ and $t$, since $\xx\in \YY_n$. By a permutation of the labels $i,j,k,\ell$ in the quadruples $q$ and $r$ we may discard the cases where $\cross_o(\xx)=0$ or $\cross_o(\xx)=\infty$. This gives \med

\cl{$\cross_r(y^s)-\cross_r(y^t)= \cross_o(\xx)\cdot[\cross_q(y^s)-\cross_q(y^t)]$}\med

with $\cross_o(\xx)\in\kk^*$ as desired. \med



\recognize



(d) {\it Variation of triples.} The simple form of the equations $E^{st}_{q_{st}}: y^s_ay^t_a=0$ defining $\HH_\xx$ has been due to very specific choices of the triples $s=(ijk)$ and $t=(ji\ell)$ defining the vertices $v$ and $w$. In this part we will look what happens if other triples are chosen instead. This knowledge will become crucial when looking at more than two vertices at the same time. The resulting formulas will be used at three places below: First, when showing in part (e) that $\HH_\xx$ is already defined by the equations $E^{st}_{q_{st}}$ for {\it adjacent} vertices $v=[x^s]$ and $w=[x^t]$. Second, when determining in part (f) the intersection pattern of the irreducible components of $\HH_\xx$. And, finally, in Section 21, when constructing the $n$ sections $\sigma_p$ of the projection $\pi_a:\YY_{n+1}\map \YY_n$.\med


{\bf Lemma 7.} \label{lemma:variationtriples} (Variation of triples) {\it Let $v$ and $w$ be two distinct vertices of the tree $\Gamma_\xx$ of a string $\xx\in \YY_n$. Let $I$ and $J$ be the two destination sets of the path from $v$ to $w$, and let $i,\ell\in I$ and $j,k\in J$ with triples $s=(ijk)$ and $t=(ji\ell)$ defining $v$ and $w$ be chosen as in Lemma 6. Let $g,h,m,o$ be further labels in $J$, with triples $t'=(ijg)$, $t''=(ihk)$ and $r=(imo)$ defining vertices $w'$, $w''$ and $u$, see Fig.~201. Fix the quadruple $q=q_{st}=(ijka)$. The equations resulting from the equality of cross-ratios are as follows.} \med\gb

\hs 1.5cm $s=(ijk)$, $t=(ji\ell)$: \hs 1.15cm $E^{vw}=E^{st}_{q}: \hs 1.2cm y^s_ay^t_a=0$,\med

\hs 1.5cm $s=(ijk)$, $t'=(ijg)$: \hs 1cm $E^{vw'}=E^{st'}_{q}: \hs 1cm (y^s_a-1)y^{t'}_a=0$, \med

\hs 1.5cm $s=(ijk)$, $t''=(ihk)$: \hs 0.9cm $\ds E^{vw''}=E^{st''}_{q}: \hs .85cm {1\over y^s_a-1}y^{t''}_a=0$,\med

\hs 1.5cm $s=(ijk)$, $r=(imo)$: \hs 0.9cm $\ds E^{vu}=E^{sr}_{q}: \hs 1.25cm (y^s_a-x^s_m)y^r_a=0$, with $x^s_m\neq 0,1,\infty$.\med

 
\cl{\viertelvorachttreelemmaseven}

\cl{{\it Figure 201.} The tree {\it Viertel-vor-Acht} with central vertex $v$ adjacent to four vertices $w,w',w'',u$.}\med


{\it Remark.} The third equation has to be understood as an equation in $\P^1$, i.e., it is equivalent to requiring that either $y^{t_*}_a=0$ or that $y^s_a=\infty$.
And, again, by Lemma 4, the choice of the quadruple does not matter.\med


{\it Proof.} The first equation $E^{vw}:y^s_ay^t_a=0$ was the content of Lemma 6 in part (c). For $E^{vw'}$ we choose first $s'=(jik)$, take $q'=(jika)$, apply Lemma 6 and replace in the resulting equation $E^{s't}_q: y^{s'}_ay^t_a=0$ the variable $y^{s'}_a$ by $y^s_a-1$ as indicated in Lemma 5 to get $E^{vw'}=E^{st}_{q'}:(y^s_a-1)y^t_a=0$. For $E^{vw''}$ we choose first $s''=(kij)$, take $q''=(kija)$, apply Lemma 6 and replace in the resulting equation $E^{s''t''}_{q''}: y^{s''}_ay^{t''}_a=0$ the variable $y^{s''}_a$ by $1-y^{s^*}_a$ as indicated in Lemma 5, and then choose $s^*=(ikj)$ to replace $y^{s^*}_a$ by $\ds {y^s_a\over y^s_a-1}$ in order to get eventually $\ds E^{vw''}_{q''}=E^{st''}_{q''}:{1\over y^s_a-1}y^{t''}_a=0$.\med

We are left with $E^{vu}$. We will apply twice the transformation rules from  Lemma 5. Set $s'=(mik)$ and $s''=(imk)$, both triples defining again $v$, and let $r=(imo)$ define $u$. Take $q=(mika)$ and get from Lemma 6 the equation $E^{s'r}_q: y^{s'}_ay^r_a=0$. The first birational transformation from Lemma 5 is $y^{s'}_a=1-y^{s''}_a $, giving $E^{s''r}_q: (y^{s''}_a-1)y^r_a=0$. The second transformation involves the entry $x^s_m$ of $x^s$. It is $y^{s''}_a={1\over x^s_m}y^s_a$, giving $E^{sr}_q: (y^s_a-x^s_m)y^r_a=0$, with  $x^s_m\neq 0,\infty$. From the position of the leaf with label $m$ in $\Gamma_\xx$ it follows that $x^s_m\neq 1$. So the equations for $\FF_\xx$ are exactly as indicated in the lemma.\qed


(e) {\it Minimal equations for $\HH_\xx$.} In the next step we will show that to define $\HH_\xx$ it suffices to consider merely the equations $E^{st}=E^{vw}$ where $v=[x^s]$ and $w=[x^t]$ are {\it adjacent} vertices of $\Gamma_\xx$. As there are exactly $d-1$ edges in $\Gamma_\xx$ these $d-1$ equations will cut out a curve $\wt\HH_\xx$ in $(\P^1)^d$ containing $\HH_\xx$, provided they form a complete intersection (which will, a posteriori, be the case). It will then remain to show that actually $\wt \HH_\xx=\HH_\xx$, with dual graph the skeleton of $\Gamma_\xx$, see part (f).\med

The argument is again based on a look at the phylogenetic tree: If $u$ and $w$ are not adjacent in the tree $\Gamma_\xx$, consider the path from $u$ to $w$. By induction, it suffices to consider a path of length $2$ with vertices $u,v,w$, see Fig.~202, and to show that the equation $E^{uw}$ is a linear combination of the equations $E^{uv}$ and $E^{vw}$ associated to the adjacent vertices $u,v$, respectively, $v,w$.\med

Deleting the path from $\Gamma_\xx$ but keeping its endpoints $u$ and $w$ produces two connected components containing $u$, respectively, $w$, and at least one component which was attached to $v$ before the deletion. Call the first two the {\it left} and the {\it right} components, and an arbitrarily chosen component at $v$ a {\it middle} component. We may then choose five leaves with labels $i,j,k,\ell,m$ belonging to the left component for $i$ and $m$, to the right component for $j$ and $\ell$, and to the middle component for $k$, see Fig.~202. In this way, the triples $r=(jim)$, $s=(ijk)$ and $t=(ji\ell)$ will define $u$, $v$, and $w$, respectively. 


\cl{\choicetriples}

\cl{{\it Figure 202.} The choice of triples $r,s,t$ for the simple path in $\Gamma_\xx$ going from $u$ to $w$ passing through $v$.}\med\gb


The equations associated to the edges $e=(uv)$ and $f=(vw)$ are computed as follows: Use Lemma 6 for the first equation, and take for the second equation $s'=(jik)$, $t'=(ij\ell)$, $q_{s't'}=(jika)$ with $E^{s't'}_{q_{s't'}}: y^{s'}_a y^{t'}_a=0$ by Lemma 6, in order to get from Lemma 5 via $y^{s'}_a=1-y^s_a$ and $y^{t'}_a=1-y^t_a$ the equation $E^{st}_{q_{st}}: (y^s_a-1)(y^t_a-1)=0$. This gives the two equations\med

\cl{$E^{rs}_{q_{rs}}: y^r_ay^s_a=0$,}\med

\cl{$E^{st}_{q_{st}}: (y^s_a-1)(y^t_a-1)=0$.}\med

The equation for the path $(uw)$ is obtained by taking $t''=(ij\ell)$ and $q_{rt''}=(jima)$ with $E^{rt''}_{q_{rt''}}: y^r_a y^{t''}_a=0$ by Lemma 6, and then using Lemma 5 to know that $y^{t''}_a=1- y^t_a$ because of $t=(ji\ell)$. This gives a third equation\med

\cl{$E^{rt}_{q_{rt}}: y^r_a(y^t_a-1)=0$.}\med

Multiplying the first equation with $y^t_a-1$ and subtracting from it the second equation multiplied with $y^r_a$ gives the third equation\med

\cl{$(y^t_a-1) y^r_ay^s_a-y^r_a(y^s_a-1)(y^t_a-1) =y^r_a(y^t_a-1)=0$.}\med

This is what we wanted to show. \med


Let us recapitulate what we have proven so far: The fiber $\FF_\xx\subseteq \XX_{n+1}\subseteq (\P^1)^{(n+1){n+1\choose 3}}$ is isomorphic to the subvariety $\GG_\xx\subseteq (\P^1)^{n\choose 3}$, which, in turn is isomorphic to $\HH_\xx\subseteq (\P^1)^d$. And we have shown that the ideal defining $\HH_\xx$ in $(\P^1)^d$ is generated by $d-1$ equations $E_{vw}=E_e$, given by the equality of cross-ratios of $n$-gons $x^s$ and $x^t$ associated to the endpoints $v$ and $w$ of the inner edges $e=(vw)$ of $\Gamma_\xx$. And each of these equations is of the form $(y^s_a-\delta)(y^t_a-\varepsilon)=0$ with distint constants $\delta,\varepsilon\in\P^1$. It thus defines a normal crossings hypersurface.  \med


(f) {\it Dual graph of $\HH_\xx$.} Let us now investigate the intersection pattern of the irreducible components of $\HH_\xx$. The goal is to show that the dual graph of $\HH_\xx$ equals the skeleton $\sk(\Gamma_\xx)$ of the tree $\Gamma_\xx$ of $\xx$. To this end, we have to associate an irreducible component $\CC_v$ to each inner vertex $v$ of $\Gamma_\xx$, and to show that two components $\CC_v$ and $\CC_w$ intersect if and only if $v$ and $w$ are adjacent to each other. This will also show that $\HH_\xx$ is connected. From the equations it will become clear that $\HH_\xx$ is a union of projective lines in $(\P^1)^d$ intersecting transversally. Actually, each component $\CC_v$ will be a translate of the $v$-coordinate axis in $(\P^1)^d$.\med



Let us write $(y^v)_{v\in V(\Gamma_\xx)}$ for the coordinates in $(\P^1)^d$ (in some order). The description of the component $\CC_v$ associated to the vertex $v$ results from a careful inspection of the equations defining $\HH_\xx$ in $(\P^1)^d$. to make things precise, we will have to choose for each vertex of $\Gamma_\xx$ a special triple in ${N\choose 3}$ defining it. We distinguish two cases.\med


{\it Non-intersection.} We first prove that components corresponding to non-adjacent inner vertices of $\Gamma_\xx$ do not intersect. Let $u$ and $w$ be two inner vertices of $\Gamma_\xx$, and assume that they are not adjacent. Choose a third inner vertex $v$ on the path joining $u$ with $w$. Then choose triples $r$, $s$, $t$ defining $u$, $v$, $w$ in the following manner: $r=(jim)$, $s=(ijk)$, $t=(ji\ell)$, say\med

\hs 5cm $u =[x^{(jim)}]$,\med

\hs 5cm $v =[x^{(ijk)}]$,\med

\hs 5cm $w =[x^{(ji\ell)}]$.\med

A possible configuration is depicted in Fig.~202. Note here that this choice of triples $s$ and $t$ for $v$ and $w$ fulfills the condition in Lemma 6, whereas $r$ and $t$, respectively $r$ and $t$ do not fulfill this condition (actually, there exists no simultaneous choice of triples making all pairs of triples satisfy this condition). \med

Pick any quadruple $q$ in $(N^{+1})^{\ul 4}\sm N^{\ul 4}$, and recall from Lemma 4 that the choice does not affect the form of the equations of $\HH_\xx$ given by the equality of cross-ratios. To remember the choice of triples we will write now $y^u=y^r$, $y^v=y^s$, and $y^w=y^t$. From Lemmata 5 and 6 we get, as in part (e), the following two equations\med\gb

\hs 5cm $E^{rs}_q: y^r_ay^s_a = 0$,\med

\hs 5cm $E^{st}_q: (y^s_a-1)(y^t_a-1) = 0$,\med

while the third equation,\med

\hs 5cm $E^{rt}_q: y^r_a(y^t_a-1)  = 0$,\med

is a linear combination of the first two and thus redundant. The three components $\CC_u=\CC_r$, $\CC_v=\CC_s$ and $\CC_w=\CC_t$ of $\HH_\xx$ are thus obtained as translates of the respective coordinate axes in $(\P^1)^d$. Their images under the projection from $(\P^1)^d$ to $(\P^1)^3$ are therefore (taking on $(\P^1)^3$ the coordinates $y^r_a$, $y^s_a$, $y^t_a$ in this order)\med

\hs 4.5cm image of $\CC_u$: $\P^1\times\{0\}\times\{0\}$,\med

\hs 4.5cm image of $\CC_v$: $\{1\}\times \P^1\times\{0\}$,\med

\hs 4.5cm image of $\CC_w$: $\{1\}\times\{1\}\times\P^1$.\med

This shows that the intersection $\CC_u\cap\CC_w$ is empty as claimed.
\med 


{\it Intersection.} We now prove that components corresponding to adjacent inner vertices of $\Gamma_\xx$ do intersect. So let $u$ and $v$ be two adjacent inner vertices of $\Gamma_\xx$ with defining triples $s$ and $t$, say $u=[x^r]$, $v=[x^s]$. To see that $\CC_u\cap\CC_v\neq \emptyset$, we cannot just project down from $(\P^1)^d$ or $(\P^1)^{n\choose 3}$ to $(\P^1)^2$ with coordinates $y^r_a$ and $y^s_a$. We have to work inside $(\P^1)^d$, considering the components of our subvariety $\HH_\xx$. \med

Let us choose the triples $r$ and $s$ defining $u$ and $v$ such that $r=(jim)$ and $s=(ijk)$ have two equal entries and satisfy the condition in Lemma 6.  Let $w=[y^t]$ be any other vertex of $\Gamma_\xx$ such that $v$ lies on the path connecting $u$ with $w$ as in Fig.~202 (if $w$ lies on the path connecting $u$ with $v$, the argument will be symmetric). Choose a label $\ell$ such that $t=(ji\ell)$ as in Fig.~202. Take further any quadruple $q\in (N^{+1})^{\ul 4}$. The system of equations for $\HH_\xx$ in $(\P^1)^d$ is then given by the equation\med

\hs 5cm $E^{rs}_q: y^r_ay^s_a = 0$,\med

together with the following equations associated to all vertices $w=[x^t]$ distinct from $u$ and $v$: if $v$ lies on the path from $u$ to $w$, the additional equations are
\med

\hs 5cm $E^{st}_q: (y^s_a-1)(y^t_a-1) = 0$;\med

if $u$ lies on the path from $v$ to $w$, the additional equations are (use Lemma 6 and $r'=(ijm)$ to get the equation $E^{r't}_q: y^{r'}_ay^t_a = 0$, and then Lemma 5 for $y^{r'}_a=1-y^r_a$)\med

\hs 5cm $E^{rt}_q: (y^r_a-1)y^t_a = 0$.\med

The points on the component $\CC_u=\CC_r$ of $\HH_\xx$ satisfy the equation $y^r_a=0$, whereas those on $\CC_v$ satisfy $y^s_a=0$. The point of $(\P^1)^d$ with coordinates $y^r_a=0$ and $y^s_a=0$ for the chosen $r$ and $s$, while $y^t_a=1$ for all vertices $w=[x^t]$ for which $v$ lies on the path from $u$ to $w$, and $y^t_a=0$ for all vertices $w=[x^t]$ for which $v$ lies on the path from $u$ to $w$, belongs to the intersection $\CC_v\cap \CC_w$. This shows that $\CC_v\cap \CC_w$ is non-empty.\med

This proves that the dual graph of $\HH_\xx$ and hence also of the fiber $\FF_\xx$ of the projection $\pi_a$ equals the skeleton of $\Gamma_\xx$. The equations also show that the components of $\FF_\xx$ are smooth rational curves intersecting transversally. This concludes the proof of the proposition describing the fiber $\FF_\xx$.\qed\gb


{\it Example.} \label{example:ribisl} (Ribisl) We will determine the fiber $\FF_\xx$ of the string $\xx\in\YY_6$ Ribisl whose tree $\Gamma_\xx$ has four vertices and six leaves, see Fig.~200. We will choose the following triples for the four vertices: $v=[x^s]$ with $s=(ijk)$, $w=[x^t]$ with $t=(ji\ell)$, $w'=[x^{t'}]$ with $t'=(ijg)$, $w''=[x^{t''}]$ with $t''=(ihk)$. As seen in Lemma 6 the equations for $\FF_\xx$ (or, more accurately, for $\HH_\xx$) are \med

\hs 3cm $E^{vw}_q:$ \hs 1.15cm $y^s_ay^t_a=0$,\med

\hs 3cm $E^{vw'}_q:$ \hs 1cm $(y^s_a-1)y^{t'}_a=0$,\med

\hs 3cm $E^{vw''}_q:$ \hs 1cm $\ds{ 1\over y^s_a-1}y^{t''}_a=0$.\med

Ordering the variables in $(\P^1)^d=(\P^1)^4$ according to $v,w,w',w''$ the components of $\FF_\xx$ have the form\med

\hs 3cm $\CC_v=\P^1\times 0\times 0\times 0$,\med

\hs 3cm $\CC_w=0\times \P^1\times 0\times 0$,\med

\hs 3cm $\CC_{w'}=1\times 0\times \P^1 \times 0$,\med

\hs 3cm $\CC_{w''}=\infty\times 0\times 0\times \P^1$.\med

This gives intersections $\CC_v\cap \CC_w=(0,0,0,0)$, $\CC_v\cap \CC_{w'}=(1,0,0,0)$, $\CC_v\cap \CC_{w''}=(\infty,0,0,0)$, while all other intersections are empty. The pattern corresponds exactly to the skeleton $\sk(\Gamma_\xx)$ of $\Gamma_\xx$.\med\med\gb



{\Bf 21.~The sections $\sigma_p$ of $\pi_a:\XX_{n+1}\map \XX_n$} \label{section:sections}\med

We have seen in Section 20 that the fibers $\FF_\xx=\pi_a\inv(\xx)$ of $\pi_a:\XX_{n+1}\map\XX_n$ are stable curves whose dual graph $\Gamma_{\FF_\xx}$ equals the skeleton $\sk(\Gamma_\xx)$ of the phylogenetic tree $\Gamma_\xx$ of the string $\xx\in\XX_n$. We are left to construct $n$ strings $\yy_1,...,\yy_n$ on each fiber in order to turn them into $n$-pointed stable curves. For these it then has to be shown that the $p$-th string $\yy_p$ lies in the component $\CC_v$ of $\FF_\xx$ corresponding to the vertex $v$ of $\Gamma_\xx$ to which the leaf with label $p$ is attached.\med\gb


{\bf Proposition.} \label{proposition:sections} (Sections) {\it  Let $N^{+1}=N\cup\{a\}$ and $N=\{1,...,n\}$ be the label sets of $\YY_{n+1}$ and $\YY_n$, respectively, and let $\pi_a:\YY_{n+1}\map \YY_n$ be the associated projection map. There are $n$-sections $\sigma_1,...,\sigma_n:\YY_n\map\YY_{n+1}$ of $\pi_a$ with disjoint images whose values $\yy_p=\sigma_p(\xx)\in\FF_\xx$ at $\xx$ turn $\FF_\xx$ into an $n$-pointed stable curve with augmented dual graph equal to the phylogenetic tree $\Gamma_\xx$ of $\xx$.}\med


{\it Proof.}  In the preceding Section 20 it was shown that the fibers $\FF_\xx$ are isomorphic to certain curves $\HH_\xx$ in $(\P^1)^d$, where $d$ is the number of inner vertices of $\Gamma_\xx$. It therefore suffices to construct, for each label $p\in N$, morphisms $\rho_p:\XX_n\map (\P^1)^d$ such that $\rho_p(\xx)\in\HH_\xx$ for all $\xx\in\XX_n$, and to then show that $\rho_p(\xx)$ belongs to the component $\CC_v$ of $\HH_\xx$ corresponding to the vertex $v$ of $\Gamma_\xx$ to which the leaf with label $p$ is attached.\med

Fix $p$ between $1$ and $n$ and set $\zz_p=\rho_p(\xx)\in (\P^1)^d$. We will construct $\zz_p$ by prescribing its entries $\zz_p^v$, for $v$ a vertex of $\Gamma_\xx$. Recall that in the proof of Section 20 we have selected for each vertex $v$ of $\Gamma_\xx$ a distinguished triple $t$ with $[x^t] =v$. We may thus write $\zz_p$ as the vector $\zz_p=(z^v_p)$, where $v=[x^t]$ runs over the set of vertices of $\Gamma_\xx$. We then define $\rho_p(\xx)$ and hence the section $\sigma_p$ through\med

\cl{$z^v_p:=x^t_p$,}\med

for all selected triples $t\in {N\choose 3}$ and all labels $p=1,...,n$. Note that $\yy=\sigma_p(\xx)\in\FF_\xx$ is then of the form $\yy=(y^t)_{t\in {N+1\choose 3}}$ with $y^t=(x^t, x^t_p)\in (\P^1)^{n+1}$. With these definitions, it remains to show that $\rho_p(\xx)$ belongs to the component $\CC_v$ of $\FF_\xx$ corresponding to the vertex $v$ of $\Gamma_\xx$ to which the leaf with label $p$ is attached.\med


For simplicity of the exposition, we show this claim only for the fiber $\FF_\xx$ over the string $\xx\in\XX_9$ called {\it Viertel-vor-Neun}, see Fig.~21 (the label $n$ is chosen only due to the lack of letters). This should give a sufficiently convincing idea of how the general argument works.\med

 
\cl{\viertelvorneuntree}

\cl{{\it Figure 21.} The tree {\it Viertel-vor-Neun} with five vertices, four edges,  and nine leaves.}\med


The chosen triples for the five vertices of Viertel-vor-Neun are: $v=[x^s]$ with $s=(ijk)$, $w=[x^t]$ with $t=(ji\ell)$, $w'=[x^{t'}]$ with $t'=(ijg)$, $w''=[x^{t''}]$ with $t''=(ihk)$, $u=[x^r]$ with $r=(imn)$. \med


The equations $E^{vw}$, $E^{vw'}$, and $E^{vw"}$ are identical with those for the example Ribisl at the end of Section 20, so we are only concerned with $E^{vu}$. Again we resort to Lemma 7, now using the last formula from there. So the equations for $\FF_\xx$ are\med

\hs 3cm $E^{vw}:$ \hs 1.15cm $y^s_ay^t_a=0$,\med

\hs 3cm $E^{vw'}:$ \hs 1cm $(y^s_a-1)y^{t'}_a=0$,\med

\hs 3cm $E^{vw''}:$ \hs 1cm $\ds{ 1\over y^s_a-1}y^{t''}_a=0$,\med

\hs 3cm $E^{vu}:$ \hs 1.15cm $(y^s_a-x^s_m)y^r_a=0$, with $x^s_m\neq 0,1,\infty$.\med

Ordering the variables $y^s_a, y^t_a, y^{t'}_a, y^{t''}_a, y^r_a$ in $(\P^1)^d=(\P^1)^5$ according to $v =[x^s],w=[x^t],w'=[x^{t'}],w''=[x^{t''}],u=[x^r]$ the components of $\FF_\xx$ have the form\med

\hs 3cm $\CC_v=\P^1\times 0\times 0\times 0\times 0$,\med

\hs 3cm $\CC_w=0\times \P^1\times 0\times 0\times 0$,\med

\hs 3cm $\CC_{w'}=1\times 0\times \P^1 \times 0\times 0$,\med

\hs 3cm $\CC_{w''}=\infty\times 0\times 0\times \P^1\times 0$,\med

\hs 3cm $\CC_u=x^s_m\times 0\times 0\times  0\times \P^1$.\med

As $x^s_m\neq 0,1,\infty$, it is readily shown that the intersection pattern is precisely the skeleton $\sk(\Gamma_\xx)$ of $\Gamma_\xx$: The component $\CC_u$ intersects all others, while these latter do not intersect among themselves.\med

Once again, the proof goes by case distinctions: If the leaf with label $p=o$ is attached to the central vertex $v=[x^s]$ of $\Gamma_\xx$, we know that $x^s_p\neq x^s_{p'}$ for all $p'\neq p$ in $N$. In particular, $x^s_p\neq 0,1,\infty, x^s_m$. It follows from inspection of the first entry $y^s_a= x^s_p$ of $\rho_p(\xx)$ that $\sigma_p(\xx)$ only belongs to $\CC_v$. \med

Let us now consider the case where the leaf with label $p$ is attached to one of the vertices $w,w',w'',u$. The values of the entries of $\sigma_p(\xx)$ can be read off from the tree $\Gamma_\xx$ depicted in Fig.~21.\med

(1) If the leaf with label $p=i$ is attached to $w=[x^t]$ with $t=(ji\ell)$, we get $y^t_a=x^t_p=x^t_i=1$ and $\sigma_p(\xx)=(0,1,0,0,0)$ only belongs to $\CC_w$. Similarly, for $p=\ell$, $\sigma_p(\xx)=(0,\infty,0,0,0)\in\CC_w$. \med

(2) If the leaf with label $p=j$ is attached to $w'=[x^{t'}]$ with $t'=(ijg)$, we get $\sigma_p(\xx)=(1,0,1,0,0)\in\CC_{w'}$. Similarly, for $p=g$, $\sigma_p(\xx)=(1,0,\infty,0,0)\in\CC_{w'}$.\med

(3) If the leaf with label $p=k$ is attached to $w''=[x^{t''}]$ with $t''=(ihk)$, we get $\sigma_p(\xx)=(\infty,0,0,\infty,0)\in\CC_{w''}$. Similarly, for $p=h$, $\sigma_p(\xx)=(\infty,0,0,1,0)\in\CC_{w''}$.\med

(4) If the leaf with label $p=m$ is attached to $u=[x^r]$ with $r=(imn)$, we get $\sigma_p(\xx)=(x^s_m,0,0,0,1)\in\CC_u$. Similarly, for $p=n$, $\sigma_p(\xx)=(x^s_m,0,0,0,1)\in\CC_u$. \med

This proves the assertion of the proposition for the string Viertel-vor-Neun.\qed\med\gb

\ignore 

{\it Example 1} \label{example:fiberfivegon} (Five leaves) [\ooo still to be revised or to be omitted \ooo] We illustrate the situation for strings $\xx\in \YY_5$ and their fibers $\FF_\xx$ in $\YY_6$. There are two possible trees, aside from the generic one, see Fig.~203.\med


\cl{\phylogenetictreefiveleaves}

\cl{{\it Figure 203.} The two phylogenetic trees with $5$ leaves and $2$, respectively, $3$ inner vertices.}\med


Assume that $\xx\in \YY_5$ has phylogenetic tree $\Gamma_\xx$ with three inner vertices $u$, $v$, $w$ forming a bamboo, see the right hand side of Fig.~203. Let $i$ and $k$ be the labels of the most left vertex $u$, let $j$ and $\ell$ be the labels of the most right vertex $w$, and let $m$ be the label of the leaf of the middle vertex $v$. The associated sets of triples defining $u$, $v$ and $w$ as their meeting points are (up to permutations of their entries)\med\gb

\hs 3cm Triples of $u$: $\{(ikj), (ikm), (ik\ell)\}$, \med

\hs3cm Triples of $v$: $\{(mk\ell), (mkj), (mij), (mi\ell)  \}$,\med

\hs3cm Triples of $w$: $\{(ji\ell),  (jk\ell), (jm\ell)\}$.\med\gb

In this list, the triples are already listed in a way such that any two consecutive triples share precisely two entries: for instance, the last triple $(ik\ell)$ of $\triples(u)$ is followed by the first triple $(mi\ell)$ of $\triples(v)$, and both share two equal entries. We have ${n\choose 3} - 1 =9$ pairs $(s,t)$ of consecutive triples, and accordingly as many distinguished equations. The totality of equations for $\wt \GG_\xx$ will be ${10\choose 2}= 45$.\med

To get a feeling of what is going on, we will determine the hypersurfaces $H^{st}$ for two pairs of triples $(s,t)$, once for $s$ and $t$ defining the same vertex $u$, once for $s$ and $t$ defining two adjacent vertices $u$ and $v$ as their meeting point. It is easy to see that these two cases are already representative for all other cases.\med

(i) Let us compute first the equations $E^{st}_{q_{st}}$ for both $s=(ikj), t=(ikm)$ defining $u$, where we choose the distinguished quadruple $q_{st}=(ikma)$ as earlier. We assume without loss of generality that $i<k<j$ and $i<k<m$. We get from the computations above, replacing there $\ell$ by $m$, \med

\cl{$E^{st}_{q_{st}}: y^s_a(y^t_a-1)=x^s_m(y^t_a-y^s_a)$.}\med

To know whether the hypersurfaces $H^{st}$ is smooth or not, let us check the value of $x^s_m$. But as $s=(ikj)$ and $t=(ikm)$ have both $u = [x^s] = [x^t]$ as their meeting point, it follows that $x^s_m\neq x^s_i=0$ and $x^s_m\neq x^s_k=1$. Therefore the hypersurface $H^{st}$ is smooth and irreducible.\med

(ii) Let us compute next the equations $E^{st}_{q_{st}}$ for $s=(i\ell k)$ defining $u$ and $t=(i\ell m)$ defining $v$, where we choose the distinguished quadruple $q_{st}=(i\ell ma)$ as earlier. We assume without loss of generality that $i<\ell<k$ and $i<\ell<m$. We get from the above computations, replacing there $\ell$ by $m$, \med

\cl{$E^{st}_{q_{st}}: y^s_a(y^t_a-1)=x^s_m(y^t_a-y^s_a)$.}\med

Let us check again the value of $x^s_m$.  But $s=(ikj)$ has $u=[x^s]$ as its meeting point whereas $t=(ikm)$ has $v=[x^t]$ as its meeting point. It follows that $x^s_m=x^s_\ell=1$. Therefore the hypersurface $H^{st}$ is reducible with two smooth transversal components. Altogether we get precisely the three lines in $(\P^1)^3$ described as the images of $\CC_u, \CC_v, \CC_w$ in case 1 of part (f). The claim of the proposition is confirmed in the example.\med\med\gb

\med\med

\recognize



{\Bf 22. The equality of $\XX_n$ with $\YY_n$} \label{section:equalityxxnyyn}\med

We are now finally at the point to be able to prove that the Zariski-closure $\XX_n$ equals the closed subvariety $\YY_n$ of $(\P^1)^{n{n\choose 3}}$.\med


{\bf Proposition.} \label{proposition:equalityxxnyyn} (Equality $\XX_n=\YY_n$) {\it Let $\XX_n=\ol{\Sigma_n(\UU_n)}$ be the Zariski-closure in $(\P^1)^{n{n\choose 3}}$ of the image of $\UU_n=((\P^1)^n\sm \Delta_n)/\PGL_2$ under the symmetrization sending orbits $[x]$ of $n$-gons $x\in (\P^1)^n$ to strings $\xx\in (\P^1)^{n{n\choose 3}}$, and let $\YY_n\subseteq \TT_n\subseteq(\P^1)^{n{n\choose 3}}$ be the closed subvariety of strings $\xx=(x^t)_{t\in {N\choose 3}}$ with equal cross-ratios. Then} \med

\cl{$\XX_n=\YY_n$.}\med


{\it Proof.} We have already shown that $\YY_n$ contains $\XX_n$ and is smooth of dimension $n-3$. Moreover, it is covered by open charts $O_\xx$, and each such chart contains the  stratum $\UU_n$ of generic strings as an open dense subset. This shows that $\YY_n$ is irreducible. Since $\YY_n$ is also closed in the projective ambient variety $(\P^1)^{n{n\choose 3}}$, it follows that it equals the Zariski closure $\XX_n$ of ${\cal U}_n$. \qed\med\med\gb


{\Bf 23. The isomorphism between $\XX_n$ and $\ol\MM_{0,n}$} \label{section:isomorphismxxnmmzeron}\med

%
When we start with a string $\xx\in\XX_n$, its associated stable curve $C=C_\xx$ is obtained by the constructions in Sections 20 and 21. Via Section 15 we recover from $C$ the original string $\xx$ we started with. Conversely, when we start with an $n$-pointed stable curve $C$, associate to it a string $\xx=\xx_C$ as in Section 15. Via Sections 20 and 21 we recover from $\xx$ an $n$-pointed curve $C_\xx$ isomorphic to $C$. This implies that there is a set-theoretic bijection between the (rational) points of $\XX_n$ and $\ol\MM_{0,n}$. \med

Going beyond set theory, it can also be shown that there is an isomorphism of varieties, or, say, rather schemes, between $\XX_n$ and $\ol\MM_{0,n}$. First, to construct a scheme morphism from $\ol\MM_{0,n}$ to $\XX_n$, we can use the cross-ratio functions $\cross_q$ for the quadruples $q\in N^{\ul 4}$. The cross-ratio functions are scheme morphisms because they are compositions of projection maps $\ol\MM_{0,n} \map \ol\MM_{0,n-1}$, up to a permutation of indices. Then the desired map $\ol\MM_{0,n}\map\XX_n$ can be obtained as a map into cartesian products of projective lines, with cross-ratio functions in each component. \med


Second, to construct a scheme morphism from $\XX_n$ to $\ol\MM_{0,n}$, we can use the well-known fact that $\ol\MM_{0,n}$ is a fine moduli space for isomorphism classes of $n$-pointed stable curves (which we do not prove here). Since $\pi_a:\XX_{n+1}\map\XX_n$ is a proper and flat morphism whose fibers are stable curves, it follows that there is a unique map $\Xi_n:\XX_n\map \ol\MM_{0,n}$ such that $\pi_a$ is the pull-back under $\Xi_n$ of the universal family $\ol\MM_{0,n+1}\map \ol\MM_{0,n}$ given by forgetting the $(n+1)$-st point on the (isomorphism classes of) stable curves in $\ol\MM_{0,n+1}$. One can then even prove that the image by this morphism of any string $\xx\in\XX_n$ is the isomorphism class of the stable curve $C_\xx$ constructed in Section 20 as the fiber $\FF_\xx$ of $\pi_a$. \qed\med\gb




{\Bf 24. R\'esum\'e of the proof of the Main Theorem} \label{section:resumeproof}\med

Let us recall the items of the theorem and indicate the place where we proved them.\med

(1) The compactification $\XX_n$ is a smooth, closed and irreducible subvariety of $(\P^1)^{n{n\choose 3}}$: It is clear that $\XX_n$ is irreducible as the closure of the connected open set $\UU_n$ of orbits of generic $n$-gons, and that $\XX_n$ is contained in the subvariety $\YY_n$ of $(\P^1)^{n{n\choose 3}}$ defined by the equality of cross-ratios. In Section 17 it was proven that $\YY_n$ has an open covering by sets $O_\xx$, for $\xx\in\YY_n$ an extremal string, isomorphic to open subvarieties $\VV_\xx$ of $(\P^1)^{n-3}$. Thus $\YY_n$ is smooth. In Section 22 it was shown that $\XX_n=\YY_n$. This gives the required three properties for $\XX_n$.\med

(2) The stratification of $\XX_n$ by strata $\SS_T$ consisting of strings with the same phylogenetic tree $T$ was described in Section 19.\med

(3) That the boundary $\BB_n=\XX_n\sm\UU_n$ of $\XX_n$ is a normal crossings divisor was proven in Section 19, together with the description of its irreducible components.\med

(4) That the projection morphism $\pi_a:\XX_{n+1}\map\XX_n$ is flat follows from a general flatness criterion, using the fact that its fibers are of constant dimension $1$ as proven in Section 20.\med

(5) The proof that the fibers $\FF_\xx$ of $\pi_a:\XX_{n+1}\map\XX_n$ are $n$-pointed stable curves via the choice of appropriate sections $\sigma_1,...,\sigma_n$ of $\pi_a$, and that their augmented dual graph equals $\Gamma_\xx$, covered Sections 20 and 21.\med\gb

(6) \& (7) A bijection between $\XX_n$ and $\ol{\MM_{0,n}}$ was described explicitly in Section 22. That this is also an isomorphism of algebraic varieties could not be proven in this text, nor that $\pi_a:\XX_{n+1}\map\XX_n$ is a universal family for (isomorphism classes of) $n$-pointed stable curves. \med

It should be emphasized here again that our compactification $\XX_n$ of the space $\UU_n=\MM_{0,n}$ of $\PGL_2$-orbits $[x]$ of generic $n$-gons $x$ in $(\P^1)^n\sm \Delta_n$ (as well as the equivalent compactification ${\frak M}_{0,n}$ in {\cite{\[Br\]}}) does not add isomorphism classes of objects (in the present case of $n$-pointed stable curves) to the boundary but picks instead (sufficiently symmetric) representatives (namely, strings) of the generic orbits in a suitable projective variety and then adds as boundary just the points of the ambient variety which belong to the Zariski-closure of the image of $\UU_n$ under the embedding. No equivalence classes had to be considered so far. Only a posteriori isomorphism classes of $n$-pointed stable curves arise, viz, as fibers of $\pi_a:\XX_{n+1}\map\XX_n$.\qed


This concludes our presentation of an alternative approach to $\MM_{0,n}$ and $\ol\MM_{0,n}$. We claim no originality, nor to have completely revealed the subtleties of the constructions and arguments of Deligne-Mumford-Knudsen. Nevertheless, at that point, we very much hope that the reader was able to enjoy wandering with us through the landscape of strings, cross-ratios and phylogenetic trees. The excursion is finished. But wait! What about the moduli space of $n$ points in the projective plane $\P^2$ with respect to the action of the next projective group $\PGL_3$? And indeed, more complicated and also more fascinating combinatorial geometries show up - {\it kite graphs}!\med\med\gb\gb


{\Bf Notation}\med

\hs 1.5cm $\kk$ ground field,\med

\hs 1.5cm $\P^1=\P^1_\kk=\kk\cup\{\infty\}$, projective line over $\kk$,\med

\hs 1.5cm $(\P^1)^n=(\P^1_K)^n$, $n$-fold cartesian product of projective lines,\med

\hs 1.5cm $\PGL_2=\PGL_2(\kk) = \GL_2(\kk)/\kk^*$, general projective linear group,\med

\hs 1.5cm $z\map {az+b\over cz + d}$, M\"obius transformation on $\P^1$,\med

\hs 1.5cm $\xi=(\xi_1,...,\xi_n)$, vector of variables $\xi_i$,\med

\hs 1.5cm $x=(x_1,...,x_n)\in (\P^1)^n$, $n$-gon with entries $x_i\in\P^1$,\med

\hs 1.5cm $[x]$, $\PGL_2$-orbit of $n$-gon $x$,\med

\hs 1.5cm $\Delta_n=\{x\in(\P^1)^n),\, x_i=x_j$ for at least one pair $i\neq j\}$, big diagonal,\med

\hs 1.5cm $(\P^1)^n\sm\Delta_n$, open set of generic $n$-gons with pairwise distinct entries,\med

\hs 1.5cm $\UU_n=((\P^1)^n\sm \Delta_n)/\PGL_2$, orbit space of generic $n$-gons,\med

\hs 1.5cm $\Sigma_n:\UU_n\map (\P^1)^{n{n\choose 3}}$, $[x]\map \xx=(x^t)_{t\in {N\choose 3}}$, symmetrization map,\med

\hs 1.5cm $\TT_n\subseteq(\P^1)^{n{n\choose 3}}=\{\xx=(x^t)_{t\in{N\choose 3}},\, x^t_i=0, x^t_j=1, x^t_k=\infty$ for $t=(ijk)\in{N\choose 3}\}$,\med

\hs 1.5cm $\XX_n=\ol{\Sigma(\UU_n)}\subseteq \YY_n\subseteq\TT_n$, Zariski-closure,\med

\hs 1.5cm  $\YY_n=\{\xx\in\TT_n,\,\cross_q(x^s)=\cross_q(x^t),\,$ for all $s,t\in {N\choose 3},\, q\in N^{\ul 4}\}\subseteq \TT_n$,\med

\hs 1.5cm  $O_\xx\subset\YY_n$, open neighborhood of extremal string $\xx$,\med

\hs 1.5cm  $\BB_n=\XX_n\sm\UU_n$, boundary divisor,\med

\hs 1.5cm  $\DD_{I,J}\subset \BB_n$, strata of boundary divisor,\med

\hs 1.5cm  $\pi_a: \XX_{n+1} \map \XX_n$, projection map,\med\gb

\hs 1.5cm  $\sigma_p: \XX_n \map \XX_{n+1}$, $i=1,...,n$, disjoint sections of $\pi_a$,\med\gb

\hs 1.5cm  $t=(ijk)\in {N\choose 3}$, $q=(ijk\ell)\in N^{\ul 4}$, sets of triples and quadruples,\med

\hs 1.5cm  $[ijk\ell]= {(\xi_i-\xi_k)(\xi_j-\xi_\ell)\over (\xi_i-\xi_\ell)(\xi_j-\xi_k)}\in \kk(\xi_1,...,\xi_n)$, formal cross-ratio,\med

\hs 1.5cm  $\cross_q:(\P^1)^n\map\P^1$, cross-ratio function wrt $q=(ijk\ell)$,\med

\hs 1.5cm  $\DP_v$, $\IP_x$, destination partition of vertex $v$, incidence partition of $n$-gon $x$,\med

\hs 1.5cm  $\Gamma_\xx$, $\Gamma_C$, phylogenetic tree of string, augmented dual graph of $n$-pointed stable curve,\med

\hs 1.5cm  $V(\Gamma_\xx)$, $E(\Gamma_\xx)$, (inner) vertices and edges of $\Gamma_\xx$,\med

\hs 1.5cm  $\sk(\Gamma_\xx)$, skeleton (= edges and inner vertices) of $\Gamma_\xx$,\med

\hs 1.5cm  $N=\{1,...,n\}$, labels of the leaves of $\Gamma_\xx$, $\xx\in \XX_n$,\med

\hs 1.5cm  $N^{+1}=\{1,...,n,a\}$, labels of the leaves of $\Gamma_\yy$, $\yy\in \XX_{n+1}$,\med

\hs 1.5cm  ${N\choose 3}=\{t=(ijk)\}$, (increasingly listed) triples in $N$,\med

\hs 1.5cm  $N^{\ul 4}=\{q=(ijk\ell)\}$, (ordered) quadruples in $N$.\med

\med\med\gb


\tex  
{
{\Bf Bibliography} \med 
{\baselineskip 12pt {\parindent 1.3cm \smallTimes
%
%
\litem{\it Directly related papers}
\litem{[Br]} Brown, F.: Multiple zeta values and periods of moduli spaces ${\frak M}_{0,n}$. Ann.~\'Ec.~Norm.~Sup.~42 (2009), 371-489. \comment{[In Section 2, a compactification $\ol{\frak M}_{0,n}$ of $\MM_{0,n}$ is constructed via the embedding of $\MM_{0,n}$ into the projective variety $(\P^1)^{n\choose 4}$ by associating to an $n$-gon the vector of all its cross-ratios and taking the Zariski-closure of the image. Using dihedral coordinates and labelled polygons it is proven that $\ol{\frak M}_{0,n}$ is smooth and irreducible, and claimed without proof that $\ol{\frak M}_{0,n}$ is isomorphic to $\ol\MM_{0,n}$.]}
\litem{[DM]} Deligne, P., Mumford, D.: The irreducibility of the space of curves of given genus. Inst. Hautes \'Etudes Sci.~Publ.~Math.~36 (1969), 75-109. \comment{[Define in Section 1 the concept of stable curve, according to {\cite{\[KM]}} this is the first appearance, but see also SGA7-I. They only consider curves of genus $g\geq 2$ and thus can dispense with marked points. They prove in Thm.~5.2 that $\ol\MM_g$ is a smooth algebraic stack.]}
\litem{[FMP]} Fulton, W., MacPherson, R.: The compactification of configuration spaces. Ann.~Math.~139 (1994), 183-225. \comment{[Very inspiring paper discussing the moduli problem for $n$ points in higher dimensional projective spaces.]}
\litem{[GHP]} Gerritzen, L., Herrlich, F., van der Put, M.: Stable $n$-pointed trees of projective lines. Proc.~Koninkl.~Nederl.~Akad.~van Wetenschappen A 91(2)  (1988), 131-163. \comment{[The authors define $\ol\MM_{0,n}$ as the Zariski-closure of the image of $\MM_{0,n}$ under a certain embedding into a large projective space using cross-ratios. The marked points on each stable curve are sent to a vector of cross-ratios. Many explicit calculations, often hard to follow.]}
\litem{[Kapr1]} Kapranov, M.: Chow quotients of Grassmannians I. Gelfand Seminar, 29-110, Adv.~Soviet Math.~16.2. Amer.~Math.~Soc. 1993. arXiv:alg-geom/9210002. \comment{[Section 4 treats $\ol\MM_{0,n}$ as the Chow quotient $G(2,n)//H$ of the Grassmannian $G(2,n)$ of lines in $P^{n-1}$ and represents it as a blowup.]}
\litem {[Kapr2]} Kapranov, M.: Veronese curves and Grothendieck-Knudsen moduli space $\MM_{0,n}$. J.~Algbraic Geom.~2 (1993), 239-262. 
\litem{[Knu1]} Knudsen, F.: The projectivity of the moduli space of stable curves. II. The stacks $M_{g,n}$. Math.~Scand.~52 (1983), 161-199. \comment{[Defines the concept of $n$-pointed stable curves and proves, among other things, that $\ol \MM_{0,n}$ is smooth. Uses quite a bit of machinery from algebraic geometry.]}
\litem{[KT2]} Keel, S., Tevelev, J.: Equations for $\ol\MM_{0,n}$.  Intern.~J.~Math.~20 (2009), 1159-1184. \comment{[From {\cite{\[GML]}}: Keel and Tevelev study the image of the particular embedding of $\ol\MM_{0,n}$ into a product of projective spaces given by the complete linear series of the very ample divisor $\kappa$.]}
\litem{[Sin]} Singh, D.: The moduli space of stable $n$-pointed curves of genus zero. Thesis, University of Sheffield 2004.
\med\gb
%
\litem{\it Further papers on moduli}
\litem{[CT]} Castravet, A., Tevelev, J.: $\ol\MM_{0,n}$ is not a Mori dream space. Duke Math.~J.~164 (2015), 1641-1667.
\litem{[FPZ]} Faber, X., Pardue, K., Zelinksy, D.: Cross-ratios of scheme-valued points. arxiv: 2012.03073v1. \comment{[Defines and studies the concept of cross-ratios and the automorphism of the projective lines over rings, i.e., scheme-theoretically.]}
\litem{[GML]} Gibney, A., Maclagan, D.: Equations for Chow and Hilbert quotients. arXiv:0707.1801v2. \comment{[We study the action of an $(n-1)$-dimensional torus $T^{n-1}$ on the Grassmannian $G(2, n)$. Here we can take $U$ to be the points with nonvanishing Pl\"ucker coordinates, and the quotient $U/T^{n-1}$ is the moduli space $\MM_{0,n}$  of
smooth $n$-pointed genus zero curves. In this case the desired compactification is the celebrated moduli space $\ol \MM_{0,n}$ of stable $n$-pointed genus zero curves. Kapranov {\cite{\[Kap2]}} has shown that $\ol \MM_{0,n}$ is isomorphic to both the Chow and Hilbert quotients of $G(2, n)$ by the $T^{n-1}$-action. We give explicit equations for $\ol \MM_{0,n}$ as a subvariety of a smooth toric variety $X_\Delta$ whose fan is the well-studied space of phylogenetic trees. We show that the equations for $\ol \MM_{0,n}$ in the Cox ring $S$ of $X_\Delta$ are generated by the Pl\"ucker relations homogenized with respect to the grading of $S$.]}
\litem{[GMP]} Gelfand, I.~M., MacPherson, R.: Geometry in Grassmannians and a generalization of the dilogarithm. Adv.~Math.~44 (1982), 279-312.
\litem{[Hack]} Hacking, P.: Compact moduli of plane curves. Duke Math.~J.~124 (2004), 213-257. 
%
%
\litem{[Hass]} Hassett, B.: Moduli spaces of weighted pointed stable curves. Adv.~Math.~173 (2003), 316-352.
\litem{[HK]} Hubbard, J., Koch, S.: An analytic construction of the Deligne-Mumford compactification of the moduli space of curves. arXiv:1301.0062.
\litem{[HKT1]} Hacking, P., Keel, S., Tevelev, J.: Compactification of the moduli space of hyperplane arrangements. J.~Algebraic Geom.~15 (2006), 657-680.*
\litem{[HKT2]} Hacking, P., Keel, S., Tevelev, J.: Stable pair, tropical, and log canonical compact moduli of del Pezzo surfaces. Invent.~Math.~178 (2009), 173-227.
\litem{[HMSV1]} Howard, B., Millson, J., Snowden, A., Vakil, R.:  The equations for the moduli space of $n$ points on the line. Duke Math.~J.~146 (2009), 175-226. 
\litem{[HMSV2]} Howard, B., Millson, J., Snowden, A., Vakil, R.: The ideal of relations for the ring of invariants of $n$ points on the line. J.~Eur.~Math.~Soc.~14 (2012), 1-60.
\litem{[HMSV3]} Howard, B., Millson, J., Snowden, A., Vakil, R.: The ideal of relations for the ring of invariants of $n$ points on the line: integrality results. Comm.~Algebra 40 (2012), 3884-3902. 
\litem{[Hos]} Hoskins, V.: Moduli spaces and geometric invariant theory: Old and new perspectives. arXiv:2302.14499.
\litem{[Keel]} Keel, S.: Intersection theory of moduli space of stable $n$-pointed curves of genus zero. Trans.~Amer. Math.~Soc.~330 (1992), 545-574. \comment{[Describes aspects of Knudsen's proof. Provides alternative construction of $\ol\MM_{0,n}$ as the iterated blowup of smooth varieties along codimension $2$ subvarieties. This allows him to deduce results of the intersection theory of this moduli space.]}
\litem{[KM]} Knudsen, F., Mumford, D.: The projectivity of the moduli space of stable curves I: Preliminaries on ``det'' and ``div''. Math.~Scand.~39 (1976), 19-55.
\litem{[Knu2]} Knudsen, F.: A closer look at the stacks of stable pointed curves. J.~Pure Appl.~Algebra 216 (2012), 2377-2385. \comment{[Relies on the functorial description of moduli spaces using stacks and fibered categories.]}
\litem{[KT1]} Keel, S., Tevelev, J.: Geometry of Chow quotients of Grassmannians. Duke Math.~J.~134 (2006), 259-311.
\litem{[MR]} Monin, L., Rana, J.: Equations for $\MM_{0,n}$. In: Combinatorial Algebraic Geometry (G.~Smith and B.~Sturmfels, eds.), Springer 2017, pp.~113- 132.
\med
%
\litem{\it Introductory papers about moduli spaces}
\litem{[Behr]} Behrend, K.: Introduction to algebraic stacks. In: Moduli Spaces (eds. L.~Brambila-Paz, O.~Garcia-Prada, P.~Newstead, R.~Thomas). London Math.~Soc.~Lecture Note Series vol.~411, 2013. \comment{[The first section illustrates the concept of moduli space and universal family in the example of classifying triangles in the real plane.]}
\litem{[Ben]} Ben-Zvi, D.: Moduli Spaces. Princeton Companion to Mathematics (ed.~T.~Gowers). Princeton University Press 2008. 
\litem{[Cap]} Caporaso, L.: Compactifying moduli spaces. Bulletin Amer.~Math.~Soc.~57 (2020), 455-482.
\litem{[Cav]} Cavalieri, R.: Moduli spaces of pointed rational curves. In: Graduate Student School, Combinatorial Algebraic Geometry program at the Fields Institute 2016. \comment{[Gentle introduction to moduli spaces, good motivations, sometimes vague.]}
\litem{[HM]} Harris, J., Morrison, I.: Moduli of curves. Springer 1998. 
\litem{[Kol]} Koll\'ar, J.: Mumford's influence on the moduli theory of algebraic varieties. arXiv: 1809.10723v1. 
\litem{[KV]} Kock, J., Vainsencher, I.: Kontsevich's Formula for Rational Plane Curves. Book 2002. \comment{[The first two sections treat cross-ratios and stable curves. Constructs $\ol\MM_{0,n}$ iteratively, treating explicitly the case $n=5$.]}
\litem{[Maz]} Mazur, B.: Perturbations, deformations and variations (and ``near-misses'') in geometry. Physics and number theory. Bulletin Amer.~Math.~ Soc.~41 (2004), 307-336. \comment{[Provides a very appealing discussion about the problem of ``classification'' in mathematics.]}
\litem{[MFK]} Mumford, D., Fogarty, J., Kirwan, F.: Geometric Invariant Theory. Springer 1994.
\litem{[MS]} Mumford, D., Suominen, K.: Introduction to the theory of moduli. In: Algebraic Geometry, Oslo, 1970: Proceedings of the Fifth Nordic Summer School in Mathematics (ed.~F.~Oort), pp. 171-2333. Wolters- Noordhoff 1972.
\litem{[Tem]} Temkin, M.: Notes on moduli spaces of curves, stacks and de Jong's theorem. Course lecture notes available at: www.math.huji.ac.il/$\sim$temkin.
\litem{[Vor]} Voronov, A.: Lecture 10: Moduli space of algebraic curves of genus zero. Lecture 13: The Fulton-MacPherson compactification. Manuscript 2001.
%
%
\med
\litem{\it Classical papers after 1960}
\litem{[FMP]} Fulton, W., MacPherson, R.: A compactification of configuration spaces. Ann.~Math.~139 (1994), 183-225. \comment{[Describes with much detail the problem of moduli spaces of points in $n$-dimensional projective space.]}
\litem{[Gro]}  Grothendieck, A.: Techniques de construction en g\'eom\'etrie analytique.~X: Construction de l'espace de Teichm\"uller. S\'eminaire Henri Cartan, 13me ann\'ee: 1960/61.
\litem{[Keel1]} Keel, S.: Intersection theory of moduli space of stable $n$-pointed curves of genus zero. Trans.~Amer. Math.~Soc.~330 (1992), 545-574. \comment{[Provides alternative construction of $\ol\MM_{0,n}$ as the iterated blowup of smooth varieties along codimension $2$ subvarieties. This allows him to deduce results of the intersection theory of this moduli space.]}
\litem{[Mat]} Matsumura, H.: Commutative Ring Theory. Cambridge Univ.~Press 1989. \comment{[Flatness miracle]}
\litem{[Mum]} Mumford, D.: Notes on seminar on moduli problems. Supplementary mimeographed notes printed at the Amer.~Math.~Soc. Woods Hole Summer Institute 1964. \comment{[Not looked at yet, cited in {\cite{\[DM\]}}.]}
\litem{[SGA7-I]} Grothendieck, A.: S\'eminaire de g\'eom\'etrie alg\'ebrique du Bois-Marie 1967-1969. Groupes de Monodromie en g\'eom\'etrie alg\'ebrique. Expos\'e I, r\'edig\'e par P.~Deligne. Lecture Notes  Math.~288, Springer 1972. \comment{[Written about the same time as {\cite{\[DM\]}}, the notion of a family of $n$-pointed stable curves of any genus is introduced in Section 5, p.~17, whereas in {\cite{\[DM\]}} only the case of genus $\geq 2$ and hence $n=0$ is considered. Knudsen and Mumford {\cite{\[KM\]}}, p.~19, attribute the introduction of stable curves to Deligne and Mumford {\cite{\[DM\]}}, requiring that every non-singular rational component intersects at least three others to ensure stability. Apparently, they don't consider yet marked points. According to Kapranov {\cite{\[Kap1\]}}, Section 4.1, it was Knudsen who constructed in {\cite{\[Knu1\]}} the moduli space $\ol\MM_{0,n}$ of $n$-pointed stable curves of genus $0$ and proved that it is a smooth variety.]}
%
%
\med\gb
\litem{\it Classical older texts}
\litem{[Riem]} Riemann, B.: Theorie der Abel'schen Functionen. J.~Reine Angew.~Math.~54 (1857), 115-155. \comment{[According to Koll\'ar [Kol], this is the first paper studying moduli problems for algebraic curves, though not in families. This latter viewpoint was only introduced by Grothendieck [Gro] in the Cartan Seminar 1960/61.]}
\litem{[Schr]} Schr\"oder, E.: Vier combinatorische Probleme. T.~Math.~Physik 15 (1870), 361-367. \comment{[Proposes the counting of number of phylogenetic trees.]}
\litem{[Sem]} Semple, J.C.: The triangle as a geometric variable. Mathematika 1 (1954), 80-88. \comment{[Studies the moduli space of triangles in the plane.]}
\litem{[Sev]} Severi, F.: Vorlesungen \"uber algebraische Geometrie (transl.~by E.~L\"offler). Teubner 1924.
%
\med\gb
\litem{\it References for phylogenetic trees}
\litem{[Comt]} Comtet, L.: Advanced Combinatorics, Reidel, 1974.
\litem{[DH.al]} Dress, A., Huber, K., Koolen, J., Moulton, V., Spillner, A.: Basic Phylogenetic Combinatorics. Cambridge Univ.~Press 2012. \comment{[Extensive study of phylogenetic trees, giving various equivalent characterizations.]}
\litem{[Fels]} Felsenstein, J.: The number of evolutionary trees. Syst.~Zool.~27 (1978), 27-33.
\litem{[Schr]} Schr\"oder, E.: Vier combinatorische Probleme. Z.~Math.~Phys.~15 (1870), 361-376.
\litem{[Stan]} Stanley, R.: Enumerative Combinatorics, vol.~II, 1999. 
\litem{[Sull]} Sullivant, S.: Algebraic Statistics. Amer.~Math.~Soc.~2018.
%
%
\med\gb
\litem {\it Papers related to points in $\P^2$}
\litem{[GP]} Gerritzen, L., Piwek, M.: Degeneration of point configurations in the projective line. Indag.~Math.~2 (1991), 39-56. \comment{[Studies the classification of $n$ points in the projective plane $\P^2$.] }
\litem{[SchT]} Schaffler, L., Tevelev, J.: Compactifications of moduli of points and lines in the projective plane. Intern.~Math.~Res.~Notices 21 (2022), 17000-17078. \comment{[Studies the classification of $n$ points in the projective plane $\P^2$, starting from {\cite{\[GP]}}.]}\med
}}
}{}  

\tex
{}
{
\nocite{Br,Del,DM,FMP,GHP,Kapr1,Kapr2,Knu1,KT2,Sin,CT,FPZ,GML,GMP,Hack,
Hass,HKT1,HKT2,HMSV1,HMSV2,HMSV3,Hos,Keel,KM,Knu2,KT1,MR,Behr,Ben,
Cap,Cav,HM,Kol,KV,Maz,MFK,MS,Tem,Vor,Mat,Mum,Riem,Schr,Sem,Sev,Comt,
DHal,Fels,Stan,Sull,GP,SchT,Gro,SGA7-I}
}


\tex{} 
{{\baselineskip 12pt {
\newcommand{\etalchar}[1]{$^{#1}$}
}} 
%
%
\bibliographystyle{alpha} 
}

\vs.5cm


{

H.H.: University of Vienna, {\it herwig.hauser@univie.ac.at}\med

J.Q.: University of Linz, {\it jqi@risc.uni-linz.ac.at}\med

J.S.: University of Linz,  {\it josef.schicho@risc.jku.at}\med

}


\end{document}